\documentclass[12pt,a4paper,twoside]{amsbook}
\usepackage{super-al}
\makeatletter
\renewcommand{\tocappendix}[3]{%
  \indentlabel{#1\@ifnotempty{#2}{#2}}#3}% Uberem ih (VM)
\makeatother

\begin{document}

\title{The analogs of Riemann and Penrose tensors on supermanifolds}

\author{Elena Poletaeva}

\begin{abstract} The Spencer cohomology of certain $\Zee$-graded Lie
superalgebras are completely computed. This cohomology is interpreted
as analogs of Riemann and Penrose tensors on supermanifolds. The
results make it manifest that there is no simple generalization of
Borel-Weil-Bott's theorem for Lie superalgebras.
\end{abstract}

\address{\begin{center}Department of Mathematics, University of California, Riverside, USA; E-mail: elena@math.ucr.edu\end{center}}

\thanks{I warmly thank  IAS (1989), SFB-170 (1990), and
MPIM, Bonn during my several visits in the past decades for
hospitality.  I thank D.~Leites for raising the problem; to him,
A.~Onishchik, J.-L.~Brylinsky, V.~Serganova and A.~Goncharov I am
thankful for help.}



\maketitle

\tableofcontents

This paper expounds the results that appeared mainly in not very
accessible [P1-P8] listed below and complements [LPS] which should
have been published in 1991. Meanwhile there appeared Grozman's
package SuperLie that confirmed our --- not very easy ---
calculations; for some of Grozman's independently obtained results
(the cases of exceptional algebras), see [LPS].

[P1] Poletaeva, E. I. Spencer cohomology connected with some Lie
superalgebras. In: Onishchik A. (ed.) {\em Problems in group theory
and homological algebra} (Russian), 162--167, Matematika, Yaroslav.
Gos. Univ., Yaroslavl, 1988. MR 1 175 010

[P2] Poletaeva, E. Structure functions on the usual and exotic
symplectic and periplectic supermanifolds. In: Bartocci C. et al.
(eds.) {\em Differential geometric methods in theoretical physics}
(Rapallo, 1990), 390--395, Lecture Notes in Phys., 375, Springer,
Berlin, 1991. MR 92h:58207

[P3] Poletaeva, E. I. Spencer cohomology of Lie superalgebras of
vector fields.In: Onishchik A. (ed.) {\em Problems in group theory
and homological algebra} (Russian), 168--169, Matematika, Yaroslav.
Gos. Univ., Yaroslavl, 1990. MR 93c:17036

[P4] Leites, D.; Poletaeva, E. Analogues of the Riemannian structure
for classical superspaces. Proceedings of the International
Conference on Algebra, Part 1 (Novosibirsk, 1989), 603--612, Contemp.
Math., 131, Part 1, Amer. Math. Soc., Providence, RI, 1992. MR
93k:58008

[P5] Poletaeva, E. Analogues of Riemann tensors for the odd metric on
supermanifolds. Acta Appl. Math. 31 (1993), no. 2, 137--169. MR
94d:58166

[P6]Poletaeva, E. Penrose's tensors on super-Grassmannians. Math.
Scand. 72 (1993), no. 2, 161--190. MR 94m:58008a

[P7] Poletaeva, E. Penrose's tensors. II. Math. Scand. 72 (1993), no.
2, 191--211.MR 94m:58008b

[P8] MR96d:58005 Leites, D.; Poletaeva, E. Supergravities and contact
type structures on supermanifolds. Second International Conference on
Algebra (Barnaul, 1991), 267--274, Contemp. Math., 184, Amer. Math.
Soc., Providence, RI, 1995.

Our results also make it clear why one can not just \lq\lq superize"
metric in order to get Einstein-Hilbert's equations. Their true
superizations --- various SUGRAs --- correspond to $\Zee$-graded Lie
superalgebras of depth $d>1$. Such Lie superalgebras were discussed
by Yu. Manin [M] but the corresponding structure functions were not
calculated yet except in the cases considered by Grozman and Leites in

[GL1] Grozman, P.; Leites, D. From supergravity to ballbearings. {\em
Supersymmetries and quantum symmetries} (Dubna, 1997), 58--67,
Lecture Notes in Phys., 524, Springer, Berlin, 1999
MR 2000j:83090

There is, however, a paper where approach similar to the one
described in what follows is applied to $\Zee$-graded Lie
superalgebras of depth $d=1$ and the results are interpreted as
supergravity since the tensor obtained after deleting all that
depends on odd parameters are exactly the standard Riemannian tensor:

[GL2] Grozman, P.; Leites, D.,  An unconventional supergravity.  {\em
Noncommutative structures in mathematics and physics} (Kiev, 2000), 41--47,
NATO Sci. Ser. II Math. Phys. Chem., 22,  Kluwer Acad. Publ.,
Dordrecht, 2001 MR 1 893 452

%%%%%%%%%%%%%%%%%%%%%%%%%%%%%%%%%%%%%%%%%%%%%%%%%%%%%%%%%%%%%%%%%%%%%%%%%%%%%%%%%%%%%%%%%%%%%%%%%%%%%%%%%%%%%%%%%%%%%%%%%%%%%%%%%%%%%%%%%%%%%%%%%%%%%%%%%%%%%%%%%%%%%%%%%%%%%%%%%%%%%%%
%Introduction
%%%%%%%%%%%%%%%%%%%%%%%%%%%%%%%%%%%%%%%%%%%%%%%%%%%%%%%%%%%%%%%%%%%%%%%%%%%%%%%%%%%%%%%%%%%%%%%%%%%%%%%%%%%%%%%%%%%%%%%%%%%%%%%%%%%%%%%%%%%%%%%%%%%%%%%%%%%%%%%%%%%%%%%%%%%%%%%%%%%%%%%%%%%%%%%%%%%%%%%%%%%

\chapter*{Introduction}

\section*{Structure functions}

The main object of the study of Riemannian geometry is the
properties of the Riemann tensor, which in turn splits into the
Weyl tensor, the traceless Ricci tensor, and the scalar curvature.
All these tensors are obstructions to the possibility of
``flattening'' the manifold on which they are considered. The word
``splits'' above means that at every point of the Riemannian
manifold $M^n$ for $n\ne 4$ the space of values of the Riemann
tensor constitutes an $O(n)$-module which splits into the sum of
three irreducible components (for $n = 4$ there are four of them,
because the Weyl tensor splits additionally in this case) [ALV,
Kob].

More generally, let $G\subset\GL(n)$ be any Lie group, not
necessarily  $O(n)$. A reduction of the principal $\GL(n)$-bundle
on $M$ to the principal $G$-bundle is called a
\emph{$G$-structure} on $M$.

Recall that on a manifold with a $G$-structure there is a
canonical connection. For a Riemannian manifold this is the
Levi-Civita connection. The so-called \emph{structure functions}
(SFs) constitute the complete set of obstructions to integrability
of the canonical connection or, in other words, to the possibility
of local flattening of a manifold with $G$-structure. The Riemann
tensor is an example of  a SF. Among the most known other examples
of SFs are the following ones:

\begin{itemize}
\item a \emph{conformal structure}, $G = O(n)\times \Ree^*$, SFs
are called the \emph{Weyl tensor};

\item Penrose's \emph{twistor} theory, $G = S(U(2)\times
U(2))\times \Cee^*$, SFs-Penrose's tensors --- split into two
components
  called the ``$\alpha$-\emph{forms}'' and ``$\beta $-\emph{forms}'';

\item an almost \emph{complex structure}, $G = GL(n;\Cee)\subset
GL(2n;\Ree)$, SFs are called the \emph{Nijenhuis tensor};

\item an almost \emph{symplectic structure}, $G = \Sp(2n)$, no
accepted name for SFs.
\end{itemize}

\section*{Spencer cohomology groups}

Recall necessary definitions [St, Gu].

The simplest $G$-structure is the \emph{flat} $G$-structure
defined as follows. Let $V$ be $\Kee ^n$ with a fixed frame.
Consider the bundle over $V$ whose fiber over $v\in V$ consists of
all frames obtained from the fixed one under the $G$-action, $V$
being identified with $T_vV$.

Obstructions to identification of the $(k + 1)$-st infinitesimal
neighborhood of a point $m\in M$ on a manifold $M$ with
$G$-structure and that of a point of the flat manifold V with the
above $G$-structure are called \emph{structure functions of order
k}. The identification is performed inductively and is possible
provided the obstructions of lesser orders vanish. At each point
of a manifold $M$ SFs take values in  certain cohomology groups,
called \emph{Spencer cohomology groups}. The corresponding complex
is defined as follows. Let $S^iV$ denote the  $i$-th symmetric
power of a vector space $V$ and $Lie(G)$ denote the Lie algebra of
the Lie group $G$. Set $\fg_{-1} = T_m M, \fg_0 = \fg = Lie(G)$
and for $i>0$ put:
$$
\eqalign {
&\fg_i = \lbrace X\in \Hom(\fg_{-1}, \fg_{i-1}): X(v)(w,\ldots ) =
X(w)(v, \ldots ) \hbox{ for any } v,w\in \fg_{-1}\rbrace \cr &=
(\fg_0\otimes S^i(\fg_{-1})^*)\cap (\fg_{-1}\otimes
S^{i+1}(\fg_{-1})^*).\cr }
$$

Now set $\fg_*(\fg_{-1}, \hbox{ }\fg_0) = \oplus _{i\geq {-1}}
\fg_i$. Suppose that the $\fg_0$-module  $\fg_{-1}$ is faithful .

Then $\fg_*(\fg_{-1},\hbox{ } \fg_0)\subset  \fvect(n) =
\emph{der} \Kee[[x_1,\ldots ,x_n]]$, where $n = \dim\fg_{-1}$. It
can be verified that the Lie algebra structure on  $\fvect(n)$
induces such a structure on $\fg_*(\fg_{-1},\hbox{ } \fg_0)$. The
Lie algebra $\fg_*(\fg_{-1}, \hbox{ }\fg_0)$, usually abbreviated
$\fg_*$, will be called the \emph{Cartan prolongation}  of the
pair
  $(\fg_{-1}, \hbox{ }\fg_0)$.

Let $E^iV$ be the $i$-th exterior power of a vector space $V$. Set
$$
C^{k,s}_{\fg_0} = \fg_{k-s}\otimes E^s(\fg_{-1}^*).
$$
Define the differentials $\partial^{k,s}_{\fg_0} :
C^{k,s}_{\fg_0}\longrightarrow C^{k,s+1}_{\fg_0}$ as follows: for
any $g_1, \ldots ,g_{s+1}\in \fg_{-1}$
$$
(\partial^{k,s}_{\fg_0}f)(g_1,\ldots , g_{s+1}) = \sum_i (-1)^i[f(g_1, \ldots ,
\hat {g_{s+1-i}}, \ldots , g_{s+1}),\hbox{ } g_{s+1-i}] \eqno(1)
$$
As expected, $\partial^{k,s}_{\fg_0}\partial^{k,s+1}_{\fg_0} = 0$.
The cohomology of  bidegree $(k,\hbox{ } s)$ of this complex is
called the $(k, \hbox{ }s)$-th Spencer cohomology group
$H^{k,s}_{\fg_0}$. It turns out that structure functions  of order
$k$ on a manifold $M$ with $G$-structure are sections of certain
vector bundles over $M$ with fiber over a point $m \in M$
isomorphic to $H^{k,2}_{\fg}(T_mM)$, where $\fg = Lie(G)$.

\section*{Generalized conformal structures}

A generalization of the notion of conformal structure is a
$G$-structure of type $X$, where $X$ is a classical space, i.e.,
an irreducible compact Hermitian symmetric space (CHSS). These
$G$-structures were introduced and intensively studied by A.
Goncharov, who calculated the corresponding structure functions
[G1, G2]. In his examples $G$ is the reductive part of the
stabilizer of a point of $X$. The usual conformal structure is the
one that corresponds to  $X = Q_n$, a quadric in the projective
space. The complex grassmannian $X = Gr_2^4$ corresponds to
Penrose's twistors.

Recall that  Penrose's idea is to embed the Minkowski space $M^4$
into the complex Grassmann manifold $Gr_2^4$ of planes in $\Cee^4$
(or straight lines in $\Cee\Pee^3$) and to express the conformal
structure on $M^4$ in terms of the incidence relation of the
straight lines in $\Cee\Pee^3$ [Pe].

The conformal structure on $M^4$ is given by a field of quadratic
cones in the tangent spaces to the points of $M^4$. In  Penrose's
case these cones possess two families of two-dimensional flat
generators, the so-called ``$\alpha$-planes'' and
``$\beta$-planes.'' The geometry of these families is vital for
Penrose's considerations. In particular, the Weyl tensor gets a
lucid description in terms of these families.

It is interesting to include  4-dimensional Penrose theory into a
more general theory of geometric structures. A. Goncharov has
shown that there is an analogous field of quadratic cones for any
irreducible compact Hermitian symmetric space $X$ of rank greater
then one [G2].

Let $S$ be a simple complex Lie group, $P$ its parabolic subgroup
with the Levi decomposition $P = GN$, i.e., $G$ is reductive and
$N$ is the radical of $P$. As one knows [He], $N$ is Abelian if
and only if $X = S/P$ is a CHSS, and in this case
  $G = G_0 \times \Cee^*$, where $G_0$ is semisimple.

Let $P_x = G_xN_x$ be the Levi decomposition of the stabilizer of
$x \in X$ in $S$. Denote by $C_x$ the cone of highest weight
vectors in the $G_x$-module $T_xX$, i.e., each element in $C_x$ is
highest with respect to some Borel subgroup in $G_x$. Since $s \in
S$ transforms $C_x$ to $C_{sx}$, then with $X$ there is associated
the cone $C(X) \subset T_{\bar e}X$, where $\bar e$ is the image
of the unit $e \in S$ in $X$.

Let $\rk(X) > 1$, i.e., $X \not = \Cee\Pee^n$. Then on a manifold
$M$ a \emph{generalized conformal structure} of type $X$ is given
if $M$ is endowed with a family of cones $C_m$ and $\Cee$-linear
isomorphisms $A_m: T_{\bar e}X \longrightarrow T_mM$ such that
$A_m(C(x)) = C_m.$

Goncharov has shown that a manifold $M$ with generalized conformal
structure of type $X$ is a manifold with a $\tilde G$-structure,
where $\tilde G$ is a group of linear automorphisms of the cone
$C(X)$ and the connected component of the identity of this group
is precisely $G$ [G2].

\section*{The case of a simple Lie algebra ${\bf \fg_*}$ over $\Cee$}

The following remarkable fact, though known to experts, is seldom
formulated explicitly [LRC, KN].

\begin{Proposition} Let $\Kee = \Cee,\hbox{ } \fg_* = \fg_*(\fg_{-1},
\hbox{ } \fg_0)$ be simple. Then only the following cases are
possible:

\emph{1)} $\fg_2 \not= 0$, then $\fg_*$ is either $\fvect(n)$ or
its special subalgebra $\fsvect(n)$ of divergence-free vector
fields, or its subalgebra $\fh(2n)$ of Hamiltonian vector fields.

\emph{2)} $\fg_2 = 0,\hbox{ } \fg_1 \not=0$, then $\fg_*$ is the
Lie algebra of the complex Lie group of automorphisms of a CHSS
(see \S3).

\end{Proposition}

  Let $R(\sum_i a_i\pi_i)$ be the irreducible $\fg_0$-module with the
highest weight
$\sum_i a_i\pi_i$, where $\pi_i$ is the $i$-th fundamental weight.

\begin{Theorem}[Serre {[St]}] In case 1) of Proposition SFs can only
be of order
  1. More precisely:
for $\fg_* = \fvect(n)$ and $\fsvect(n)$ SFs vanish, for $\fg_* =
\fh(2n)$ nonzero SFs are $R(\pi_1)$ for n = 2, and $R(\pi_1)
\oplus R(\pi_3)$ for $n > 2$.
\end{Theorem}

When $\fg_*$ is a simple finite dimensional Lie algebra over
$\Cee$ computation of SFs becomes an easy corollary of the
Borel-Weil-Bott (BWB) theorem in a form due to W. Shmid [Sh], cf.
work of A. Goncharov [G2]. Indeed, by definition,
$$
\oplus_k H^{k,2}_{\fg_0} = H^2(\fg_{-1},\hbox{ } \fg_*).
$$
The BWB theorem implies that as a $\fg_0$-module, $H^2(\fg_{-1},
\hbox{ }\fg_*)$ has as many components as $H^2(\fg_{-1})$. Thanks
to commutativity of $\fg_{-1}$ one has $H^2(\fg_{-1}) =
E^2\fg_{-1}^*$, which facilitates the count of components. The BWB
theorem also gives the formula for the highest weights of these
components.

\section*{Reduced structures}

Let $X = S/P$, where $Lie(S) =\fg_* = \fg_{-1} \oplus \fg_0 \oplus
\fg_1$ and  $Lie(P) = \fg_0\oplus \fg_1$, be a CHSS. Let $\hat
{\fg_0}$  be the semisimple part of $\fg_0 = Lie(G)$. A $\hat
G$-structure, where $Lie(\hat G) = \hat {\fg_0}$,
  will be referred to  as a
\emph{Riemannian structure } of type X. To reduce the structure
group $G$ to its semisimple part $\hat G$ is an action similar to
distinguishing a metric from a conformal class on a conformal
manifold.

The structure functions of the $\hat G$-structures form an
analogue of the Riemann tensor for the metric. They include the
structure functions of the $G$-structure and several other
irreducible components, some of which are analogues of the
traceless Ricci tensor or the scalar curvature.

More precisely, the structure functions of the $G$-structure are
defined as the part of the structure functions of the $\hat
G$-structure obtained by a reduction of the $G$-structure that
does not depend on the choice of reduction. In other words, this
is a generalized conformally invariant part of the structure
functions of the $\hat G$-structure.

Since in the case of the Riemannian structure
$\fg_*(\fg_{-1},\hbox{ } \hat{\fg_0}) = \fg_{-1} \oplus
\hat{\fg_0}$,then there only exist  SFs of orders 1 and 2. Though
the BWB theorem doesn't work in this case, SFs are describable
thanks to the following proposition:

\begin{Proposition}{}[G2]. \emph{1)} $H^{1,2}_{\hat{\fg_0}} =
H^{1,2}_{\fg_0}$;

\emph{2)} $H^{2,2}_{\hat{\fg_0}} = H^{2,2}_{\fg_0} \oplus
S^2(\fg_{-1}^*)$.
\end{Proposition}

The Riemannian structure in the classical case of Riemannian geometry will be
considered next.

\section*{Einstein equations}

Let $G = O(n)$. In this case $\fg_1 = \fg_{-1}$ and  a
1-dimensional subspace is distinguished in $S^2(\fg_{-1})^*$. The
sections through this subspace constitute a Riemannian metric $g$
on $M$. The usual way to determine a metric on $M$ is to define a
matrix-valued function, but actually this function with values in
symmetric matrices  depends only on one functional parameter. The
values of the Riemann tensor at a point
  of $M$ constitute an $O(n)$-module $H^2(\fg_{-1},\hbox{ } \fg_*)$,
which contains a
  trivial
component. Let a section through it be denoted by $R$. This
trivial component is naturally realized as a submodule in a module
isomorphic to $S^2(\fg_{-1})^*$.

Thus, there exist two matrix-valued functions: $g$ and $R$, both
preserved by $O(n).$ Now let $R$ correspond to the Levi-Civita
connection. The process of restoring $R$ from $g$ involves
differentiations and in this way one gets a nonlinear pde, which
constitutes one of the two conditions called  \emph{Einstein
equations} [L4, LSV, LPS]:
$$
R = \lambda  g, \hbox{ where } \lambda \in R.\eqno (EE_0)
$$
The other condition is that the other component belonging to
$S^2(\fg_{-1})^*$, the traceless Ricci tensor Ric, vanishes:
$$
\hbox{ Ric }= 0.\eqno (EE_{ric})
$$

There is a close relation between $G$-structures and so-called $F$-structures,
which are also of interest, in particular, because of their application to
Penrose's geometry. This relation will now be explained.

\section*{${\bf F}$-structures and their structure functions}

Recall that the notion of $F$-structure is a generalization of the
notion of distribution, i.e., a subbundle in $TM$ and the SFs of
an $F$-structure generalize the notion of the Frobenius form [G2].

Let $V = T_mM$, $F\subset Gr_k(V)$ be a manifold with a transitive
action of a subgroup $G_F\subset GL(V)$, $\Fee(M)$ be a subbundle
of $Gr_k(TM)$, where the fiber of $Gr_k(TM)$ is $Gr_k(T_mM)$. The
bundle $\Fee(M) \longrightarrow M$ is called an \emph{F-structure
} on $M$, if for any point $m$ of $M$ there is a linear
isomorphism $I_m: V\longrightarrow T_mM$, which induces a
diffeomorphism $I_m(F) = \Fee(m)$. A submanifold $Z\subset M$ of
dimension $k$ such that $T_zZ \subset \Fee(z)$ for any $z \in Z$
is called an integral submanifold. An $F$-structure is integrable
if for any $z \in Z$ and for any subspace $V(z) \subset \Fee(z)$
there is an integral manifold $Z$ with $T_zZ = V(z)$.

SFs of an $F$-structure are defined as follows. For $f \in F$ let
$V_f \subset V$ be the subspace corresponding to $f$. Set
$$
(T_fF)_{-1} = V/V_f,\quad (T_fF)_0 = T_fF.
$$
Define
$$(T_fF)_s = ((T_fF)_{s-1} \otimes V_f^*)\cap ((T_fF)_{s-2}
\otimes S^2V_f^*)
$$
for $s > 0$, and
$$
C^{k,s}_{T_fF} = (T_fF)_{k-s} \otimes E^sV_f^*
$$

Define the differentials as in (1). Then the cohomology groups
$H^{k,s}_{T_fF}$ are naturally defined. It turns out that the
obstruction to integrability of order $k$ of an $F$-structure on a
manifold $M$ is a section of a certain vector bundle over
$\Fee(M)$ with fiber over a point $\psi \in \Fee(m)$
  isomorphic to $H^{k,2}_{T_{I_m^{-1}(\psi )}F}$.
Moreover, there exists a map $H^{k,s}_{\fg_F} \longrightarrow
H^{k,s}_{T_fF}$, where $\fg_F = Lie(G_F)$ [G2].

The relation between SFs of a $G_F$-structure and the obstructions
to integrability of an $F$-structure generalizes a theorem of
Penrose, which states that the anti-selfdual part of the Weyl
tensor on a 4-dimensional manifold with a conformal structure
vanishes if and only if $\alpha$-surfaces exist, in other words,
the metric is $\alpha$-integrable [AHS, Gi].

More precisely, for a generalized conformal structure of type $X$,
where  $X = Gr^{m+n}_m(\Cee)$, there exist  two families of $m$
and $n$-dimensional flat generators--analogues of Penrose's
$\alpha$-planes and $\beta$-planes. When neither $m$ nor $n$ is
equal to 1, i.e., the grassmannian is not a projective space, SFs
decompose into the direct sum of two components, which are
analogues of the self-dual and anti-self-dual parts of the Weyl
tensor  on a 4-dimensional manifold  with a conformal structure.
The integrability of each  of two families of generators is
equivalent to the vanishing of the corresponding component of the
SFs.

\section*{Structure functions on supermanifolds}

The necessary background on Lie superalgebras and supermanifolds is gathered
in  [L1, L2, L3, K1, M].

The classical superspaces (homogeneous compact Hermitian symmetric
superspaces), which are the super analogues of CHSS, considered by
Goncharov, are listed in [S1].

The above definitions of SFs are generalized to Lie superalgebras via the sign
  rule. However, in the super case  new phenomena appear, which have no
analogues in the classical case:

\begin{itemize}

\item Cartan prolongations of $(\fg_{-1},\hbox{ } \fg_0)$ and of
$(\Pi \fg_{-1},\hbox{ } \fg_0)$ are essentially different;

\item faithfulness of the $\fg_0$-action on $\fg_{-1}$ is violated
in natural examples of supergrassmannians of subsuperspaces in an
$(n, \hbox{ }n)$-dimensional superspace when the center
\hfil\break $\hbox{ }$ $z$ of $\fg_0$ acts trivially. This  will
be explained in \S 9.

\item the formulation of Serre's theorem and of the Proposition of
\S4 fails to be literally true for Lie superalgebras.
\end{itemize}

\section*{Description of results}

In Chapter 1 I compute the SFs for the odd analogue of the metric
on the supermanifolds and for several related $G$-structures (see
\S2 of Chapter 1). In this case $\fg_0 = Lie(G)$ is the
periplectic Lie superalgebra, the special periplectic Lie
superalgebra, or their central extensions. It turns out that
unlike  the classical case of Riemannian geometry, the
$\fg_0$-module $H^{k,2}_{\fg_0}$ is not completely reducible, and
I describe the Jordan-H\"{o}lder series for this module. Thus, my
computations show that there is no analogue of $(EE_0)$ for the
odd metric.

In Chapter 2 and Chapter 3 I obtain an explicit description of the
Spencer cohomology groups $H^{k,2}_{\fg_0}$ for simple
finite-dimensional complex classical Lie superalgebras endowed
with $\Zee$-grading of depth 1: $\fg = \oplus _{i\geq {-1}}
\fg_i$, where $\fg_0$ is the zero-th part of the grading.

It is known [K2, S2] that all such $\Zee$-gradings are of the form
$\fg = \fg_{-1} \oplus \fg_0 \oplus \fg_1$, except for the case
when $\fg$ is the special periplectic  superalgebra considered in
Chapter 1. Thus, the cohomology groups  $H^{k,2}_{\fg_0}$
constitute the space of values of SFs
  of $G$-structures corresponding to homogeneous compact Hermitian
symmetric superspaces, where $G$ is a reductive complex Lie
supergroup of classical type and $\fg_0 = Lie(G)$. The groups
$H^{k,2}_{\hat{\fg_0}}$ correspond to  structures of Riemannian
type.

An important particular case is $\fg = \fsl(m|n)$, where $m \not=
n$, corresponding to general supergrassmannians.

In Chapter 2 I consider a $\Zee$-grading of $\fg$ for which
$\fg_0$ is a reductive Lie algebra. Thus, the $\fg_0$-module
$H^{k,2}_{\fg_0}$ is completely reducible, and for $m, n > 2$
decomposes into the direct sum of two irreducible components
--- super analogues of Penrose's tensors for the usual complex
grassmannians (see Theorem 3.1 of Chapter 2).

The case $\fg =
\fsl(n|n)$ is also interesting, because  I discovered a phenomenon
which has no an analogue in the classical case. Indeed, the center $z$ of
$\fg_0$ acts trivially on $\fg_{-1}$. If one retains the same
definition of the Cartan prolongation, then it has the form of the
semidirect sum $S^*(\fg_{-1}^*)\subplus \fg_*(\fg_{-1},\hbox{ }
\fg_0/z)$ (the ideal is $S^*(\fg_{-1}^*)$) with the natural
$\Zee$-grading and Lie superalgebra structure, but this Lie
superalgebra is not a subsuperalgebra of $\fvect(\dim \fg_{-1}$)
anymore (see Theorem 1.1 of Chapter 2).

In Chapter 3 I describe the Spencer cohomology groups for the
other $\Zee$-gradings of depth 1 of $\fsl(m|n)$ and $\fpsl(n|n)$
(see Theorem 1.4 and Theorem 1.5 of Chapter 3). These theorems
show that the superspace of SFs can be not completely reducible,
and I get the answer in terms of nonsplit exact sequences of
$\fg_0$-modules.

Finally, in the cases when $m$ or $n$ are equal to 1, I get SFs of
the Lie superalgebra of vector fields  $\fvect(m|n)$ or of
divergence-free vector fields  $\fsvect(m|n)$ (see Theorem 2.1 of
Chapter 2 and Theorem 1.3 of Chapter 3).

Theorem 2.3 of Chapter 3 shows that the SFs for queer
grassmannians constitute a module looking exactly the same as that
for grassmannians of generic dimensions.

The case $\fg = \fosp(m|2n)$ is similar to Riemannian geometry.
The SFs constitute an irreducible $\fg_0$-module, which is an
analogue of the Weyl tensor, and  the superspace of the SFs for
the reduced structure  decomposes into the direct sum of three
components--the super analogues of the Weyl tensor, the traceless
Ricci tensor, and the scalar curvature. I find the highest weights
of these components (see Theorem 3.3 of Chapter 3).

Finally, I describe the Spencer cohomology groups for exceptional
Lie superalgebras $D(\alpha )$ and $AB_3$ (see Theorem 4.3 and
Theorem 5.3 of Chapter 3, respectively).

\section*{Algebraic methods}

As in the classical case (Lie theory), computation of Spencer
cohomology groups reduces  to certain problems of  representation
theory. However, in the super case computations become much more
complicated, because of the absence of complete reducibility. I
could not directly apply the usual tools for computing
(co)homology (spectral sequences and the Borel-Weil-Bott theorem)
to superalgebras and had to retreat a step  and apply these tools
to the even parts of the considered Lie superalgebras. Then, using
certain necessary conditions, I verified whether two modules over
a Lie superalgebra that could be glued into an indecomposable
module were glued or not.

My method of computing  the structure functions is based on the
Hochcshild-Serre spectral sequence [Fu]. Let $\fg = \fg_0\oplus
\fg_1$ be a Lie superalgebra and $M$ be a $\fg$-module. On the
superspace of $k$-dimensional cochains $C^k = C^k (\fg, M)$ define
a filtration:
$$
F^0C^k = C^k\supset F^1C^k\supset F^2C^k\supset \ldots
F^jC^k\supset \ldots \supset F^{k+1}C^k = 0,
$$
where $F^jC^k = \lbrace c\in C^k| c(g_1,\ldots ,g_i,\ldots ,g_k)
  = 0$ if $k-j+1$ arguments belong to $\fg_0\rbrace,\quad 0\leq j
  \leq k + 1.$ Using this filtration  define the usual corresponding
spectral sequence $E^{p,q}_r$ [GM].  Thus, $H^2(\fg, M) = \oplus
_{p+q=2}E^{p,q}_\infty = E^{2,0}_3\oplus E^{1,1}_3 \oplus
E^{0,2}_4$.

In particular, $E^{p,q}_1 = H^q(\fg_0,\hbox{ } M\otimes
S^p\fg_1^*)$ [Fu]. Since, in the case of Spencer cohomology, $\fg
= \fg_{-1} = (\fg_{-1})_0 \oplus (\fg_{-1})_1$ is a commutative
Lie superalgebra, then
$$
E^{p,q}_1 = H^q((\fg_{-1})_0,\hbox{ } \fg_*\otimes S^p(\fg_{-1})_1^*) =
H^q((\fg_{-1})_0,\hbox{ } \fg_*)\otimes  S^p(\fg_{-1})_1^*.
$$
Then in special cases I use the BWB theorem to compute
$H^q((\fg_{-1})_0,\hbox{ } \fg_*)$ as a module over $(\fg_0)_0$.

\section[Plans for the future]{Plans for the future; structure functions of Lie
superalgebras of Cartan type}

It is interesting to compute the Spencer cohomology groups
$H^{k,s}_{\fg_0}$ for simple finite-dimensional Lie superalgebras
of Cartan type with nonstandard $\Zee$-gradings of depth 1. For
example, let $\fg = \fvect(0|n)$. Define the $\Zee$--grading as
follows: $\fg = \fg_{-1} \oplus \fg_0 \oplus \fg_1$, where $\fg_0
= \Lambda (n - 1)\subplus \fvect(0|n - 1)$, $\Lambda (n - 1)$ is a
commutative ideal in $\fg_0$ isomorphic to the Grassmann algebra
as the $\fvect (0|n - 1)$-module, $\fg_1 \cong \Pi (\fvect (0|n -
1))$ as the $\fg_0$-module, $\fg_{-1} \cong \Pi (\Lambda (n - 1))$
as the $\fg_0$-module, and the action of $\Lambda (n - 1)$ on
$\fg_{-1}$ is the usual multiplication in $\Lambda (n - 1)$.

It is much more difficult to compute the Spencer cohomology groups for
Cartan Lie superalgebras than for classical ones, because in this case
the number of the irreducible quotient modules in the Jordan-H\"{o}lder series
does depend on $n$. Moreover, even for small $n$ computations seem to be very
complicated, because of the absence of complete reducibility [LPS].

\section*{Terminological conventions}

$\fcg$ is the trivial central extension of the Lie superalgebra
$\fg$;

$\subplus$ stands for the semidirect sum of the algebras with the
ideal to the left;

$\Pi $ is the functor of the change of parity;

$p(X)$ is the parity of homogeneous element $X$ of a superspace;

$S^iV$ is the $i$-th symmetric power of a vector (super)space $V$;

$S^*V = \oplus_{k\geq 0}S^kV$;

$E^iV$ is the $i$-th exterior power of a vector (super)space $V$;

$\langle 1\rangle$ is a trivial module over a Lie superalgebra;

I will consider the following Lie superalgebras  [L2, L3]:

$\fgl(m|n) = Mat(m|n;\hbox{ }\Cee)$,

$\fsl(m|n) = \lbrace X\in \fgl(m|n) |\hbox{ } str X = 0\rbrace $,

$\fpsl(n|n) = \fsl(n|n)/\langle 1_{2n} \rangle$, where $1_{2n}$
is the unit matrix;

$\fvect(m|n) = \fder \Cee[[x]]$ is the Lie superalgebra of vector
fields;

$\fsvect(m|n) = \lbrace D\in \fvect(m|n) |\hbox{ } \div(D) =
0\rbrace$ is the Lie superalgebra of divergence-free vector
fields;

$\fh(0|n) = \lbrace D\in \fvect(0|n) |\hbox{ } Dw = 0$ for the
Hamiltonian form $w = \sum_{i=1}^n(d\xi _i)^2\rbrace$ is the Lie
superalgebra of Hamiltonian vector fields;

Other Lie superalgebras will be defined in  Chapter 1 and Chapter
3.

%%%%%%%%%%%%%%%%%%%%%%%%%%%%%%%%%%%%%%%%%%%%%%%%%%%%%%%%%%%%%%%%%%%%%%%%%%%%%%%%
%%%%%%%%%%%%%%%%%%%%%%%%%%%%%%%%%%%%%%%%%%%%%%%%%%%%%%%%%%%%%%%%%%%%%%%%%%%%%%%%%%%%%%%%%%%%%%%%%%%%%%%%%%%%%%%%%
%Chapter 1
%%%%%%%%%%%%%%%%%%%%%%%%%%%%%%%%%%%%%%%%%%%%%%%%%%%%%%%%%%%%%%%%%%%%%%%%%%%%%%%%%%%%%%%%%%%%%%%%%%%%%%%%%%%%%%%%%%%%
%%%%%%%%%%%%%%%%%%%%%%%%%%%%%%%%%%%%%%%%%%%%%%%%%%%%%%%%%%%%%%%%%%%%%%%%%%%%%%%%
\chapter[The analogues of the  Riemannian tensors]
{The analogues of the  Riemannian tensors for the
odd metric on supermanifolds}

\section*{Periplectic superalgebras and their Cartan prolongations}

Let $z = 1_{2n}$ be the unit matrix and $\tau = \diag (1_n,\hbox{
}-1_n)$.

Let $P$ be a nondegenerate supersymmetric
odd bilinear form on a superspace $V$.  Clearly, $\dim V =
(n,\hbox{ } n)$. Define the odd analogue of the symplectic Lie
algebra, the periplectic Lie superalgebra $\fpe(n)$, and its
special subsuperalgebra $\fspe(n)$, setting
$$
\eqalign{
&\fpe(n)=\{ X\in \fgl(n|n)|\hbox{ } X^{st}P+(-1)^{p(X)}PX=0\},
\cr &\fspe(n)=\fpe(n)\cap \fsl(n|n). }
$$

Thus,
$$
\fpe(n)= \fspe(n)\subplus \langle\tau \rangle.
$$

Denote by $\eps_1,\ldots,\eps_n$ the standard basis
  of the space dual to the space of diagonal matrices in
  $\fgl(n)\subset \fpe(n).$ Denote by $V_\lambda$
the irreducible $\fgl(n)$-module with highest weight $\lambda$
and highest vector $v_\lambda$ and by $X_\lambda$ the irreducible
$\fpe(n)$-module with highest weight $\lambda$ and an even highest
vector.

  Let $V=V_0\oplus V_1$
be the standard (identity) $\fpe(n)$-module, $e_1,\ldots,e_n$ be a
basis of $V_0$,
  and $f_1,\ldots,f_n$ be a basis of $V_1$
with respect to which the form $P$ on
$V$ takes the form $P = \antidiag (1_n, 1_n$). With respect to this
basis the elements
$X\in \fpe(n)$ are represented by matrices of the standard format
$(n,\hbox{ }n)$ [L2]:
$$
X =\diag(A,-A^t) +\antidiag(B,C), \hbox{ where } A\in
\fgl(n),\hbox{ }
  B^t = B,\hbox{ }  C^t = -C.$$

  In what follows we will often use
  a natural abbreviation: e.g.,
$B_{1,n}$ stands for the matrix $X$
  whose components $A$ and $C$ are zero and
  all the entries of $B$ are also zero except
  for $(1,\hbox{ }n)$-th and $(n,\hbox{ }1)$-st.

  Denote by $\tilde e_1,\ldots,\tilde e_n$ and
  $\tilde f_1,\ldots,\tilde f_n$ the basis of
  $V^*$ dual to the above basis of $V$,
  i.e., $\tilde f_i(e_j) = \tilde e_i(f_j)
  = \delta _{ij}.$ Since the form $P$ preserved
by $\fpe(n)$ is odd, then $V^*$ and $\Pi (V)$ are isomorphic as
$\fpe(n)$-modules.
  Notice that as $\fpe(n)$-modules,
$$
\fpe(n) \cong \Pi (E^2V^*).
$$

\ssbegin{1.1}{Lemma} a) There exists a $\Zee$-grading of the Lie superalgebra
$\fpe(n+1)$ of the form
$$\fg_{-1}\oplus \fg_0\oplus \fg_1\oplus \fg_2,$$
where
$$
\eqalign {&\fg_{-1} = V,\cr
&\fg_0 = \fcpe(n),\cr &\fg_1 = V^* = \Pi (V),\cr &\fg_2 = \Pi
(\langle 1\rangle).\cr}
$$

b) There exists a $\Zee$-grading of the Lie
superalgebra $\fspe(n+1)$ of the form
$$
\fg_{-1}\oplus \fg_0\oplus \fg_1\oplus \fg_2,
$$
where
$$
\eqalign {
&\fg_{-1} = V,\cr &\fg_0 = \fspe(n)\subplus \langle \tau
+nz\rangle,\cr &\fg_1 = V^*
= \Pi (V),\cr &\fg_2 = \Pi (\langle 1\rangle).\cr}
$$
\end{Lemma}

\begin{proof}
Let $W = W_0\oplus W_1$ be the standard (identity)
$\fpe(n+1)$-module, $e_1,\ldots, e_{n+1}$ be a basis of
$W_0$, and $f_1, \ldots, f_{n+1}$ be a basis of $W_1$ with respect
to which the form $P$ on $W$ takes the form $P = \antidiag
(1_{n+1},\hbox{ }1_{n+1}).$ Denote by $\tilde e_1,\ldots,\tilde
e_{n+1}$ and $\tilde f_1\ldots,\tilde f_ {n+1}$ the basis of $W^*$
dual to the above basis of $W$, e.g., $\tilde f_i(e_j) = \tilde
e_i(f_j) = \delta _{ij}.$

Note that
$$
\fpe(n+1) = \Pi (E^2W^*) = \Pi (E^2W_1 \oplus W_0 \wedge W_1 \oplus S^2W_0).
$$
Thus,
$$
\fpe(n+1) = \langle e_i\tilde e_j, e_i \wedge \tilde f_j, f_i \wedge
\tilde f_j\rangle
(1\leq i,\hbox{ } j \leq n + 1), \hbox{ where }
$$
$$
\begin{array}{l}
   e_i\tilde e_j = (1/2)(e_i \otimes \tilde e_j + e_j \otimes\tilde e_i), \\
   e_i \wedge \tilde f_j = (1/2)(e_i \otimes \tilde
   f_j - f_j\otimes \tilde e_i),\\
   f_i \wedge \tilde f_j =
   (1/2)(f_i \otimes \tilde f_j - f_j \otimes \tilde f_i).
\end{array}
\eqno{(1.1.1)}
$$
Note that the commutator in $\fpe(n+1)$ is defined as follows:
$$
[w_i \wedge \tilde w_j,\hbox{ } w_s \wedge \tilde w_t] =
(1/2)(\tilde w_j(w_s)(w_i \wedge \tilde w_t) - (-1)^{p(w_s)}\tilde w_t(w_i)
(w_s \wedge \tilde w_j)), \eqno (1.1.2)
$$
where $w_l \in \lbrace e_1, \ldots , e_{n+1}; f_1, \ldots , f_{n+1} \rbrace
\hbox{ for } 1 \leq l \leq n + 1.$

Let $V = V_0 \oplus V_1 = \langle e_1, \ldots , e_n; f_1, \ldots ,
f_n \rangle$ be
$(n,\hbox{ }n)$-dimensional subsuperspace in $W$. Then
$$
\begin{array}{l}
   \Pi (E^2W_1) = \Pi (E^2V_1) \oplus V_1 \wedge \langle \tilde
f_{n+1} \rangle,\\
   \Pi (W_0 \wedge W_1) = \Pi (V_0 \wedge V_1) \oplus V_0 \wedge
   \langle \tilde f_{n+1} \rangle \oplus V_1 \wedge \langle \tilde
e_{n+1} \rangle \oplus
   \langle e_{n+1} \wedge \tilde f_{n+1} \rangle,\\
   \Pi (S^2W_0) = \Pi (S^2V_0) \oplus V_0 \wedge \langle \tilde
e_{n+1} \rangle\oplus
   \langle e_{n+1}\tilde e_{n+1} \rangle.
\end{array}
$$

Set
$$
\begin{array}{l}
   \fg_{-1} = V_0 \wedge \langle \tilde f_{n+1} \rangle \oplus V_1 \wedge
   \langle\tilde f_{n+1} \rangle, \\
   \fg_0 = \Pi (E^2V_1 \oplus V_0 \wedge V_1
   \oplus S^2V_0) \oplus \langle e_{n+1} \wedge \tilde f_{n+1}
   \rangle, \\
   \fg_1 = V_0 \wedge \langle \tilde e_{n+1} \rangle \oplus V_1
   \wedge \langle \tilde e_{n+1} \rangle, \\
   \fg_2 = \langle e_{n+1}\tilde e_{n+1}
   \rangle.
\end{array}
\eqno{(1.1.3)}
$$

According to (1.1.2), formulas (1.1.3) indeed define a
$\Zee$-grading of $\fpe(n+1)$, described in Lemma 1.1.

In order to define a $\Zee$-grading of $\fspe(n + 1)$ we set
$$
\begin{array}{l}
   \fg_0 = \Pi (E^2V_1 \oplus S^2V_0) \oplus \langle
   \sum_{i,j=1,n}a_{ij}e_i \wedge \tilde f_j | \sum_{i=1,n}a_{ii} = 0
   \rangle\oplus\\
   \langle (\sum_{i=1,n}e_i \wedge \tilde f_i) - ne_{n+1}
   \wedge \tilde f_{n+1} \rangle.
\end{array}
$$
Note that by (1.1.2) we have
$2[e_{n+1} \wedge \tilde f_{n+1}, \hbox{ }\fg_i] = i\fg_i$ for
$(-1 \leq i \leq 2).$ Hence, $-2e_{n+1} \wedge \tilde
f_{n+1} = z$. Since by (1.1.1) $2\sum_{i=1,n}e_i \wedge\tilde f_i
= \tau $, then
$$
2(\sum_{i=1,n}e_i \wedge \tilde f_i - ne_{n+1} \wedge \tilde f_{n+1}) =
\tau + nz.
$$
Thus, $\fg_0 = \fspe(n) \subplus \langle \tau + nz \rangle$.
This proves Lemma 1.1.
\end{proof}

\ssbegin{1.2}{Theorem} Let $\fg_{-1} = V.$ Then

  a) If $\fg_0 =
\fspe(n)$, $\fpe(n)$, $\fcspe(n)\hbox { or } \fspe(n)\subplus \langle a\tau
+bz\rangle$, where $a$, $b \in \Cee$ are such that $a$, $b \not= 0$ and
$b/a \not= n$, then $\fg_*(\fg_{-1},\hbox{ }\fg_0)
  = \fg_{-1}\oplus \fg_0.$

b) If $\fg_0 = \fcpe(n)\hbox { or } \fspe(n)\subplus \langle\tau
+nz\rangle$,
  then $\fg_*(\fg_{-1},\hbox{ }\fg_0 )$ is
either $\fpe(n+1)$ or $\fspe(n+1)$, respectively, in the
$\Zee$-grading described in Lemma 1.1.
\end{Theorem}

\begin{proof}
Let us consider the case where $\fg_0 = \fcpe(n).$ By Lemma 1.1
we have
$$
\fpe(n+1) = \fg_{-1}\oplus \fg_0\oplus
  \fg_1\oplus \fg_2,
$$
where
$$
\eqalign{&\fg_{-1} = V,\cr
&\fg_0 = \fcpe(n),\cr &\fg_1 = V^* = \Pi (V),\cr &\fg_2 = \Pi
(\langle 1\rangle).}
$$
Therefore,
$$
\fpe(n+1) \subset \fg_*(\fg_{-1},\hbox{ } \fg_0). \eqno(1.2.1)
$$
In fact, since $\fspe(n+1)$ is a simple Lie superalgebra, then it
is transitive, (i.e., if there exists $g \in \fg_i (i \geq 0)$
such that $[\fg_{-1},\hbox{ }g] = 0$, then $g = 0$). It follows
that $\fg_i \subset \fg_{i-1} \otimes \fg_{-1}^*$. The Jacobi
identity implies $\fg_i \subset \fg_{i-2} \otimes S^2\fg_{-1}^*$.
\end{proof}

Let us find $\fg_1$.

\ssbegin{1.3}{Lemma} As a
$\fgl(n)$-module, $\fg_0\otimes \fg_{-1}^*$ is the direct sum of
irreducible $\fgl(n)$-submodules whose highest weights and
highest vectors are listed in Table 1.
\end{Lemma}

\begin{rem*}{Convention} Let $v, w$ be elements of a vector space. Set
$$
vw = (v\otimes w + (-1)^{p(v)p(w)}w\otimes v)/2,\hbox{ } v\wedge w =
(v\otimes w - (-1)^{p(v)p(w)}w\otimes v)/2.
$$
\end{rem*}
\begin{proof}[Proof of Lemma 1.3] consists of:

a) a
verification of the fact that vectors $v$ from Table 1 are indeed
highest with respect to $\fgl(n)$, i.e., $A_{i,j}v = 0 \hbox{ for
}i < j$,

b) a calculation of dimensions of the
corresponding irreducible  submodules by the formula from the
Appendix.
\end{proof}

Let us show with the help of Table 1 that if $\lambda \not= \eps
_1,\hbox{ }-\eps _n$, then $v_\lambda  \not\in \fg_{-1}\otimes
S^2\fg_{-1}^*.$ Indeed, if $\lambda = -\eps _{n-1}-2\eps _n$, then
$$v_\lambda (e_n)(e_{n-1}) = -f_n/2, v_\lambda(e_{n-1})(e_n) =0;$$
if $\lambda = -\eps_ {n-2}-\eps _{n-1}-\eps _n$, then
$$v_\lambda (e_n)(e_{n-1}) = f_{n-2}/2, v_\lambda (e_{n-1})(e_n) =
-f_{n-2}/2;$$
if $\lambda = \eps _1-2\eps _n $, then
$$v_\lambda (e_n)(f_1) = f_n/2, v_\lambda (f_1)(e_n) = 0;$$
if $\lambda = \eps _1+\eps _2-\eps _n $, then
$$v_\lambda (f_2)(e_n) = -e_1/2, v_\lambda (e_n)(f_2) = 0;$$
if $\lambda = 3\eps _1$, then
$$v_\lambda (f_1)(f_1) = e_1\not= 0;$$
if $\lambda = 2\eps_ 1+\eps _2$, then
$$
v_\lambda (f_1)(f_1) = e_2/2\not= 0.
$$
Let $\lambda = \eps _1-\eps _{n-1}-\eps _{n}.$ According to Table
1, $\fg_0\otimes \fg_{-1}^*$ contains two highest vectors of
weight $\lambda$. Let
$$
v_\lambda = k_1f_{n-1}\wedge \tilde f_n\otimes \tilde e_1
+k_2(f_{n-1}\wedge \tilde e_1\otimes \tilde f_n-f_n\wedge\tilde
e_1 \otimes\tilde f_{n-1}), \hbox{ where } k_1,\hbox{ } k_2\in
\Cee,
$$
be a linear combination of these vectors. The condition
$$
v_\lambda (e_n)(e_{n-1}) = v_\lambda (e_{n-1})(e_n)
$$
implies $k_2 = 0 $. Then the condition
$$
v_\lambda (e_n)(f_1) = v_\lambda (f_1)(e_n)
$$
implies $k_1 = 0$.

Let $\lambda = 2\eps  _1-\eps _n$.
Let
$$
v_\lambda = k_1f_n\wedge \tilde e_1\otimes \tilde e_1
+ k_2e_1\tilde e_1\otimes\tilde f_n, \hbox{ where } k_1,\hbox{ }
k_2\in \Cee,
$$
be a linear combination of the highest vectors of
weight $\lambda$ which belong to $\fg_0\otimes \fg_{-1}^*$.
The condition $v_\lambda(f_1)(f_1) = 0$ implies $k_1
= 0$. Then the condition $v_\lambda(e_n)(f_1) =
v_\lambda(f_1)(e_n)$ implies $k_2 = 0$. Therefore, if $\lambda
\not=\eps _1,\hbox{ }-\eps _n $, then $v_\lambda\not\in
\fg_{-1}\otimes S^2\fg_{-1}^*$, hence $v_\lambda\not\in \fg_1.$

Let $\lambda = \eps _1$. According to Table 1,
$\fg_0\otimes \fg_{-1}^*$ has four highest vectors of weight
$\lambda$. Let
$$
v_\lambda = k_1\sum_{i=1}^n f_i\wedge \tilde e_i\otimes \tilde e_1
+ k_2\sum_{i=1}^n f_i\wedge \tilde e_1\otimes \tilde e_i +
k_3\sum_{i=1}^n e_1\tilde e_i\otimes \tilde f_i + k_4\sum_ {i=1}^n
e_i\tilde f_i\otimes \tilde e_1,
$$ where $k_1$, $k_2$, $k_3$, $k_4
\in \Cee,$ be their linear combination. Note that $v_\lambda \in
\fg_{-1}\otimes S^2 \fg_{-1}^*$ if and only if the following
conditions are satisfied:
$$
\eqalign{
&v_\lambda(f_1)(f_1) = 0 ,\cr
&v_\lambda (f_1)(f_i) = -v_\lambda (f_i)(f_1) \hbox{ for } i\not= 1,\cr
&v_\lambda (f_1)(e_1) = v_\lambda (e_1)(f_1),\cr
&v_\lambda (f_1)(e_i) = v_\lambda (e_i)(f_1)\hbox{ for } i\not= 1,\cr
&v_\lambda (f_i)(e_i) = v_\lambda (e_i)(f_i)\hbox{ for } i\not= 1,\cr
}
$$
which determine, respectively, the following system of linear equations:
$$
\eqalign{
&k_1 + k_2 + k_4 = 0, \cr
&k_1 + k_4 = -k_2,\cr
&(-k_1 -k_2 + k_4)/2 = k_3,\cr
& -k_1 + k_4 = k_3,\cr
&-k_2 = k_3.\cr
}
$$
The solution of this system is
$$
k_1 = 0,\hbox{ } k_3 = -k_2 = k_4.
$$
Therefore,
$$
v_{\eps _1} = -\sum_{i=1}^n f_i\wedge \tilde e_1\otimes \tilde e_i
+ \sum_{i=1}^n e_1\tilde e_i\otimes\tilde f_i + \sum_{i=1}^n
e_i\tilde f_i\otimes \tilde e_1 \in \fg_1.  \eqno  (1.3.1)
$$ Since
$\fg_1$ is a $\fpe(n)$-module and $v_{\eps _1}$ is an odd vector,
$\fg_1 = V^*.$

Let us find $\fg_2$.

\ssbegin{1.4}{Lemma} There exist
the following nonsplit sequences of $\fpe(n)$-modules:
$$
0\longrightarrow X_{2{\eps _1}} \longrightarrow E^2V^*
\longrightarrow  \Pi(\langle1\rangle) \longrightarrow 0,  \eqno (1.4.1)
$$
$$
0\longrightarrow  \Pi (\langle1\rangle)\longrightarrow S^2V^*\longrightarrow
X_{\eps _1+ \eps_2}\longrightarrow 0.\eqno (1.4.2)
$$
\end{Lemma}

\begin{proof} First of all recall that $\Pi (E^2V^*)$ and $\fpe(n)$
itself are isomorphic $\fpe(n)$-modules and there exists the
following nonsplit sequence of $\fpe(n)$-modules
$$
0\longrightarrow \fspe(n)\longrightarrow \fpe(n)\longrightarrow
\langle\tau\rangle
\longrightarrow 0.\eqno (1.4.3)
$$
Note that as a
$\fgl(n)$-module, $E^2V^*$ is isomorphic to
$$
E^2V_0^*\oplus V_0^*\wedge V_0\oplus S^2V_0,\hbox{ where }
S^2V_0 = V_{2{\eps _1}}, v_{2{\eps_1}} = \tilde e_1^2. \eqno
(1.4.4)
$$
Since $B_{i,j}v_{2{\eps _1}} = 0,$ then $v_{2{\eps _1}}$
is a $\fpe(n)$-highest vector. From the simplicity of $\fspe(n)$
we get (1.4.1) after the change of parity.

Let us prove (1.4.2). Notice that as $\fgl(n)$-modules,
$$
S^2V^* = S^2V_0^*\oplus V_0^*\cdot V_0\oplus E^2V_0,\hbox { where }
$$
$$
\begin{array}{l}
   S^2V_0^* = V_{-2{\eps_n}}, v_{-2{\eps_n}} = \tilde f_n^2,\\
   V_0^*\cdot V_0 = V_{\eps _1-\eps_n}\oplus \langle v_0\rangle, \hbox  { }
   v_{\eps_1-\eps_n} = \tilde f_n \tilde e_1, \hbox  {  } v_0 =
   \sum_{i=1}^n \tilde f_i \tilde e_i,\\
   E^2V_0 =
   V_{\eps_1 + \eps_2},\hbox {  } v_{\eps_1 + \eps_2} = \tilde
   e_1\wedge \tilde e_2.
\end{array}
\eqno{(1.4.5)}
$$
Note that $\langle v_0\rangle$ is the trivial
1-dimensional $\fpe(n)$-module. Indeed,
$$
B_{i,j}(v_0) = -\tilde e_j\wedge \tilde e_i +
-\tilde e_i\wedge \tilde e_j = 0,\hbox{ } C_{i,j}(v_0) = \tilde
f_i \tilde f_j + \tilde f_j(-\tilde f_i) = 0.
$$
Since
$B_{i,j}v_{\eps_1 + \eps_2} = 0$, then $v_{\eps_1 + \eps_2}$ is a
$\fpe(n)$-highest vector.

Let us prove that as $\fgl(n)$-modules,
$$
X_{\eps_1 + \eps_2}  \cong V_{\eps_1 + \eps_2}
\oplus V_{\eps_1-\eps_n}\oplus V_{-2\eps_n}.
$$
Indeed,
$$
\eqalign {&C_{2,n}(v_{\eps_1 + \eps_2}) = v_{\eps_1-\eps_n},
\hbox{ } -B_{2,n}(v_{\eps_1-\eps_n}) = v_{\eps_1 + \eps_2},\cr
&C_{1,n}(v_{\eps_1-\eps_n}) = v_{-2\eps_n}, \hbox{ }
B_{1,n}(v_{-2\eps_n}) = v_{\eps_1-\eps_n}.\cr}
$$
Finally, we get
$$
v_0 =\frac12\sum_{i=1}^{n-1}B_{i,n}A_{n,i} + nB_{n,n}
v_{-2\eps_n}.
$$
This proves (1.4.2) and Lemma 1.4.
\end{proof}

Let us prove that if $v_\lambda$ from $\fg_1\otimes \fg_{-1}^*$ is
a $\fpe(n)$-highest vector of weight either $\lambda = 2\eps_1$ or
$\eps_1 + \eps_2$, then $v_\lambda \not\in \fg_0\otimes
S^2\fg_{-1}^*$. In fact, if $\lambda = 2\eps_1$, then by (1.3.1)
and (1.4.4) we have

$$
v_\lambda = (-\sum_{i=1}^n f_i\wedge \tilde e_1\otimes \tilde e_i
  + \sum_{i=1}^n e_1\tilde e_i\otimes\tilde f_i
+ \sum_{i=1}^n e_i\tilde f_i\otimes\tilde e_1)\otimes\tilde e_1.
$$
Then
$$
v_\lambda(f_1)(f_1) = -f_1\wedge\tilde e_1 + \sum_{i=1}^n e_i\tilde f_i
\not= 0.
$$
Therefore, $v_\lambda \not\in \fg_0\otimes
S^2\fg_{-1}^*$.

If $\lambda = \eps_1 + \eps_2$, then
by (1.3.1) and (1.4.5)
$$
\eqalign{
&v_\lambda = (-\sum_{i=1}^n f_i\wedge\tilde e_1\otimes \tilde e_i
+ \sum_{i=1}^n e_1\tilde e_i\otimes\tilde f_i + \sum_{i=1}^n
e_i\tilde f_i\otimes\tilde e_1)\otimes\tilde e_2 \cr
&-(-\sum_{i=1}^n f_i\wedge\tilde e_2\otimes \tilde e_i +
\sum_{i=1}^n e_2\tilde e_i\otimes \tilde f_i + \sum_{i=1}^n
e_i\tilde f_i\otimes\tilde e_2)\otimes\tilde e_1.\cr }
$$
Thus,
$$
v_\lambda(f_2)(e_i) = e_1\tilde e_i \not= 0  \hbox{ and }
v_\lambda(e_i)(f_2) = 0.
$$
Hence, $v_\lambda \not\in \fg_0\otimes
S^2\fg_{-1}^*.$

Let $\lambda = 0$. According to Lemma 1.4, the Jordan-H\"{o}lder series
  of the $\fpe(n)$-module $\fg_1\otimes \fg_{-1}^*$ contains two
  $\fpe(n)$-modules with highest weight 0.
By Lemma 1.4 the sequence (1.4.1) is nonsplit and we have already
proved that $\fg_2$ has no irreducible $\fpe(n)$-module with
highest weight $2\eps_1$. Therefore, either $\fg_2$ consists of
one trivial $\fpe(n)$-module or $\fg_2 = 0.$ But by (1.2.1)
$$
\fpe(n+1)\subset  \fg_*(\fg_{-1},\hbox{ }\fg_0).
$$
Hence, $\fg_2 = \Pi (\langle 1\rangle)$.

Finally, let us show that $\fg_3 = 0$. By definition
$$
\fg_3 = (\fg_2\otimes \fg_{-1}^*)\cap (\fg_1\otimes S^2\fg_{-1}^*).
$$
Note that
$$
\fg_2\otimes \fg_{-1}^* = \Pi (\langle 1\rangle)\otimes V^* \cong V,
$$
as $\fpe(n)$-modules. By (1.1.3) the $\fpe(n)$-highest vector in
$\fg_2\otimes \fg_{-1}^*$ is
  $v = e_{n+1}\tilde e_{n+1}\otimes\tilde e_1.$

By the explicit formula (1.1.2) of multiplication in
$\fpe(n+1)$ we have
$$
v(f_1)(e_1) = [e_{n+1}\tilde e_{n+1},e_1\wedge\tilde
f_{n+1}] = e_1\tilde e_{n+1} \in \fg_1.
$$
On the other hand,
$v(e_1)(f_1) = 0.$ Therefore, $v\not\in \fg_1\otimes S^2
\fg_{-1}^*.$ Hence, $\fg_3 = 0$.

Thus, Theorem 1.2 is
proved for $\fg_0 = \fcpe(n).$ This result and part b) of Lemma
1.1 imply the statement of Theorem 1.2 for $\fg_0 =\fspe(n)\subplus
\langle\tau + nz\rangle.$

Let us prove that $\fg_1 = 0$ for $\fg_0 = \fspe(n)$, $\fpe(n)$,
$\fcspe(n), \hbox{ or }\fspe(n)\subplus \langle a\tau + bz\rangle$,
where  $a, \hbox{ }b\in \Cee$, $a,\hbox{ } b\not=0$, and
$b/a\not=n$.

Indeed, as has been shown, $\fg_1 = \Pi
(V) \hbox{ for } \fg_0 = \fcpe(n) \hbox{ or } \fspe(n)\subplus \langle\tau
+nz\rangle$, and by $(1.3.1)$ the corresponding $\fspe(n)$-highest
vector is
$$
v_{\eps_1} = -\sum_{i=1}^n f_i\wedge \tilde e_1\otimes\tilde e_i
+ \sum_{i=1}^ne_1\tilde e_i\otimes\tilde f_i
+ \sum_{i=1}^n e_i\tilde f_i\otimes\tilde e_1.
$$
Then
$$
v_{\eps_1}(f_1) = -f_1\wedge \tilde e_1 + \sum_{i=1}^n  e_i
\tilde f_i \in \fspe(n)\subplus \langle\tau + nz\rangle.
$$
Note that
$$
-f_1\wedge\tilde e_1 + \sum_{i=1}^n e_i\tilde f_i \not\in  \fspe(n).
$$
Hence $v_{\eps_1} (f_1) \not\in \fg_0 \hbox{ for } \fg_0 =
\fspe(n)\subplus \langle a\tau + bz\rangle, \hbox{ where }a,\hbox{
}b\in \Cee,\hbox{
} a\not=0 \hbox{ and }b/a\not=n, \hbox{ or } a=0.$
Therefore,
$$
\fg_1 = 0 \hbox{ for } \fg_0 =\fspe(n), \hbox{ } \fpe(n),\hbox{ }
\fcspe(n), \hbox{ or }
\fspe(n)\subplus\langle a\tau + bz\rangle,
$$
  where $a,\hbox{ } b \in \Cee,\hbox{ }
a,\hbox{ } b\not= 0 \hbox{ and }b/a \not= n.$ Thus, in these cases
$$
\fg_*(\fg_{-1},\hbox{ } \fg_0) = \fg_{-1}\oplus \fg_0.
$$

\section*{The main theorem}

%%%%%%%%%%%%%%%%%%%%%%%%%%%%%%%%%%%%%%%%%%%%%%%%%%%%
%%%%%%%%%%%%%%%%%%%%%%%%%%%%%%%%%%%%%%%%%%%%%%%%%%%%%%%%%%%%%
%\S2%%%%%%%%%%%%%%%%%%%%%%%%%%%%%%%%%%%%%%%%%%%%%%%%%%%%%%%%%%%%%%%%%%%%%
%%%%%%%%%%%%%%%%%%%%%%%%%%%%%%%%%%%%%%%%%%%%%%%%%%%%%%%%%

The following theorem describes SFs for the odd analogues of
Riemannian metric and various conformal versions.

\ssbegin{2.1}{Theorem} For the $G$-structures with the following
Lie $(G)  = \fg_0$ the nonzero SFs are of orders not exceeding 2
and as follows:

\emph{order 1:} if $\fg_0
= \fspe(n)$ or $\fspe(n)\subplus \langle\tau + nz\rangle $ we have $V^*$;
\emph{order 2:}if $\fg_0 = \fspe(n)$,
where $n>3$, we have the following nonsplit exact sequence of
$\fspe(n)$-modules:
$$
0\longrightarrow X_{\eps_1 + \eps_2}\longrightarrow H^{2,2}_{\fg_0}
\longrightarrow \Pi (X_{2\eps_1 + 2\eps_2})\longrightarrow 0;
$$
if $\fg_0 = \fspe(3)$  another space  is added to the SFs: we have
the following nonsplit exact sequence of $\fspe(3)$-modules:
$$
0\longrightarrow X\longrightarrow H^{2,2}_{\fg_0}\longrightarrow \Pi (X_
{3\eps_1})\longrightarrow 0,
$$
where $X$ is determined from the
following nonsplit exact sequence of $\fspe(3)$-modules:
$$
0\longrightarrow X_{\eps_1 + \eps_2}\longrightarrow X
\longrightarrow \Pi (X_{2\eps_1 + 2\eps_2}) \longrightarrow 0;
$$
if $\fg_0 = \fspe(n)\subplus \langle a\tau + bz\rangle$, where
$a,\hbox{ } b \in
\Cee$ are such that $a = 0,\hbox{ } b\not= 0$ or $a\not= 0,\hbox{
} b/a \not= n$, then for $n>2$ we have the following nonsplit
exact sequence of $\fspe(n)$-modules:
$$
0\longrightarrow H^{2,2}_{\fspe(n)}\longrightarrow H^{2,2}_{\fg_0}
\longrightarrow X_{2\eps_1}\longrightarrow 0;
$$
if $\fg_0 =
\fspe(n)\subplus \langle\tau + nz\rangle$, then for $n>3\hbox{ }
H^{2,2}_{\fg_0} =
  \Pi (X_{2\eps_1 + 2\eps_2})$ is an irreducible $\fspe(n)$-module,
for $n = 3$ we have the following nonsplit exact
sequence of $\fspe(3)$-modules:
$$
0\longrightarrow \Pi (X_{2\eps_1 + 2\eps_2})
\longrightarrow H^{2,2}_{\fg_0} \longrightarrow \Pi (X_{3\eps_1})
\longrightarrow  0;
$$
if $\fg_0 = \fcpe(n),\hbox{ } n>2$, then
$$
H^{2,2}_{\fcpe(n)} =
\Pi (S^2(E^2V/\Pi (\langle 1\rangle))/E^4V),
$$
more precisely, we have the
following nonsplit exact sequence of $\fspe(n)$-modules:
$$
0\longrightarrow H^{2,2}_{\fspe(n)\subplus \langle\tau + nz\rangle}
\longrightarrow H^{2,2}_{\fg_0} \longrightarrow X_{2\eps_1}
\longrightarrow 0.
$$

\end{Theorem}

%%%%%%%%%%%%%%%%%%%%%%%%%%%%%%%%%%%%%%%%%%%%%%%%%%%%%%%%%%%%%%%%%%%%%%%%%%%%%%%%%%%%%%%%%%%%%%%%%%%%%%%%%%%
%\S3
%%%%%%%%%%%%%%%%%%%%%%%%%%%%%%%%%%%%%%%%%%%%%%%%%%%%%%%%%%%%%%%%%%%%%%%%%%%%%%%%%%%%%%%%%%%%%%%%%%%%%%%%%%%%%%%%%%%%%%%%%%%%%%%%
\section*{Proof of the main theorem}

\ssec{3.1. Calculation of SFs of order 1.} Recall that the bidegree of the
differentials in the Spencer complex is $(-1,\hbox{ } 1)$.
  We will often refer to the following
\begin{Lemma}
Let $(\fg_{-1}, \hbox{ }\fg_0)$ be an
arbitrary pair, where $\fg_{-1}$ is a faithful module over a Lie
superalgebra $\fg_0$, and let
$$
\fg_{k-1}\otimes \fg_{-1}^* \buildrel{\partial^{k+1,1}_{\fg_0}}
\over \longrightarrow\fg_{k-2} \otimes E^2\fg_{-1}^* \buildrel
{\partial^{k,2}_{\fg_0}} \over \longrightarrow \fg_{k-3} \otimes E
^3\fg_{-1}^* \eqno (k \geq 1)
$$
be the  corresponding Spencer cochain sequence. Then
$$
\Im \partial^{k+1,1}_{\fg_0} \cong (\fg_{k-1} \otimes \fg_{-1}^*) /\fg_k.
\eqno(3.1.1)
$$
\end{Lemma}
\begin{proof} By the definition of the Cartan prolongation
$$
\fg_k = (\fg_{k-1}\otimes \fg_{-1}^* )\cap( \fg_{k-2} \otimes S^2\fg_{-1}^*).
$$
Let $c \in \fg_{k-1} \otimes \fg_{-1}^*$. Then $c\in \fg_k$ if and
only if $c(g_1)g_2 = (-1)^{p({g_1})p({g_2})} c(g_2)g_1$ for any
(homogeneous) $g_1, g_2 \in \fg_{-1}$. On the other hand
$$
\partial^{k+1,1}_{\fg_0} c(g_1, g_2) = -(-1)^{p({g_1})p({g_2})}
c(g_2)g_1 + c(g_1)g_2.
$$
Hence, $\Ker \partial^{k+1,1}_{\fg_0} =
\fg_k$. This proves the Lemma.
\end{proof}

In particular, to define $H^{1,2}_{\fg_0}$ we have the following
Spencer cochain sequence:
$$
\fg_0\otimes \fg_{-1}^*\buildrel {\partial^{2,1}_{\fg_0}}\over \longrightarrow
\fg_{-1}\otimes E^2\fg_{-1}^*\buildrel
{\partial^{1,2}_{\fg_0}}\over \longrightarrow 0, \hbox{ where }\Ker
\partial^{2,1}_{\fg_0} = \fg_1.
$$

Let us prove that
$$
H^{1,2}_{\fg_0} = 0\hbox{ if either } \fg_0 = \fcpe(n) \hbox { or }
\fg_0 = \fspe(n)\subplus \langle a\tau + bz\rangle,
$$
$\hbox { where } a=0,\hbox{ }b\not= 0\hbox { or } a\not= 0,\hbox{ }
b/a\not= n.$

Let $\fg_0 =\fspe(n)\subplus \langle\tau\rangle = \fpe(n).$ Since
$\fpe(n)$ and $\Pi
(E^2V^*)$ are isomorphic $\fg_0$-modules, the $\fpe(n)$-module
$\fg_0\otimes \fg_{-1}^*$ is isomorphic to
$$
\Pi (E^2V^*)\otimes V^* \cong E^2V^*\otimes V \cong \fg_{-1}\otimes
E^2\fg_{-1}^*
= \Ker \partial^{1,2}_{\fg_0}.
$$
By part a) of Theorem 1.2 we have $\fg_1 = 0$. Therefore,  by
(3.1.1)
$$
\Im \partial^{2,1}_{\fg_0} \cong \fg_0\otimes \fg_{-1}^* \cong
\Ker \partial^{1,2}_{\fg_0},
$$
i.e., $H^{1,2}_{\fpe(n)} = 0$.

Let $\fg_0 = \fspe(n)\subplus \langle a\tau + bz\rangle$, where $a = 0,\hbox{ }
b\not= 0$ or $a\not= 0,\hbox{ } b/a \not= n$. By Theorem 1.2,
$\fg_1 = 0$ for such $\fg_0$. Note that $\dim \fg_0 =  \dim
\fpe(n)$. Therefore,
$$
\dim  \Im \partial^{2,1}_{\fg_0} =  \dim  \fg_0\otimes \fg_{-1}^* =  \dim
\fpe(n)\otimes V^* =  \dim E^2V^*\otimes V =  \dim
\Ker\partial^{1,2}_{\fg_0}.
$$
Hence, $H^{1,2}_{\fg_0} = 0$.

Let $\fg_0 =\fcpe(n).$ Then by part b) of Theorem 1.2 $\fg_1 =
V^*$. Note that
$$
\fg_0\otimes \fg_{-1}^* = (\fpe(n)\oplus \langle z\rangle)\otimes V^*
= \fpe(n)\otimes V^*
\oplus \langle z\rangle\otimes V^* \cong V\otimes
E^2V^*\oplus V^*.
$$
Therefore,
$$
\Im \partial^{2,1}_{\fg_0} = V\otimes E^2V^* =  \Ker \partial^{1,2}_{\fg_0}.
$$
Hence, $H^{1,2}_{\fcpe(n)} = 0$.

Let us prove that
$$
H^{1,2}_{\fg_0} = V^*\hbox{ if } \fg_0 = \fspe(n)\hbox{ or } \fspe(n)\subplus
\langle \tau + nz\rangle.
$$
Let $\fg_0 = \fspe(n)$. By part a) of Theorem 1.2
$\fg_1 = 0$. Therefore,
$$
\Im \partial^{2,1}_{\fg_0} \cong \fg_0\otimes \fg_{-1}^* = \fspe(n)\otimes V^*.
$$
As has been shown for the case $\fg_0 =\fcspe(n) = \fspe(n)\oplus
\langle z\rangle,$ we have
$$
\fg_0\otimes \fg_{-1}^* = (\fspe(n)\oplus \langle z\rangle)\otimes
V^* \cong V\otimes E^2V^*.
$$
Therefore,
$$
\Ker \partial^{1,2}_{\fspe(n)}/\Im \partial^{2,1}_{\fspe(n)} \cong V^*.
$$
Hence $H^{1,2}_{\fspe(n)} = V^*$.

Finally, let $\fg_0
= \fspe(n)\subplus \langle \tau + nz\rangle$. By part b) of Theorem
1.2 $\fg_1 =
V^*$. Since
$$
\fg_0\otimes \fg_{-1}^* = (\fspe(n)\subplus \langle \tau +
nz\rangle)\otimes V^*,
$$
the  Jordan-H\"{o}lder series for the $\fspe(n)$-module $\Im
\partial^{2,1}_{\fg_0}$ contains the same irreducible quotient
modules as that for the $\fspe(n)$-module $\fspe(n) \otimes V^*$.
Since
$$
\Ker \partial^{1,2}_{\fg_0} = V\otimes E^2V^* \cong (\fspe(n)\oplus
\langle z\rangle)\otimes V^*,
$$
then $H^{1,2}_{\fg_0} = V^*$. This proves Theorem 2.1 in the case
of the SFs  of the first order.

\ssec{3.2. Hochschild-Serre spectral sequence.} The continuation of the proof
of  Theorem 2.1 is based on the Hochschild-Serre spectral
sequence. Let us recall the corresponding formulations in a form
convenient for us, since the case of Lie superalgebras is hardly
reflected in the literature (one might think that the union of
[Fu] and [GM] should suffice, but the sign rule
  applied to the Lie algebra case does not completely solve the problem).

Let $\fg = \fg_0\oplus \fg_1$ be a Lie superalgebra and $M$ be a
$\fg$-module. On the superspace of $k$-dimensional cochains $C^k =
C^k (\fg, M)$ define a filtration:
$$
F^0C^k = C^k\supset F^1C^k\supset F^2C^k\supset \ldots\supset
F^jC^k\supset \ldots
\supset F^{k+1}C^k = 0,
$$
where
$$
F^jC^k = \lbrace c\in C^k|\hbox{ } c(g_1,\ldots ,g_i,\ldots ,g_k)
  = 0 \hbox{ if }k-j+1 \hbox { arguments belong to } \fg_0\rbrace,
$$
$0\leq j\leq k+1.$
Set
$$
Z^{p,q}_r = \lbrace c\in F^pC^{p+q}|\hbox{ } dc \in F^{p+r}C^{p+q+1}\rbrace.
\eqno  (3.2.1)
$$
Finally, set
$$
E^{p,q}_r = Z^{p,q}_r/(Z^{p+1,q-1}_{r-1} + dZ^{p-r+1,q+r-2}_{r-1}).
\eqno (3.2.2)
$$
Notice that the differential $d$ induces the differentials
$$
d^{p,q}_r: E^{p,q}_r\longrightarrow E^{p+r,q-r+1}_r  \eqno(3.2.3)
$$
and $E^{p,q}_{r+1} = H^{p,q}(E_r)$ [GM].

Since $d(F^jC^k)\subset F^jC^{k+1}$, on $H^k = H^k(\fg, M)$ we get
the induced filtration such that $F^pH^k/F^{p+1}H^k =
E^{p,q}_\infty $, where $p+q=k$.

We want to compute the group $H^2(V, \hbox{ }\fg_*)$, where $V =
V_0\oplus V_1 $ is the standard $\fpe(n)$-module, $\fg_* =
\fg_*(V, \hbox{ }\fcpe(n)).$ The
  Hochschild-Serre
spectral sequence corresponding to the subalgebra $V_0$ converges
to $H^2(V,\hbox{ }\fg_*).$ Thus, $H^2(V,\hbox{ }\fg_*) = \oplus
_{p+q=2}E^{p,q}_\infty $ and in order to compute
  the limit
terms of the spectral sequence $E^{p,q}_\infty $ we have to consider three
  cases:

\emph{1)} $p = 2$, $q = 0.$ Then by formula (3.2.3) we have
$$
E^{1,0}_1\buildrel {d^{1,0}_1}\over \longrightarrow
E^{2,0}_1\buildrel {d^{2,0}_1}\over \longrightarrow
E^{3,0}_1,
$$
$$
E^{0,1}_2\buildrel {d^{0,1}_2}\over \longrightarrow
E^{2,0}_2\buildrel {d^{2,0}_2}\over \longrightarrow 0, \eqno (3.2.4)
$$
$$
0\buildrel {d^{-1,2}_3}\over \longrightarrow
E^{2,0}_3\buildrel {d^{2,0}_3}\over \longrightarrow 0.
$$
Therefore, $E^{2,0}_\infty = E^{2,0}_3.$

\emph{2)} $p=1$, $q=1.$ Then by formula (3.2.3) we have
$$
E^{0,1}_1\buildrel {d^{0,1}_1}\over \longrightarrow
E^{1,1}_1\buildrel {d^{1,1}_1}\over \longrightarrow
E^{2,1}_1,
$$
$$
0\buildrel {d^{-1,2}_2}\over \longrightarrow
E^{1,1}_2\buildrel {d^{1,1}_2}\over \longrightarrow
E^{3,0}_2, \eqno (3.2.5)
$$
$$
0\buildrel {d^{-2,3}_3}\over \longrightarrow
E^{1,1}_3\buildrel {d^{1,1}_3} \over \longrightarrow 0.
$$
Therefore, $E^{1,1}_\infty = E^{1,1}_3.$

\emph{3)} $p = 0, q = 2.$ Then by formula (3.2.3) we have
$$
0\buildrel {d^{-1,2}_1}\over \longrightarrow
E^{0,2}_1\buildrel {d^{0,2}_1}\over \longrightarrow
E^{1,2}_1,
$$
$$
0\buildrel {d^{-2,3}_2}\over \longrightarrow
E^{0,2}_2\buildrel {d^{0,2}_2}\over \longrightarrow
E^{2,1}_2, \eqno (3.2.6)
$$
$$
0\buildrel {d^{-3,4}_3}\over \longrightarrow
E^{0,2}_3\buildrel {d^{0,2}_3}\over \longrightarrow
E^{3,0}_3,
$$
$$
0\buildrel {d^{-4,5}_4}\over \longrightarrow
E^{0,2}_4\buildrel {d^{0,2}_4}\over \longrightarrow 0.
$$
Therefore, $E^{0,2}_\infty = E^{0,2}_4.$

\ssec{3.3. Continuation of the proof.} Notice that by [Fu] $E^{p,q}_1 =
H^q(\fg_0,\hbox{ } M\otimes S^p\fg_1^*)$. Since in our case $V =
V_0\oplus V_1$ is a commutative Lie superalgebra,  then
$$
E^{p,q}_1 = H^q(V_0,\hbox{ } \fg_*\otimes S^pV_1^*) = H^q(V_0,\hbox{
} \fg_*)\otimes S^pV_1^*.
\eqno (3.3.1)
$$
Let us calculate $H^q(V_0, \hbox{ }\fg_*)$ for $q
= 0, 1, 2.$ By Lemma 1.1
$$
\begin{array}{l}
   \fg_* = \fpe(n+1) = \fg_{-1}\oplus \fg_0\oplus \fg_1 \oplus
   \fg_2,\hbox{ where}\\
   \fg_{-1} = V_0\oplus V_1, \fg_0 =
   \fcpe(n) =
   \fsl(n)\oplus \Pi  (S^2V_0)\oplus \Pi (E^2V_1)\oplus \langle
   \tau\rangle\oplus \langle z\rangle,\\
   \fg_1 = \Pi (V_1)\oplus \Pi (V_0),\\
   \fg_2 = \Pi(\langle 1\rangle).
\end{array}
\eqno{(3.3.2}
$$

Recall that as an $\fsl(n+1)$-module,
$$
\fpe(n+1) \cong \fsl(n+1)\oplus E^2W_0^*\oplus S^2W_0\oplus \langle
\hbox{ d }\rangle ,
$$
where $\hbox{d}$ is $\diag(1_{n+1},\hbox{ } -1_{n+1})$ and $W_0$ is the
standard $\fsl(n+1)$-module. Clearly,
$$
\eqalign{
&E^2W_0^* =\Pi ( E^2V_1)\oplus V_1,\hbox{ }\cr &S^2W_0 =\Pi (
S^2V_0)\oplus \Pi (V_0)\oplus \fg_2,\hbox{ }\cr &\fsl(n+1) =
\fsl(n)\oplus \langle \tau + nz\rangle\oplus V_0\oplus \Pi (V_1),\cr &\hbox{ d
}= \tau -z.\cr }
$$
Let $\eps_1,\ldots, \eps_{n+1}$ be the standard basis of the dual
space to the space of the diagonal matrices in $\fgl(n+1)$ and the
ordering is performed
so that
$$
\Delta_+ = \lbrace \eps_i-\eps_j, i<j\rbrace ,\hbox{ }
\Delta_ - = \lbrace \eps_i-\eps_j, i>j\rbrace.
$$
Let  $E_{\eps_i-\eps_j} (i\not= j)$ be the corresponding root vectors.
Then $V_0$ is the subspace of $\fsl(n+1)$ generated by
$$
E_{\eps_1-\eps_{n+1}},\hbox{ } E_{\eps_2-\eps_{n+1}},\ldots ,
\hbox{ } E_{\eps_n-\eps_{n+1}}.
$$
Let $V_\lambda $ be the
irreducible $\fsl(n+1)$-module with highest weight $\lambda$. The
BWB theorem says [Kos] that there exists a 1-1 correspondence
between the irreducible components of $H^q(V_0,\hbox{ }
V_\lambda)$, considered as $\fgl(n)$-module,
  and elements
$w\in W(\fsl(n+1))$ of length $q$ from the Weyl group of the Lie
algebra $\fsl(n+1)$ such that
$$
w(\Delta_-)\cap \Delta_+\subset \lbrace \eps_1-\eps_{n+1},\hbox{ }
\eps_2-\eps_{n+1},\ldots,\eps_n-\eps_{n+1}\rbrace. \eqno (3.3.3)
$$
Moreover, the highest weight of the $\fgl(n)$-module
corresponding to $w$ is equal to $w(\lambda + \rho
)-\rho$, where $\rho =(\sum _{\alpha\in {\Delta_+}}\alpha)/2.$

Notice that $E^2W_0^*,\hbox{ } S^2W_0,\hbox{ } \fsl(n+1),\hbox{ }
\langle \tau-z\rangle$ are all irreducible $\fsl(n+1)$-modules with highest
weights, respectively,
$$
\eps_1 + \ldots +\eps_{n-1},\hbox{ } 2\eps_1,\hbox{ } 2\eps_1 + \eps_2 +
\ldots + \eps_n,\hbox{ } 0.
$$

Let us find the highest weights of irreducible
$\fgl(n)$-submodules of
$H^q(V_0,\hbox{ } V_\lambda )$ for each of the indicated $\lambda $.

\emph{1)} $q = 0$. The only element of the Weyl group of length 0 is
the unit. Hence
$w(\lambda +\rho )-\rho = \lambda$.

\emph{2)} $q = 1$. Let $ \langle \alpha_1,\ldots, \alpha_n\rangle$, where
$$
\alpha_1 = \eps_1-\eps_2,\ldots ,\hbox{ } \alpha_i = \eps_i-
\eps_{i+1},\ldots ,\hbox{ } \alpha_n = \eps_n-\eps_{n+1}
$$
be the
system of simple roots. The  elements of the Weyl group  of length
1 are reflections corresponding to the simple roots:
$$
r_{\alpha_i}:\alpha\longrightarrow \alpha -{2(\alpha_i, \alpha)\over
(\alpha_i, \alpha_i)}\alpha_i.
$$
Since the only element $r_{\alpha_i}$ satisfying (3.3.3) is $r_{\alpha_n}$,
then $w(\lambda + \rho)-\rho = r_{\alpha_n}(\lambda)-\alpha_n$. For
$$
\lambda = \eps_1 + \ldots + \eps_{n-1},\hbox{ } 2\eps_1,\hbox{ } 2\eps_1 +
\eps_2 + \ldots + \eps_n,\hbox{ } 0
$$
this expression is equal to, respectively:
$$
-2\eps_n,\hbox{ } \eps_1-\eps_2-\ldots -\eps_{n-1}-2\eps_n,\hbox{ }
  -\eps_2-\eps_3-\ldots -\eps_{n-1}-3\eps_n,\hbox{ }
-\eps_1-\eps_2-\ldots -\eps_{n-1}-2\eps_n.
$$

\emph{3)} $q = 2$. The
elements of length  2 are of the form $r_{\alpha_i}r_{\alpha_j}.$
The only such element satisfying (3.3.3) is
$r_{\alpha_n}r_{\alpha_{n-1}}.$ Then
$$
w(\lambda + \rho)-\rho = r_{\alpha_n}r_{\alpha_{n-1}}(\lambda)-\alpha_{n-1}
-2\alpha_n.
$$
For $\lambda = \eps_1 + \ldots + \eps_{n-1}$,
$2\eps_1$, $2\eps_1 + \eps_2 + \ldots + \eps_n$, $0$ this
expression is equal to, respectively:
$$
\eqalign{
&-2\eps_1-2\eps_2-\ldots -2\eps_{n-2}-4\eps_{n-1}-4\eps_n, \cr
&-3\eps_2-3\eps_3 \hbox{ }(\hbox{if } n=3) \hbox{ or }
-2\eps_2-\ldots -2\eps_{n-2}-3\eps_{n-1}-3\eps_n \hbox{ }
(\hbox{if } n>3),\cr &-\eps_1-3\eps_2-4\eps_3 \hbox{ }(\hbox{if }
n=3) \hbox{ or } -\eps_1-2\eps_2-\ldots
-2\eps_{n-2}-3\eps_{n-1}-4\eps_n \hbox{ } (\hbox{if } n>3),\cr
&-2\eps_1-\ldots -2\eps_{n-2}-3\eps_{n-1}-3\eps_n.\cr }
$$
\begin{Remark} We have obtained the weights with
respect to
  $\fgl(n) = \fsl(n)\oplus
\langle \tau+ nz\rangle$ embedded into $\fsl(n+1)$.  Now it is not difficult to
rewrite these weights as the highest ones with respect to
$\fpe(n)_0 = \fsl(n)\oplus \langle \tau\rangle$. We collect all our results in
the following
\end{Remark}

\ssbegin{3.4}{Lemma} $\fgl(n) = \fpe(n)_0 \hbox{ module }
H^q(V_0,\hbox{ } \fg_*)$ $(q = 0,\hbox{ } 1,\hbox{ } 2)$ is the
direct sum of irreducible submodules with the highest weights and
highest vectors  listed in Table 2.
\end{Lemma}

\ssbegin{3.5}{Lemma} Let $V_\lambda$ be an
irreducible  $\fgl(n)$-module with highest weight
$\lambda.$ Then $E^{p,0}_1 (p = 1,\hbox{ } 2,\hbox{ } 3)$,
$E^{p,1}_1 (p = 0,\hbox{ } 1,\hbox{ } 2)$, $E^{p,2}_1$ $(p = 0,\hbox{
} 1)$ are the direct sums  of irreducible $\fgl(n)$-submodules with
highest weights given in the corresponding columns of Tables 3, 4,
and 5, respectively.
\end{Lemma}

\begin{proof} By formula
(3.3.1) the following $\fgl(n)$-modules are isomorphic:\break %VM
$E^{p,q}_1\cong H^q(V_0,\hbox{ } \fg_*)\otimes
S^pV_0.$ Making use of the description of $H^q(V_0, \hbox{
}\fg_*)$ as a $\fgl(n)$-module given in Lemma 3.4 we find the
decomposition of the indicated tensor product into irreducible
components described in Tables 3, 4, and 5.
\end{proof}

\ssbegin{3.6}{Lemma} $E^{p,0}_0(p = 1,\hbox{ } 2,\hbox{ } 3)$ are the direct
sums of the irreducible $\fgl(n)$-modules with highest weights
described in the corresponding columns of Table 6.
\end{Lemma}
\begin{proof} By formula (3.2.2) we have
$$
E^{p,0}_0 = Z^{p,0}_0 = \fg_*\otimes S^pV_1^* = \fg_*\otimes S^pV_0.
$$
Making use of the description of $\fg_*$ as $\fgl(n)$-module
given in (3.3.2) we find the decomposition of the indicated tensor
products into direct sum of
the irreducible components described in Table 6.
\end{proof}

\ssbegin{3.7}{Lemma} $E^{2,0}_\infty $ is an
irreducible $\fgl(n)$-module with highest weight
$2\eps_1+2\eps_2.$
\end{Lemma}
\begin{proof} First, recall that $H^{1,2}_{\fcpe(n)} = 0$ by the already
  proved part of Theorem 2.1 for the case of SFs of order 1. Making
use of Lemma 3.4, we note that if $E^{2,0}_\infty $
had contained a $\fgl(n)$-submodule belonging to either
$V_{-\eps_n}\otimes S^2V_0$ or $V_{\eps_1}\otimes S^2V_0$, then
this submodule would have belonged to $H^{1,2}_{\fcpe(n)}.$
\end{proof}

Therefore, with the help of Table 3 we deduce that
$E^{2,0}_\infty $ has no $\fgl(n)$-submodules with highest
weights $2\eps_1-\eps_n$, $\eps_1$, $3\eps_1$, $2\eps_1+\eps_2$.

Let us show that $E^{2,0}_\infty $ has no irreducible
$\fgl(n)$-submodules with
highest weights $4\eps_1$, $3\eps_1+\eps_2$, and $2\eps_1$ either.
More precisely, let us show that even $E^{2,0}_2$ does not have them.

Recall that the corresponding differentials act as follows:
$$
E^{1,0}_1\buildrel {d^{1,0}_1}\over \longrightarrow
E^{2,0}_1\buildrel {d^{2,0}_1}\over \longrightarrow
E^{3,0}_1.
$$
Note that according to (3.2.2), $E^{p,0}_1 = Z^{p,0}_1$ for $p = 1, 2, 3.$
Let us show that $\Ker d^{2,0}_1$ has no components with weights
$4\eps_1$, $3\eps_1+\eps_2, \hbox{ and }2\eps_1.$
It follows from Tables 2 and 3 that the corresponding highest vectors in
$E^{2,0}_1$ are
$$
v_{4\eps_1} = e_1\tilde e_1\otimes \tilde e_1^2,\hbox{ }
v_{3\eps_1+ \eps_2} = e_1\tilde e_2\otimes \tilde e_1^2- e_1\tilde
e_1\otimes \tilde e_1\tilde e_2,\hbox{ } v_{2\eps_1} =
(\tau-z)\otimes \tilde e_1^2.
$$
We remind  the reader that the
differentials $d$ in our case  are the same as the differentials
$\partial^{k,s}_{\fg_0}$. Notice that if $c\in E^{2,0}_1$ then
$$
d^{2,0}_1c(f_1, f_1, f_1) = dc(f_1, f_1, f_1) = -3c(f_1,f_1)(f_1).
$$
Therefore,
$$
d^{2,0}_1v_{4\eps_1}(f_1, f_1) = 3e_1, \hbox{ }
d^{2,0}_1v_{3\eps_1+\eps_2}(f_1, f_1) = (3/2)e_2,\hbox{ }
d^{2,0}_1v_{2\eps_1}(f_1, f_1) = -6f_1. \eqno(3.7.1)
$$
Hence,
$$
v_{4\eps_1},\hbox{ } v_{3\eps_1+\eps_2},\hbox{ } v_{2\eps_1}\not\in
\Ker d^{2,0}_1.
$$

Finally, let us prove that in $E^{2,0}_\infty $ there is an irreducible
  $\fgl(n)$-submodule of  highest weight $2\eps_1+2\eps_2.$ Notice that
$E^{2,0}_2$ has submodule of  highest weight $2\eps_1+2\eps_2$,
since according to Table 3, this module is contained in
$E^{2,0}_1$ and is not contained in either $E^{1,0}_1$ or
$E^{3,0}_1.$  Recall that the corresponding
differentials act as follows:
$$
E^{0,1}_2\buildrel {d^{0,1}_2}\over \longrightarrow
E^{2,0}_2\buildrel {d^{2,0}_2}\over \longrightarrow 0.
$$
Therefore, $\Ker d^{2,0}_2 = E^{2,0}_2$ has a component of weight
$2\eps_1+ 2\eps_2.$ By Table 4 $E^{0,1}_1$  has no components of
weight $2\eps_1 + 2\eps_2$. Hence, neither $E^{0,1}_2$ nor $\Im
d^{0,1}_2$ have  such a component . Therefore, it must be in
$E^{2,0}_3 = E^{2,0}_\infty$.

\ssbegin{3.8}{Lemma} \emph{a)}
As a $\fgl(n)$-module, $E^{1,1}_\infty $ can only have the
irreducible submodules with the following highest weights, each of
multiplicity
  not greater then 1:
$$
2\eps_1,\hbox{ } 2\eps_1+\eps_2-\eps_n,\hbox{ }
\eps_1+\eps_2-2\eps_n, \hbox{ and } \eps_1-\eps_n;
$$

\emph{b)}
$E^{1,1}_\infty $ has an irreducible $\fgl(n)$-submodule with
highest weight $2\eps_1$.
\end{Lemma}
\begin{proof} By Theorem
2.1  for the case of SFs of order 1 and Tables 2 and 4 we see that
$E^{1,1}_\infty $ has no irreducible $\fgl(n)$-submodules of
highest weight $\eps_1-2\eps_n$ and $-\eps_n$, since they would
have corresponded to SFs of order 1.

Let us show that there are no components of weight
$2\eps_1-2\eps_n$ or 0 in $E^{1,1}_\infty $, more precisely, that
even $E^{1,1}_2$ does not have them. Recall that the
corresponding differentials act as follows:
$$
E^{0,1}_1\buildrel {d^{0,1}_1}\over \longrightarrow
E^{1,1}_1\buildrel {d^{1,1}_1}\over \longrightarrow
E^{2,1}_1.
$$
By (3.2.2) we have
$$
\begin{array}{l}
   E^{0,1}_1 = Z^{0,1}_1/(Z^{1,0}_0 + dZ^{0,0}_0),\\
   E^{1,1}_1 = Z^{1,1}_1/(Z^{2,0}_0 + dZ^{1,0}_0),\\
   E^{2,1}_1 = Z^{2,1}_1/(Z^{3,0}_0 + dZ^{2,0}_0).
\end{array}
\eqno{3.8.1}
$$

By Tables 2 and 4 the
highest vectors of weights $2\eps_1-2\eps_n$ and 0 in $E^{1,1}_1$
are, respectively,
$$
v_{2\eps_1-2\eps_n} = (e_1\wedge \tilde f_n)\otimes \tilde e_1\wedge
\tilde f_n \hbox{ and } v_0 = \sum _{i=1}^n(\tau -z)\otimes \tilde
f_i\wedge\tilde e_i.
$$
We see that
$$
dv_{2\eps_1-2\eps_n}(e_n, f_1, f_1) = -v_{2\eps_1-2\eps_n}
(f_1, f_1)(e_n) -2v_{2\eps_1-2\eps_n}(e_n, f_1)(f_1) = -f_n/2\not=
0,
$$
$$
dv_0(e_1, f_1, f_1) = -v_0(f_1, f_1)(e_1)-2v_0(e_1, f_1)(f_1) = 2f_1\not= 0.
$$
Suppose that
$$
dv_{2\eps_1-2\eps_n}\in Z^{3,0}_0 + dZ^{2,0}_0.
$$
Then there exist highest $\fgl(n)$-vectors
$v_{2\eps_1-2\eps_n}'\in Z^{3,0}_0$ and $v_{2\eps_1-2\eps_n}''\in
dZ^{2,0}_0$ of weight $2\eps_1-2\eps_n$ such that
$$
dv_{2\eps_1-2\eps_n} = v_{2\eps_1-2\eps_n}' +
v_{2\eps_1-2\eps_n}''.
$$
Since $e_n\in V_0$, then
$v_{2\eps_1-2\eps_n}'(e_n, f_1, f_1) = 0.$ Hence
$v_{2\eps_1-2\eps_n}'' \not= 0$ and therefore, $dZ^{2,0}_0$ has an
irreducible $\fgl(n)$-submodule of highest weight
$2\eps_1-2\eps_n.$

Similarly, having assumed that
$dv_0\in Z^{3,0}_0 + dZ^{2,0}_0$, we deduce that $dZ^{2,0}_0$ has
an irreducible $\fgl(n)$-submodule of weight 0. Note that
according to Table 6, $E^{2,0}_0$ has no submodules of highest
weight $2\eps_1-2\eps_n$ or 0. Since  $E^{2,0}_0 = Z^{2,0}_0$,
then $Z^{2,0}_0$ and $dZ^{2,0}_0$ have no such components either.
Therefore, $dv_{2\eps_1-2\eps_n}$ and $dv_0$ do not belong to
$Z^{3,0}_0 + dZ^{2,0}_0.$ Thanks to (3.8.1) this implies that
$$
d^{1,1}_1v_{2\eps_1-2\eps_n}\not= 0 \hbox{ and }
d^{1,1}_1v_0\not= 0.
$$
Hence, $\Ker  d^{1,1}_1$, and therefore,
$E^{1,1}_2$, have no irreducible $\fgl(n)$-submodules of highest
weight $2\eps_1-2\eps_n$ and 0.

Let us prove now that $E^{1,1}_\infty  = E^{1,1}_3$ has no
irreducible $\fgl(n)$-submodule of highest weight
$3\eps_1-\eps_n$.
Notice that $E^{1,1}_2$
has such a submodule, since by Table 4 it is contained in
$E^{1,1}_1$ and is not contained in either $E^{0,1}_1$ or
$E^{2,1}_1$. Tables 2 and 4 imply that the $\fgl(n)$-highest
vector in $E^{1,1}_1$ of weight $3\eps_1-\eps_n$ is
$v_{3\eps_1-\eps_n} = e_1\tilde e_1\otimes \tilde f_n\wedge\tilde
e_1.$

Recall that the corresponding differentials act
as follows:
$$
0\buildrel {d^{-1,2}_2}\over \longrightarrow
E^{1,1}_2\buildrel {d^{1,1}_2}\over \longrightarrow
E^{3,0}_2.
$$
By formula (3.2.2) we have
$$
\begin{array}{l}
   E^{1,1}_2 = Z^{1,1}_2/(Z^{2,0}_1 + dZ^{0,1}_1),\\
   E^{3,0}_2 = Z^{3,0}_2/dZ^{2,0}_1.
\end{array}
\eqno{(3.8.2)}
$$

Thanks to formulas (3.8.1)
we see that the $\fgl(n)$-highest vector in $E^{1,1}_2$
  of weight $3\eps_1-\eps_n$ is
$$
w_{3\eps_1-\eps_n} = v_{3\eps_1-\eps_n} +
v_{3\eps_1-\eps_n}' + v_{3\eps_1-\eps_n}'',
$$
where
$v_{3\eps_1-\eps_n}'$ and $v_{3\eps_1-\eps_n}''$ are
$\fgl(n)$-highest vectors in $Z^{2,0}_0$ and $dZ^{1,0}_0,$
respectively.

Since by (3.2.1)
$dw_{3\eps_1-\eps_n}\in Z^{3,0}_2$, then $dw_{3\eps_1-\eps_n}(e_n,
f_1, f_1) = 0.$ Since $ v_{3\eps_1-\eps_n}''\in dZ^{1,0}_0$, then
$dv_{3\eps_1-\eps_n}'' = 0.$ We have
$$
dv_{3\eps_1-\eps_n}(e_n, f_1, f_1) = -2v_{3\eps_1-\eps_n}
(e_n, f_1)(f_1) = -e_1\not= 0.
$$
Therefore,
$v_{3\eps_1-\eps_n}'\not= 0.$ Looking at Table 6 we see that the
unique highest vector of weight $3\eps_1-\eps_n$ in $Z^{2,0}_0$ is
$e_1\wedge \tilde f_n\otimes \tilde e_1^2.$ Hence,
$v_{3\eps_1-\eps_n}' = ke_1\wedge \tilde f_n\otimes \tilde e_1^2$,
where $k\in \Cee^*.$ Note that since $v_{3\eps_1-\eps_n}(f_1, f_1)
= 0$, then $dv_{3\eps_1-\eps_n}(f_1, f_1, f_1) = 0$. Therefore,
$$
\eqalign{
&dw_{3\eps_1-\eps_n}(f_1, f_1, f_1) = dv_{3\eps_1-\eps_n}' (f_1,
f_1, f_1) = \cr &-3v_{3\eps_1-\eps_n}'(f_1, f_1)(f_1) =
(-3/2)kf_n\not= 0.\cr }
$$
Note that $dw_{3\eps_1-\eps_n}\not\in
dZ^{2,0}_1$. In fact, by Table 3 $E^{2,0}_1$ has no  irreducible
$\fgl(n)$-component with highest weight $3\eps_1-\eps_n$. Since
$E^{2,0}_1 = Z^{2,0}_1$, then $Z^{2,0}_1$
  and $dZ^{2,0}_1$ have no such  component either.
Therefore, by (3.8.2) we have $d^{1,1}_2w_{3\eps_1-\eps_n}\not=
0.$ Hence, $\Ker d^{1,1}_2 = E^{1,1}_3$ has no components with
highest weight $3\eps_1-\eps_n.$

Let us prove that the irreducible component with highest weight
$\eps_1-\eps_n$ cannot be contained in $E^{1,1}_\infty $ with
multiplicity greater than 1. Note that $E^{1,1}_1$ has two
components of weight
  $\eps_1-\eps_n$.  According to Tables 2 and 4, one of the
$\fgl(n)$-highest vectors of weight $\eps_1-\eps_n$ in
$E^{1,1}_1$ is $v_{\eps_1-\eps_n} = (\tau -z)\otimes \tilde
f_n\wedge\tilde e_1.$ We see that
$$
dv_{\eps_1-\eps_n}(e_n, f_1, f_1) = -2v_{\eps_1-\eps_n}
(e_n, f_1)(f_1) = -(\tau -z)(f_1) = 2f_1\not= 0.
$$
Suppose that
$dv_{\eps_1-\eps_n}\in Z^{3,0}_0 + dZ^{2,0}_0.$ Then there exist
$\fgl(n)$-highest vectors $v_{\eps_1-\eps_n}'\in
Z^{3,0}_0$
  and $v_{\eps_1-\eps_n}''\in dZ^{2,0}_0$ of weight
  $\eps_1-\eps_n$ such that
$$
dv_{\eps_1-\eps_n} = v_{\eps_1-\eps_n}' +
v_{\eps_1-\eps_n}''.
$$
Since $e_n\in V_0$, then
$v_{\eps_1-\eps_n}' (e_n, f_1, f_1) = 0.$ Therefore,
$v_{\eps_1-\eps_n}''\not= 0.$ Now note that
$$
dv_{\eps_1-\eps_n}(e_1, f_1, f_1) =
-v_{\eps_1-\eps_n}(f_1, f_1)(e_1) -2v_{\eps_1-\eps_n}(e_1,
f_1)(f_1) = 0.
$$
Since $e_1\in V_0$, then $v_{\eps_1-\eps_n}'(e_1,
f_1, f_1) = 0.$ Therefore,
$$
v_{\eps_1-\eps_n}''(e_1, f_1, f_1) = 0. \eqno (3.8.3)
$$
By Table 6 $E^{2,0}_0 = Z^{2,0}_0$ contains a unique highest
vector of weight $\eps_1-\eps_n$, namely,
$$
\sum_{i=1}^n f_n\wedge \tilde f_i\otimes \tilde e_i\tilde e_1.
$$
Then
$$
v_{\eps_1-\eps_n}'' = kd(\sum_{i=1}^n f_n\wedge \tilde f_i
\otimes \tilde e_i\tilde e_1), \hbox{ where }k\in \Cee^*.
$$
  Note that in this case
$$
\eqalign {
&v_{\eps_1-\eps_n}''(e_1, f_1, f_1) = kd(\sum_{i=1}^n f_n\wedge
\tilde f_i\otimes \tilde e_i\tilde e_1) (e_1, f_1, f_1) =\cr
&-k(\sum_{i=1}^n f_n\wedge \tilde f_i\otimes \tilde e_i\tilde e_1)
(f_1, f_1)(e_1) = (k/2)f_n \not= 0,\cr}
$$
which contradicts
(3.8.3). Thus, $dv_{3\eps_1-\eps_n}\not\in Z^{3,0}_0
  + dZ^{2,0}_0.$
Then (3.8.1) yields $d^{1,1}_1v_{\eps_1-\eps_n} \not= 0.$
Therefore, the component of highest weight $\eps_1-\eps_n$ can not
be contained in $\Ker d^{1,1}_1 $ and hence, in $E^{1,1}_\infty$,
with multiplicity exceeding 1. So part a) of Lemma 3.8 is proved.

Let us prove that $E^{1,1}_\infty = E^{1,1}_3$ contains an irreducible
  $\fgl(n)$-submodule of highest weight $2\eps_1$.
Note that $E^{1,1}_2$ does contain such a submodule, since by
Table 4 it is contained in $E^{1,1}_1$ and is not contained in
either $E^{0,1}_1$ or $E^{2,1}_1$. \hfill\break Let $u_{2\eps_1}$
be the $\fgl(n)$-highest vector of weight
  $2\eps_1$ in $E^{1,1}_1$. By  (3.8.1) $u_{2\eps_1}$ can be
  chosen so that
$$
u_{2\eps_1}(v_1, v_2) = 0  \hbox{ for any } v_1, v_2\in V_1.
$$
According to (3.8.1), the $\fgl(n)$-highest vector of weight
$2\eps_1$ in $E^{1,1}_2$ is
$$
w_{2\eps_1} = u_{2\eps_1} + t_{2\eps_1} + s_{2\eps_1},
\eqno (3.8.4)
$$
where $t_{2\eps_1}\in Z^{2,0}_0$ and $s_{2\eps_1
}\in dZ^{1,0}_0$ are $\fgl(n)$-highest vectors. If $dw_{2\eps_1}
= 0$, then
$$
w_{2\eps_1}\in \Ker d^{1,1}_2 = E^{1,1}_3
$$
and therefore, $E^{1,1}_\infty $ has an irreducible
$\fgl(n)$-submodule with highest weight $2\eps_1$.

Suppose that $dw_{2\eps_1} \not= 0.$ Let us prove that then
$$
dw_{2\eps_1}\in dZ^{2,0}_1.   \eqno (3.8.5)
$$
Recall that $Z^{2,0}_1 = E^{2,0}_1$ and by Table 3 $E^{2,0}_1$ has
one highest vector of weight $2\eps_1$, namely, $v_{2\eps_1}.$ Let
us show that
$$
dw_{2\eps_1} = kdv_{2\eps_1}, \hbox { where } k\in \Cee^*.
\eqno (3.8.6)
$$
Note that by (3.2.1) $dw_{2\eps_1}\in Z^{3,0}_2$
and $dZ^{2,0}_1\subset Z^{3,0}_0.$ Therefore,  in order to prove
(3.8.6) it suffices to show that
$$
dw_{2\eps_1}(v_1, v_2, v_3) = kdv_{2\eps_1}(v_1, v_2, v_3),
\hbox{ where } k\in \Cee^*, \hbox{ for any } v_1, v_2, v_3\in
V_1.
$$
We have
$$
du_{2\eps_1}(v_1, v_2, v_3) = 0 \hbox{ for any } v_1, v_2, v_3\in V_1.
$$
Since $s_{2\eps_1}\in dZ^{1,0}_0$, then $ds_{2\eps_1} = 0$.
Therefore,
$$
dw_{2\eps_1}(v_1, v_2, v_3) = dt_{2\eps_1}(v_1, v_2, v_3)
\hbox { for any } v_1, v_2, v_3\in V_1. \eqno (3.8.7)
$$
Since $t_{2\eps_1}\in Z^{2,0}_0$, then $dt_{2\eps_1}\in Z^{2,1}_0.$
Hence $dt_{2\eps_1 }= t_{2\eps_1}' + t_{2\eps_1}''$,where
$t_{2\eps_1}', t_{2\eps_1}''$ are $\fgl(n)$-highest vectors from
$Z^{2,1}_0$ such that
$$
t_{2\eps_1}'(v_1, v_2, v_3) = 0 \hbox { for all } v_1, v_2, v_3\in V_1
\hbox{ and } t_{2\eps_1}''\in Z^{3,0}_0.
$$
Since by hypothesis
$dw_{2\eps_1} \not= 0$, then $t_{2\eps_1}''\not= 0.$ Since
$v_{2\eps_1}\in Z^{2,0}_1$, then $dv_{2\eps_1}\in Z^{3,0}_0.$

In Lemma 3.7  we have proved that $dv_{2\eps_1}\not=
0$ (see (3.7.1)). By Table 6 $Z^{3,0}_0 = E^{3,0}_0$ has only one
irreducible $\fgl(n)$-submodule with highest weight $2\eps_1$.
Hence $t_{2\eps_1}'' = kdv_{2\eps_1}$, where $k\in \Cee^*.$
Therefore,
$$
dt_{2\eps_1}(v_1, v_2, v_3) = kdv_{2\eps_1}(v_1, v_2, v_3)
\hbox { for any } v_1, v_2, v_3\in V_1.
$$
Thus, by (3.8.7)
$$
dw_{2\eps_1}(v_1, v_2, v_3) = kdv_{2\eps_1}(v_1, v_2, v_3)
\hbox { for any } v_1, v_2, v_3\in V_1
$$
and formula (3.8.6) is
proved. Then by (3.8.2) $d^{1,1}_2w_{2\eps_1} = 0.$ Therefore,
$w_{2\eps_1} \in \Ker d^{1,1}_2 = E^{1,1}_3.$ Thus, $E^{1,1}_\infty
$ contains an irreducible $\fgl(n)$-submodule with highest weight
$2\eps_1.$ This proves part b) of Lemma 3.8.
\end{proof}

\ssbegin{3.9}{Lemma} Only the following highest weights of irreducible
$\fgl(n)$-submodules can be encountered among those in
$E^{0,2}_\infty $:
$$
-2\eps_{n-1}-2\eps_n,\hbox{ } 2\eps_1-\eps_{n-1}-\eps_n,\hbox{ }
\eps_1-\eps_{n-1}-2\eps_n,\hbox{ } -\eps_{n-1}-\eps_n.
$$
\end{Lemma}
\begin{proof} By Table 5 $E^{0,2}_1$ is a direct sum
of irreducible $\fgl(n)$-components with the indicated highest
weights.

Note that due to Table 2 $H^2(V,\hbox{ } \fg_*)$ can only possess
first and second order SFs. By the statement of Theorem 2.1 for
SFs of order 1, $H^{1,2}_{\fcpe(n)} = 0$. Therefore, by Lemmas 7,
8, and 9 $H^{2,2}_{\fcpe(n)}$ can only contain irreducible
$\fgl(n)$-submodules with the following highest weights:
$$
\eqalign
{&2\eps_1 + 2\eps_2,\hbox{ } 2\eps_1, \hbox{ }2\eps_1 + \eps_2 -
\eps_n,\hbox{ } \eps_1 + \eps_2 - 2\eps_n,\hbox{ }
\eps_1-\eps_n,\cr &-2\eps_{n-1}-2\eps_n,\hbox{ } 2\eps_1 -
\eps_{n-1} - \eps_n,\hbox{ } \eps_1 - \eps_{n-1} - 2\eps_n,\hbox{
} -\eps_{n-1} - \eps_n\cr }
$$
each with multiplicity not greater
than one, and the components with highest weights $2\eps_1 +
2\eps_2$ and $2\eps_1$ are contained in $H^{2,2}_{\fcpe(n)}$ with
multiplicity one each.

Recall that by definition
$H^{2,2}_{\fg_0} $ is determined by the sequence
$$
\fg_1\otimes \fg_{-1}^*\buildrel {\partial^{3,1}_{\fg_0}}\over \longrightarrow
\fg_0\otimes E^2\fg_{-1}^*\buildrel {\partial^{2,2}_{\fg_0}}\over
\longrightarrow \fg_{-1}\otimes E^3\fg_{-1}^*.
$$
Thus, $\fg_0\otimes E^2\fg_{-1}^*$ for
$$
\fg_0 = \fcpe(n) \hbox{ or  }\fspe(n)\subplus \langle \tau + nz\rangle
$$
is equal to
$$
\fcpe(n)\otimes E^2\fg_{-1}^* \hbox{ or } (\fspe(n)\subplus \langle
\tau + nz\rangle)\otimes E^2\fg_{-1}^*,
$$
respectively. Note that $\Im \partial^{3,1}_{\fcpe(n)}$ and $\Im
\partial^{3,1}_{ \fspe(n)\subplus \langle \tau + nz\rangle}$ are isomorphic
$\fspe(n)$-modules by part b) of Theorem 1.2 and (3.1.1).

Note also that by Lemma 1.4 the Jordan-H\"{o}lder series for
$E^2V^*$ contains $\fspe(n)$-modules with highest weights
$2\eps_1$ and 0. Therefore, the Jordan-H\"{o}lder series of the
$\fspe(n)$-module $H^{2,2}_{\fcpe(n)}$, as compared with that of
$H^{2,2}_{\fspe(n)\subplus \langle \tau + nz\rangle}$, can
additionally contain only
the irreducible $\fspe(n)$-modules with highest weights $2\eps_1$
and 0. But we have shown that $H^{2,2}_{\fcpe(n)}$ has no trivial
$\fgl(n)$-submodule and therefore, $H^{2,2}_{\fcpe(n)}$ and
$H^{2,2}_{\fspe(n)\subplus \langle \tau + nz\rangle}$ can only differ by an
irreducible $\fspe(n)$-submodule with highest weight $2\eps_1$,
which,
  being considered as $\fsl(n)$-module, is the sum of irreducible
$\fsl(n)$-submodules with highest weights
$$
2\eps_1,\hbox{ } \eps_1-\eps_n,\hbox{ and } -\eps_{n-1}-\eps_n\hbox{
} (\hbox{see}\hbox{ } (1.4.4)).
$$
\end{proof}

\ssbegin{3.10}{Lemma} $H^{2,2}_{\fspe(n)\subplus \langle \tau +
nz\rangle}$ has no
irreducible $\fsl(n)$-submodule with highest weight $2\eps_1.$
\end{Lemma}

\begin{proof} Making use of Lemma 3.8, let us prove
that $E^{1,1}_\infty $ considered as an $\fsl(n)$-module has no
irreducible component with highest weight $2\eps_1$ in the case
where $\fg_0 = \fspe(n)\subplus \langle \tau + nz\rangle$.

Indeed, in
Lemma 3.8 we have shown that if $w_{2\eps_1}$ is the
$\fsl(n)$-highest vector of weight $2\eps_1$ in $E^{1,1}_2$ such
that $dw_{2\eps_1}\not= 0$ then $dw_{2\eps_1}\in dZ^{2,0}_1$ (see
(3.8.5)).

According to Table 2, $H^0(V_0,\hbox{}\fg_*(\fg_{-1}, \fg_0))$
has no irreducible $\fsl(n)$-submodule
with highest weight 0 when $\fg_0 = \fspe(n)\subplus \langle \tau + nz\rangle$
Therefore, by Table 3 $E^{2,0}_1$ has no irreducible
$\fsl(n)$-submodule with highest weight $2\eps_1$. Since
$Z^{2,0}_1 = E^{2,0}_1$, then $Z^{2,0}_1$ and $dZ^{2,0}_1$ have no
such  component either. Thus, by (3.8.2)
$d^{1,1}_2w_{2\eps_1} \not= 0$ and therefore, $w_{2\eps_1}\not\in
E^{1,1}_3 = E^{1,1}_\infty .$ It remains to show that
$dw_{2\eps_1}\not= 0.$ Recall that $w_{2\eps_1} = u_{2\eps_1} +
t_{2\eps_1} + s_{2\eps_1}$, where $u_{2\eps_1}\in E^{1,1}_1,
t_{2\eps_1}\in Z^{2,0}_0$ and $s_{2\eps_1}\in dZ^{1,0}_0$ are
$\fsl(n)$-highest vectors (see (3.8.4)). By Tables 2 and 4
$$
u_{2\eps_1} = 2\sum_{i=1}^n e_1\tilde e_i\otimes \tilde f_i\wedge
\tilde e_1-
(n+1)\sum_{i=1}^n e_1\tilde e_1 \otimes \tilde f_i\wedge \tilde
e_i.\eqno(3.10.1)
$$
Therefore,
$$
\eqalign{
&du_{2\eps_1}(e_1, f_1, f_1) = -u_{2\eps_1}(f_1, f_1)(e_1)
-2u_{2\eps_1}(e_1, f_1)(f_1) = \cr &-2(e_1\tilde
e_1-(n+1)e_1\tilde e_1/2)(f_1) = (n-1)e_1\not= 0.\cr }
$$
Since
$s_{2\eps_1}\in dZ^{1,0}_0$, then $ds_{2\eps_1} = 0.$ Hence, if
$dw_{2\eps_1} = 0$ then $t_{2\eps_1} \not= 0$ and since
$u_{2\eps_1}$ is an even vector, then vector $t_{2\eps_1}$ must be
even.

By Table 6 $Z^{2,0}_0 = E^{2,0}_0$ has 4 irreducible
$\fsl(n)$-submodules with highest weight $2\eps_1$. The corresponding
highest vectors are
$$
\eqalign{
&(\sum_{i=2}^n e_i\wedge \tilde f_i - (n-1)e_1\wedge \tilde f_1)
\otimes\tilde e_1^2- n\sum_{i=2}^ne_1\wedge \tilde f_i\otimes
\tilde  e_1\tilde e_i,\cr &\tau \otimes \tilde e_1^2, z\otimes
\tilde e_1^2,\hbox{ and } \fg_2\otimes \tilde e_1^2. \cr }.
$$
Only
the first three of these vectors are even. Therefore, if we
confine ourselves to the case $\fg_0 = \fspe(n)\subplus\langle \tau +
nz\rangle$, we
see that there should be two $\fsl(n)$-highest vectors of weight
$2\eps_1$ in $E^{2,0}_0.$ Let
$$
\eqalign {
&t_{2\eps_1} = k_1((\sum_{i=2}^n e_i\wedge \tilde f_i-
(n-1)e_1\wedge \tilde f_1)\otimes \tilde e_1^2
-n\sum_{i=2}^ne_1\wedge \tilde f_i\otimes \tilde e_1\tilde e_i)\cr
&+ k_2((\tau + nz)\otimes \tilde e_1^2),\hbox{ where } k_1,\hbox{
} k_2\in \Cee,\cr }
$$
be a linear combination of these vectors.
Note that
$$
\eqalign {
&du_{2\eps_1}(e_2, f_2, f_1) = -u_{2\eps_1}(f_2, f_1)(e_2) -
u_{2\eps_1}(e_2, f_1)(f_2) -u_{2\eps_1}(e_2, f_2)(f_1) = \cr
&-2(e_1\tilde e_2/2)(f_2) +  (n+1)(e_1\tilde e_1/2)(f_1) = -e_1/2
+ (n+1)e_1/2 = ne_1/2,\cr }
$$
$$
\eqalign {
&dt_{2\eps_1}(e_2, f_2, f_1) = -t_{2\eps_1}(f_2, f_1)(e_2) -
t_{2\eps_1}(e_2, f_1)(f_2) -t_{2\eps_1}(e_2, f_2)(f_1) = \cr
&(-1/2)k_1ne_1\wedge \tilde f_2(e_2) = (-1/4)k_1ne_1.\cr }
$$
Therefore, if $dw_{2\eps_1} = 0$, then $k_1 = 2.$ Observe that
$$
\eqalign {
&dt_{2\eps_1}(e_1, f_1, f_1) = -t_{2\eps_1}(f_1, f_1)(e_1)
-2t_{2\eps_1}(e_1, f_1)(f_1) = \cr &-k_1(n-1)(e_1\wedge \tilde
f_1)(e_1) + k_2(\tau +nz)(e_1) = (-k_1(n-1)/2 + k_2(n+1))e_1.\cr}
$$
Since $du_{2\eps_1}(e_1, f_1, f_1) = (n-1)e_1$, then
$$
-(n-1) + k_2(n+1) + (n-1) = 0.
$$
Hence, $k_2 = 0.$ But then
$$
\eqalign {
&dw_{2\eps_1}(f_1, f_1, f_1) = dt_{2\eps_1}(f_1, f_1, f_1) =
-3t_{2\eps_1}(f_1, f_1)(f_1) = \cr &-3k_1(n-1)e_1\wedge \tilde
f_1(f_1) = (3/2)k_1(n-1)f_1 = 3(n-1)f_1 \not= 0.\cr }
$$
This proves Lemma 3.10.
\end{proof}

Lemma 3.10 implies that $H^{2,2}_{\fcpe(n)}$ and
$H^{2,2}_{\fspe(n)\subplus \langle \tau +nz\rangle}$ differ by an irreducible
$\fspe(n)$-module with highest weight $2\eps_1$. Thus,
$H^{2,2}_{\fspe(n)\subplus \langle \tau+ nz\rangle}$ can only contain
irreducible
$\fsl(n)$-submodules with highest weights
$$
\eqalign {
&2\eps_1 + 2\eps_2,\hbox{ } 2\eps_1 + \eps_2 - \eps_n,\hbox{ }
\eps_1 + \eps_2 -2\eps_n,\cr &-2\eps_{n-1}-2\eps_n,\hbox{ }
2\eps_1-\eps_{n-1}-\eps_n,\hbox{ } \eps_1-\eps_{n-1}-2\eps_n\cr}
$$
each with multiplicity not greater than 1, and the
multiplicity of the submodule with highest weight $2\eps_1 +
2\eps_2$ is precisely 1.

\ssbegin{3.11}{Lemma} The
irreducible $\fpe(n)$-module with highest weight $2\eps_1 +
2\eps_2$ is the direct sum of irreducible $\fgl(n)$-modules with
the following highest weights:

\emph{a)} for $n>3$:
$$
\eqalign {
&2\eps_1 + 2\eps_2,\hbox{ } 2\eps_1 + \eps_2 - \eps_n,\hbox{ }
\eps_1 + \eps_2 -2\eps_n,\cr &-2\eps_{n-1}-2\eps_n,\hbox{ }
2\eps_1 - \eps_{n-1}-\eps_n,\hbox{ } \eps_1-\eps_{n-1}-2\eps_n;\cr
}
$$

\emph{b)} for $n=3$:
$$
2\eps_1 + 2\eps_2,\hbox{ } 2\eps_1 + \eps_2-\eps_3,\hbox{ } \eps_1 +
\eps_2 -2\eps_3.
$$
\end{Lemma}
\begin{proof} Let us consider the $\fpe(n)$-module $S^2(S^2V).$
Note that $v_{2\eps_1 + 2\eps_2} =
(e_1e_2)^2-(e_1^2)(e_2^2)$ is a $\fpe(n)$-highest vector. Indeed,
$B_{i,j}v_{2\eps_1 + 2\eps_2} = 0  \hbox{ for any  }i \hbox{ and }
j.$ Set
$$
\eqalign {
&v_{2\eps_1 + \eps_2 -\eps_n} =
(e_1e_2)(e_1f_n)-(e_1^2)(e_2f_n),\cr &v_{\eps_1 + \eps_2 -2\eps_n}
= (e_1f_n)\wedge (e_2f_n),\cr &v_{-2\eps_{n-1} - 2\eps_n} =
(f_{n-1}\wedge f_n)(f_{n-1}\wedge f_n),\cr
&v_{2\eps_1-\eps_{n-1}-\eps_n} = (e_1f_{n-1})\wedge
(e_1f_n)-(e_1^2)(f_{n-1}\wedge f_n),\cr
&v_{\eps_1-\eps_{n-1}-2\eps_n} = (f_{n-1}\wedge f_n)(e_1f_n).\cr}
$$
Notice that these vectors are the $\fgl(n)$-highest ones.
Moreover,
$$
\begin{array}{lr}
   C_{2,n}v_{2\eps_1 + 2\eps_2} = 2v_{2\eps_1 + \eps_2 - \eps_n},&(3.11.1)\\
   B_{2,n}v_{2\eps_1 + \eps_2 - \eps_n} = v_{2\eps_1 + 2\eps_2}, &(3.11.2)\\
   C_{1,n}v_{2\eps_1 + \eps_2 - \eps_n} = -3v_{\eps_1 + \eps_2
-2\eps_n}, &(3.11.3)\\
   B_{1,n}v_{\eps_1 + \eps_2 -2\eps_n} = -v_{2\eps_1 + \eps_2
-\eps_n}. &(3.11.4)
\end{array}
$$

If $n>3$, then additionally
$$
\eqalign {
&C_{2,n-1}v_{2\eps_1 + \eps_2 - \eps_n}= v_{2\eps_1 - \eps_{n-1} -
\eps_n},\cr &B_{2, n-1}v_{2\eps_1 - \eps_{n-1} - \eps_n} =
v_{2\eps_1 + \eps_2 -\eps_n},\cr &C_{1,n}v_{2\eps_1 - \eps_{n-1} -
\eps_n} = -3v_{\eps_1 - \eps_{n-1} -2\eps_n},\cr
&B_{1,n}v_{\eps_1-\eps_{n-1} - 2\eps_n} = -v_{2\eps_1 - \eps_{n-1}
- \eps_n},\cr &C_{1,n-1}v_{\eps_1 - \eps_{n-1} - 2\eps_n} =
v_{-2\eps_{n-1} -2\eps_n},\cr &B_{1,n-1}v_{-2\eps_{n-1} -2\eps_n}
= 2v_{\eps_1 - \eps_{n-1} -2\eps_n}.\cr }
$$
Therefore, if $n>3$,
then the irreducible $\fpe(n)$-module with highest weight $2\eps_1
+ 2\eps_2$ contains irreducible $\fgl(n)$-modules with highest
weights
$$
\eqalign {
&2\eps_1 + 2\eps_2,\hbox{ } 2\eps_1 + \eps_2 - \eps_n,\hbox{ }
\eps_1 + \eps_2 -2\eps_n,\cr &-2\eps_{n-1} -2\eps_n,\hbox{ }
2\eps_1 - \eps_{n-1} -\eps_n,\hbox{ } \eps_1 - \eps_{n-1}
-2\eps_n.\cr }
$$
We have already shown that the $\fspe(n)$-module
$H^{2,2}_{\fspe(n)\subplus \langle \tau + nz\rangle}$ does contain irreducible
$\fsl(n)$-submodules with these highest
  weights exactly, their multiplicities are not greater than 1,
and the multiplicity of $\fsl(n)$-submodule with highest weight
$2\eps_1 + 2\eps_2$ is  precisely one.

 From Tables 2 and 3 we see that the corresponding
$\fsl(n)$-highest vector is
$$
v_{2\eps_1 + 2\eps_2} =
(e_1\tilde e_1)\otimes  (\tilde e_2\tilde e_2) + (e_2\tilde
e_2)\otimes  (\tilde e_1\tilde e_1) - 2(e_1\tilde e_2)\otimes
(\tilde e_1\tilde e_2). \eqno(3.11.5)
$$
Hence, $v_{2\eps_1 + 2\eps_2}$ is an odd $\fspe(n)$-highest
vector. So part a) of Lemma 3.11  and Theorem 2.1 for the case
where $\fg_0 =  \fspe(n)\subplus \langle \tau + nz\rangle, n > 3$, are proved.

Let $n=3$. Consider the $\fspe(3)$-module $E^3V$. As an
$\fsl(3)$-module, this module is isomorphic to
$$
E^3(V_0\oplus V_0^*) = E^3V_0 \oplus (E^2V_0)(V_0^*)\oplus
V_0(S^2V_0^*)\oplus S^3V_0^*.
$$
Therefore, $E^3V$ is the direct
sum of irreducible $\fsl(3)$-modules with the highest  weights and
highest vectors listed in  Table 7.  Note that the
vectors of weights 0 and $\eps_1$ are the $\fspe(3)$-highest
  ones.
Therefore, the Jordan-H\"{o}lder series of the $\fspe(3)$-module
$E^3V$ contains as quotient modules the trivial and the standard
ones. Notice that the vector of weight $-2\eps_3$ is the
$\fspe(3)$-highest one in the corresponding quotient module which
can only contain $\fsl(3)$-submodules with highest weights
$-2\eps_3,\hbox{ } \eps_1 - 2\eps_3,\hbox{ and } -3\eps_3.$

Since $v_{2\eps_1 + 2\eps_2}$ is the
$\fsl(3)$-highest vector of weight $-2\eps_3$, then the relations
(3.11.1)--(3.11.4) imply part b) of Lemma 3.11.
\end{proof}

\ssbegin{3.12}{Lemma}  For $n = 3$ we have the following nonsplit exact
sequence of $\fspe(3)$-modules
$$
0\longrightarrow\Pi (X_{2\eps_1 + 2\eps_2})\longrightarrow
H^{2,2}_{\fspe(3)\subplus \langle \tau + 3z\rangle}\longrightarrow \Pi
(X_{3\eps_1})\longrightarrow 0.\eqno (3.12.1)
$$
\end{Lemma}

\begin{proof} By part b) of Lemma 3.11 and (3.11.5) we see that
$H^{2,2}_{\fspe(3)\subplus \langle \tau + 3z\rangle}$ contains an irreducible
$\fspe(3)$-module with highest weight
  $-2\eps_3$, which being considered as an $\fsl(3)$-module, is the sum of
irreducible $\fsl(3)$-components with highest weights $-2\eps_3$,
$\eps_1-2\eps_3,
  \hbox{ and} -3\eps_3$.

In addition to these $\fsl(3)$-components,
$H^{2,2}_{\fspe(3)\subplus \langle \tau + 3z\rangle}$ can only
contain irreducible
$\fsl(3)$-components with the following highest
  weights:
$$
3\eps_1,\hbox{ } 2\eps_1,\hbox{ and } 2\eps_1 - \eps_3.
$$
Let us show that these components are  indeed contained in
$H^{2,2}_{\fspe(3)\subplus \langle \tau + 3z\rangle}$, and that their sum is an
irreducible $\fspe(3)$-quotient module with highest weight
$3\eps_1$.

First, note that $E^{0,2}_\infty =
E^{0,2}_4$ has an irreducible $\fgl(n)$-module with highest
weight $2\eps_1 - \eps_{n-1} - \eps_n$. In fact, by Table 5 such a
submodule is contained in $E^{0,2}_1$ but is not contained in
$E^{1,2}_1$. Therefore, by (3.2.6) it is contained in $E^{0,2}_2$.

According to Table 4, in $E^{2,1}_1$ there is no
submodule with highest weight $2\eps_1 - \eps_{n-1} - \eps_n$,
hence such a submodule is not contained in $E^{2,1}_2$ either.
Therefore, by (3.2.6) the submodule with this highest weight is
contained in $E^{0,2}_3.$

According to Table 6,
$E^{3,0}_0$ has no submodule with highest weight $2\eps_1 -
\eps_{n-1}-\eps_n$, hence it is not contained in $E^{3,0}_3$
either. Therefore, by (3.2.6) it is contained in $E^{0,2}_\infty =
  E^{0,2}_4.$

By Tables 2 and 5 the $\fgl(n)$-highest vectors in
$E^{0,2}_1$ of weights
$$
-2\eps_{n-1} - 2\eps_n,\hbox{ } 2\eps_1 - \eps_{n-1} - \eps_n,
\hbox{ and } \eps_1 - \eps_{n-1} -2\eps_n
$$
are, respectively,
$$
\eqalign {
&v_{-2\eps_{n-1}-2\eps_n} = (f_{n-1}\wedge \tilde f_n)\otimes
(\tilde f_{n-1}\wedge \tilde f_n),\cr
&v_{2\eps_1-\eps_{n-1}-\eps_n} = (e_1 \tilde e_1)\otimes (\tilde
f_{n-1}\wedge \tilde f_n), \hbox{ and }\cr
&v_{\eps_1-\eps_{n-1}-2\eps_n} =
  (e_1\wedge \tilde f_n)\otimes (\tilde f_{n-1}\wedge \tilde f_n).\cr}
$$
Note that
$$
\eqalign {
&C_{1,n}v_{2\eps_1 - \eps_{n-1} - \eps_n} =
-2v_{\eps_1-\eps_{n-1}-2\eps_n},\cr
&B_{1,n}v_{\eps_1-\eps_{n-1}-2\eps_n} =
-v_{2\eps_1-\eps_{n-1}-\eps_n}- (e_1\wedge \tilde f_n)\otimes
(\tilde f_{n-1}\wedge \tilde e_1),\cr }
$$
$$
\eqalign {
&C_{1,n-1}v_{\eps_1-\eps_{n-1}-2\eps_n }=
v_{-2\eps_{n-1}-2\eps_n},\cr &B_{1,n-1}v_{-2\eps_{n-1}-2\eps_n} =
v_{\eps_1-\eps_{n-1}-2\eps_n} + (f_{n-1}\wedge \tilde f_n)\otimes
(\tilde e_1\wedge \tilde f_n).\cr }
$$

Therefore, for $n = 3$ the
components with highest weights $3\eps_1, \hbox{ } 2\eps_1,\hbox{
and } 2\eps_1 - \eps_3$ constitute an irreducible quotient module
with highest weight $3\eps_1$.

 From Tables 2 and 4 we see that the $\fgl(n)$-highest vector in $E^{1,1}_1$
of weight $2\eps_1+\eps_2-\eps_n$ is
$$
v_{2\eps_1+\eps_2-\eps_n} =
(e_1\tilde e_2)\otimes (\tilde e_1\wedge \tilde f_n)-
(e_1\tilde e_1)\otimes (\tilde e_2\wedge \tilde f_n).
$$
Observe that
$$
(B_{1,2}A_{2,1}-(1/2)B_{2,2})(v_{2\eps_1 -\eps_2-\eps_3}) =
2v_{2\eps_1+\eps_2-\eps_3}.
$$
Therefore, the sequence (3.12.1) is
nonsplit. This proves Lemma 3.12, and Theorem 2.1 in the case
where $\fg_0 = \fspe(3)\subplus \langle \tau + 3z\rangle$.
\end{proof}

Recall that the Jordan-H\"{o}lder series of $\fspe(n)$-module
$H^{2,2}_{\fcpe(n)}$, as compared to that of $H^{2,2}_{\fspe(n)\subplus
\langle \tau + nz\rangle}$, contains in addition the $\fspe(n)$-component with
highest weight $2\eps_1$. Recall that by (3.10.1) the highest
$\fgl(n)$-vector with weight $2\eps_1$ in $E^{1,1}_1$ is
$$
u_{2\eps_1} = 2\sum_{i=1}^n e_1\tilde e_i\otimes \tilde f_i\wedge
\tilde e_1-
(n+1)\sum_{i=1}^n e_1\tilde e_1 \otimes \tilde f_i\wedge \tilde e_i.
$$
Note that
$$
C_{n-1,n}(u_{2\eps_1}) = -2(n+1)v_{2\eps_1-\eps_{n-1}-\eps_n}
+2(e_1\wedge \tilde f_{n-1}\otimes \tilde f_n\wedge \tilde
e_1-e_1\wedge \tilde f_n\otimes \tilde f_{n-1}\wedge \tilde
e_1).
$$
This proves Theorem  2.1 in the case where $\fg_0 =
\fcpe(n)$.

The proof of Theorem  2.1 in the case where $\fg_0 = \fspe(n)$
follows from the fact that the Jordan-H\"{o}lder series of
$\fspe(n)$-module $H^{2,2}_{\fspe(n)}$, as compared to that of
$H^{2,2}_{\fspe(n) \subplus \langle \tau + nz\rangle}$, contains in
addition the
$\fspe(n)$-component with highest weight $\eps_1+\eps_2$.

Finally, the proof of Theorem  2.1 in the case when
  $\fg_0 = \fspe(n)\subplus \langle a\tau + bz\rangle$, where $a$, $b
\in\Cee$ are such that
$a = 0$, $b \not= 0$, or $a \not= 0$, $b/a \not= n$, follows from
the fact that the Jordan-H\"{o}lder series of $\fspe(n)$-module
$\Ker\partial ^{2,2}_{\fspe(n)\subplus \langle a\tau + bz\rangle}$,
as compared to
that of $\Ker\partial ^{2,2}_{\fspe(n)}$, contains in addition the
$\fspe(n)$-component with highest weight $2\eps_1$.

%%%%%%%%%%%%%%%%%%%%%%%%%%%%%%%%%%%%%%%%%%%%%%%%%%%%%%%%%%%%%%%%%%%%%%%%%%%%%%%%
%%%%%%%%%%%%%%%%%%%%%%%%%%%%%%%%%%%%%%%%%%%%%%%%%%%%%%%%%%%%%%%%%%%%%%%%%%%%%%%%%%%%%%%%%%%%%%%%%%%%%%%%%%%%
%Chapter 2 %%%%%%%%%%%%%%%%%%%%%%%%%%%%%%%%%%%%%%%%%%%%%%%%%%%%%%%%%%%%%%%%%%%%%%%%%%%%%%%%%%%%%%%%%%%%%%%%%%%%%%%%%%%%%%%%%%%%%%%%%
%%%%%%%%%%%%%%%%%%%%%%%%%%%%%%%%%%%%%%%%%%%%%%%%%%%%%%%%%%%%%%%%%%%%%%%%%%%%%%%%

\chapter{The analogues of Penrose's tensors}

In this chapter
$\Lambda^iV$  is the $i$-th exterior power of a vector
(super)space $V$.

\section*[Standard $\Zee$-grading of $\fsl(m|n)$]{Standard
$\Zee$-grading of $\fsl(m|n)$
and the corresponding Cartan prolongations}

Let $V
= V(m|0) $ and $U = U(0|n)$ be the standard (identity) $\fgl(m)$-
and $\fgl(n)$-modules. (Hereafter $\fgl(m) = \fgl(m|0)$,
$\fgl(n) = \fgl(0|n)$, etc.)

In what follows we will consider the standard (compatible)
$\Zee$-grading of $\fg =
\fsl(m|n)$ with $m \leq n$ and let the degrees of all even roots
be zero. This yields the $\Zee$-grading of the form:
$$
\fg = \fg_{-1} \oplus \fg_0 \oplus \fg_1,\hbox{ where }
\fg_0 = \fsl(m)\oplus \fsl(n) \oplus\Cee,\hbox{ } \fg_{-1} =
\fg_1^* = U\otimes V^*.
$$
Let $\hat{\fg_0}$ be the Levi subalgebra of $\fg_0$, i.e., $\hat{\fg_0} =
\fsl(m) \oplus \fsl(n)$. The weights are given with respect to the
bases $\eps_1, \ldots, \eps_m$ and $\delta_1,\ldots,\delta_n$ of
the dual spaces to the maximal tori of $\fgl(m|n)$. Let
$e_1,\ldots,e_m$ be the weight basis of $V$ and $f_1,\ldots,f_n$
be the weight basis of $U$. Let $\tilde e_1,\ldots, \tilde e_m$
and $\tilde f_1,\ldots,\tilde f_n$ be the bases of the dual spaces
to $V$ and $U$, respectively, normed so that $\tilde e_i(e_j) =
\tilde f_i(f_j) = \delta_{ij}$. If $\oplus_\lambda k_\lambda
V_\lambda$ is a direct sum of irreducible $\fg_0$-modules (here
  $k_\lambda $ is the multiplicity of $V_\lambda$) with highest weight
$\lambda$,
denote by $v^i_\lambda$ the highest weight vectors of the
corresponding components: $i = 1,\ldots, k_\lambda$. We will often
represent the elements of $\fgl(m|n)$ by the matrices
$$
X = \diag(A,\hbox{ } D) + \antidiag(B,\hbox{ } C)
$$
where the dimensions of the matrices $A$, $B$, $C$, and $D$ are $m
\times m$,  $m\times n$,  $n \times m$ and  $n \times n$,
respectively. Denote by $A_{i,j}$ the matrix $X$ whose components
$B, C, \hbox{ and }D$ are zero and all the entries of $A$ are also
zero except for the $(i, \hbox{ }j)$-th. The matrices $B_{i,j}$,
$C_{i,j}$, and $D_{i,j}$ are defined similarly.

\ssbegin{1.1}{Theorem}\emph{a)} If $m = 1,\hbox{ } n > 1$, then
$\fg_*(\fg_{-1},\hbox{ } \fg_0) = \fvect(0|n)$, $\fg_*(\fg_{-1},
\hbox{ }\hat {\fg_0}) = \fsvect(0|n)$;

\emph{b)} if
$m,\hbox{ }n >1$ and $m\not= n$, then $\fg_*(\fg_{-1},\hbox{ }
\fg_0) = \fg$,  $\fg_*(\fg_{-1}, \hbox{ }\hat {\fg_0}) = \fg_{-1}
\oplus\hat {\fg_0};$

\emph{c)} if $m = n = 2$, then
$\fg_*(\fg_{-1},\hbox{ } \hat{\fg_0}) = \fh(0|4)$,
$\fg_*(\fg_{-1},\hbox{ } \fg_0) = S^*(\fg _{-1}^*) \subplus \fh(0|4);$

\emph{d)} if $m = n > 2$, then $\fg_*(\fg_{-1},\hbox{ }
\hat{\fg_0}) = \fpsl (n|n)$, $\fg_*(\fg_{-1},\hbox{ } \fg_0) =
S^*(\fg_{-1}^*) \subplus \fpsl (n|n).$
\end{Theorem}
\begin{proof}
Consider all cases mentioned in Theorem 1.1.

\ssec{1.2.} $\bf m = 1,\hbox{ } n \geq 2.$ Then $\fg_0 = \fsl(n) \oplus
\Cee = \fgl(n)$ and $\hat {\fg_0} = \fsl(n)$, where $\fg_{-1}$ is
the standard $\fg_0$ (or $\hat {\fg_0}$) module. Therefore,
$\fg_*(\fg_{-1}, \hbox{ }\fg_0) = \fvect (0|n),\hbox{ }
\fg_*(\fg_{-1},\hbox{ } \hat {\fg_0}) = \fsvect (0|n).$

Notice that if $m \not= n$, then
$$
\fsl (m|n)\subset \fg_*(\fg_{-1}, \hbox{ } {\fg_0}) \eqno (1.2.1)
$$
and if $m = n$, then
$$
\fpsl (n|n)\subset  \fg_* (\fg_{-1},\hbox{ } \hat {\fg_0}).\eqno(1.2.2)
$$
Indeed, the Lie superalgebras $\fsl(m|n)$, where $m\not= n$, and
$\fpsl(n|n)$ are simple and therefore, they are transitive (i.e.,
if there exists $g\in \fg_i$ ($i \geq 0$) such that
$[\fg_{-1},\hbox{ }g] = 0$, then $g = 0$. It follows that $\fg_1$
is embedded into $\fg_0\otimes \fg_{-1}^*$
  (or $\hat {\fg_0}\otimes \fg_{-1}^*$). The Jacobi identity implies
$\fg_1\subset \fg_{-1}\otimes S^2\fg_{-1}^*.$

\ssec{1.3. Calculation of the first term of the Cartan
prolongation for $m,\hbox{ } n \geq 2,\hbox{ }
m\not=n$.} Let $\fg_1'$ be the first term of the Cartan
prolongation of the pair $(\fg_{-1},\hbox{ } \fg_0)$. Let us show
that $\fg_1' = \fg_1$. By definition,
$$
\eqalign{
& \fg_1' = (\fg_0\otimes\fg_{-1}^*)\cap (\fg_{-1}\otimes
S^2\fg_{-1}^*), \hbox{ where, as }
\fgl(m)\oplus\fgl(n)\hbox{-modules },\cr &\fg_0\otimes
\fg_{-1}^* \cong [(V\otimes V^*)/\Cee \oplus (U\otimes U^*)/\Cee
\oplus \Cee] \otimes (U^*\otimes V).\cr }
$$
Note that if  $g\in
\fg_0\otimes\fg_{-1}^*$, then
$$
g\in \fg_{-1}\otimes S^2\fg_{-1}^*\hbox{ if and only if }g(g_1)(g_2) =
-g(g_2)(g_1) \hbox{ for any } g_1, g_2\in \fg_{-1},
$$
since
$\fg_{-1}$ is purely odd.

\begin{Lemma} The
$\fgl(m)\oplus \fgl(n)$-module $\fg_0\otimes \fg_{-1}^*$ is the
direct sum of irreducible submodules whose highest weights and
highest vectors are listed in Table 8.
\end{Lemma}
\begin{proof}
The proof of the Lemma consists of:
\emph{a)} a verification
of the fact that  vectors $v_\lambda$ from Table 8 are indeed
highest with respect to $\fgl(m)\oplus \fgl(n)$, i.e.
$A_{i,j}v_\lambda = D_{i,j}v_\lambda = 0$ for $i<j$;

\emph{b)} a calculation of dimension of $\fg_0\otimes \fg_{-1}^*$ and of
dimensions of the irreducible submodules of $\fg_0\otimes
\fg_{-1}^*$ by the formula from the Appendix.

Let us show that if
$$
\lambda = 2\eps_1 - \eps_m - \delta_n,\hbox{ }
\eps_1 + \delta_1 - 2\delta_n,\hbox{ } \eps_1 + \eps_2 - \eps_m -
\delta_n \hbox{ }(\hbox{if } m\geq 3),
$$
or
$$
\lambda = \eps_1 + \delta_1 - \delta_{n-1} - \delta_n \hbox{ }(\hbox{if }
  n\geq 3),
  $$
then $v_\lambda\not \in \fg_1'$. For this it suffices to indicate
$g_1, g_2\in \fg_{-1}$ such that
$$
v_\lambda(g_1)(g_2) \not= -v_\lambda (g_2)(g_1) \eqno (1.3.1)
$$
or, perhaps, there exists just one $g\in\fg_{-1}$ such that
$$
v_\lambda(g)(g) \not = 0. \eqno (1.3.2)
$$
Let $\lambda = 2\eps_1 - \eps_m - \delta_n$. Then
$$
v_\lambda (f_n\otimes \tilde e_1)(f_n\otimes \tilde e_1) =
A_{1,m}(f_n\otimes \tilde e_1) = -f_n\otimes \tilde e_m \not= 0.
$$
If $\lambda = \eps_1 + \delta_1 -2\delta_n$, then
$$
v_\lambda(f_n\otimes \tilde e_1)(f_n\otimes \tilde e_1) =
D_{1,n}(f_n\otimes \tilde e_1) = f_1\otimes \tilde e_1 \not= 0.
$$
If $\lambda = \eps_1 + \eps_2 - \eps_m - \delta_n$ (for $m\geq
3$), then
$$
\eqalign{
&v_\lambda(f_n\otimes \tilde e_2)(f_{n-1}\otimes \tilde e_1) =
A_{1,m}(f_{n-1}\otimes \tilde e_1) = -f_{n-1}\otimes \tilde
e_m,\cr &\hbox{ but } v_\lambda(f_{n-1}\otimes \tilde
e_1)(f_n\otimes \tilde e_2) = 0.\cr }
$$
Finally, if $\lambda =
\eps_1 + \delta_1 - \delta_{n-1} - \delta_n$ (for $n\geq 3$), then
$$
\eqalign{
&v_\lambda(f_n\otimes \tilde e_1)(f_{n-1}\otimes \tilde e_2) =
D_{1,n-1}(f_{n-1}\otimes \tilde e_2) = f_1 \otimes \tilde e_2,\cr
&\hbox{ but } v_\lambda (f_{n-1}\otimes \tilde e_2)(f_n\otimes
\tilde e_1) = 0.\cr }
$$
Now, let us show that if $\lambda = \eps_1
- \delta_n$, then $\fg_1'$ contains precisely one irreducible
$\fgl(m)\oplus \fgl(n)$-module with highest weight $\lambda$.
Notice that by (1.2.1) $\fg_1'$ contains at least one such module.
Let
$$
v_\lambda = k_1 v_\lambda^1 + k_2v_\lambda^2 + k_3v_\lambda^3,
\hbox{ where } k_1, k_2, k_3\in \Cee,
$$
be a linear combination of
highest vectors of weight $\lambda$. Then the condition
$$
v_\lambda (f_n\otimes \tilde e_2)(f_{n-1}\otimes \tilde e_1) =
-v_\lambda(f_{n-1}\otimes \tilde e_1)(f_n\otimes \tilde e_2)
$$
implies
$$
mk_1 = nk_2, \eqno (1.3.3)
$$
whereas the condition
$$
v_\lambda(f_n\otimes\tilde e_1)(f_n\otimes \tilde e_1) = 0
$$
implies
$$
k_1(m-1) + k_2(1-n) + k_3(m-n) = 0.
$$
Hence,
$$
k_2 = mk_1/n \hbox{ and } k_3 = -k_1/n. \eqno (1.3.4)
$$
Thus, $\fg_1' = V_{\eps_1 - \delta_n} \hbox { and } \fg_1' =
\fg_1.$
\end{proof}
\ssec{1.4. Calculation of the second
term of the Cartan prolongation for $m,\hbox{ } n\geq
2,\hbox{ } m\not= n$.}  Let $\fg_2$ be the second term of the
Cartan prolongation of $(\fg_{-1}, \hbox{ }\fg_0)$.

Let us show that $\fg_2 = 0$. Indeed, by definition, $\fg_2 =
(\fg_1\otimes\fg_{-1}^*)\cap (\fg_0 \otimes S^2\fg_{-1}^*)$.
Notice that, as $\fg_0$-module,
$$
\eqalign{
&\fg_1\otimes \fg_{-1}^* \cong (U^*\otimes V)\otimes (U^*\otimes
V) = \cr &S^2U^*\otimes S^2V\oplus \Lambda^2U^*\otimes \Lambda^2V
\oplus \Lambda^2U^*\otimes S^2V \oplus S^2U^*\otimes
\Lambda^2V.\cr }
$$
This decomposition and Table 5 of [OV] imply
the following

\begin{Lemma} The $\fgl(m) \oplus
\fgl(n)$-module $(U^*\otimes V)\otimes (U^*\otimes V)$ is the
direct sum of irreducible submodules whose highest weights and the
corresponding highest vectors are listed in Table 9.
\end{Lemma}

Let us show that $v_\lambda \not \in \fg_2$, where  $v_\lambda$ is
any of the highest vectors listed in Table 9.
  Let us  indicate $g_1, g_2\in \fg_{-1}$ for which either (1.3.1) or
(1.3.2) holds.

Let $\lambda = 2\eps_1 - 2\delta_n.$ Then
$$
v_\lambda(f_n\otimes \tilde e_1)(f_n\otimes \tilde e_2) =
B_{1,n}(f_n\otimes \tilde e_2) = e_1\otimes \tilde e_2, \hbox{ but
} v_\lambda (f_n\otimes \tilde e_2)(f_n\otimes \tilde e_1) = 0.
$$
If $\lambda = \eps_1 + \eps_2 - 2\delta_n$, then
$$
v_\lambda (f_n\otimes \tilde e_2)(f_n\otimes \tilde e_2) =
B_{1,n}(f_n\otimes \tilde e_2) = e_1 \otimes \tilde e_2 \not= 0.
$$
If $\lambda = 2\eps_1 - \delta_{n-1} - \delta_n$, then
$$
v_\lambda (f_{n-1}\otimes \tilde e_1)(f_n\otimes \tilde e_2) =
-B_{1,n}(f_n\otimes \tilde e_2) = -e_1\otimes \tilde e_2, \hbox{
but } v_\lambda (f_n\otimes \tilde e_2)(f_{n-1}\otimes \tilde e_1)
= 0.
$$
Let $\lambda = \eps_1 + \eps_2 - \delta_{n-1} - \delta_n$.
Then if $n>2$, we have
$$
v_\lambda(f_{n-1}\otimes \tilde e_2)(f_1\otimes \tilde e_1) =
B_{1,n}(f_1\otimes \tilde e_1) = f_1 \otimes \tilde f_n,
$$
but $v_\lambda(f_1\otimes \tilde e_1)(f_{n-1}\otimes \tilde e_2) = 0.$
If $m>2$, then
$$
v_\lambda (f_{n-1}\otimes \tilde e_2)(f_n\otimes \tilde e_m) =
B_{1,n}(f_n\otimes \tilde e_m) = e_1 \otimes \tilde e_m,
$$
but
$v_\lambda(f_n\otimes \tilde e_m)(f_{n-1}\otimes \tilde e_2) = 0.$
Therefore, $\fg_2 = 0$ and $\fg_*(\fg_{-1},\hbox{ } \fg_0) = \fg$.
Note that by (1.3.4) we have $\fg_*(\fg_{-1},\hbox{ } \hat
{\fg_0}) = \fg_{-1} \oplus \hat {\fg_0}$. This proves part b) of
Theorem 1.1.

\ssec{1.5.} $\bf m = n.$ Let $m = n = 2$
Since $\hat {\fg_0} = \fsl(2) \oplus \fsl(2)
  =\fo(4)$ and $\fg_{-1}$ is the standard $\fo(4)$-module (considered as
  purely odd superspace), then $\fg_*(\fg_{-1},\hbox{ } \hat{\fg_0}) =
\fh(0|4).$

Let $m = n >2$ and $\fg'_1$ be the first term of the
Cartan prolongation of the pair  $(\fg_{-1},\hbox{ }
\hat{\fg_0})$. Let us show that $\fg_1' = \fg_1$.
Indeed, by (1.2.2) $\fg_1\subset \fg_1'$. By sec.1.3 and Table 8
we see that the only highest weights of $\fg_1'$ are all equal to
$\eps_1 - \delta_n$. Then formula (1.3.3) implies that the highest
vector of such weight in $\fg_1'$ is precisely one and therefore,
$\fg_1' = \fg_1$. By sec.1.4 the second term of the Cartan
prolongation of the pair $(\fg_{-1},\hbox{ } \hat{\fg_0})$ is
zero. Hence, for $m = n >2$ we have $\fg_*(\fg_{-1}, \hbox{
}\hat{\fg_0}) = \fpsl(n|n)$.

Let $m = n >1$ and let
$\fg_k$ be the $k$-th term of the  Cartan prolongation of the pair
$(\fg_{-1},\hbox{ } \fg_0)$. Recall that $\fg_k = (\fg_0 \otimes
S^k\fg_{-1}^*) \cap (\fg_{-1}\otimes S^{k+1}\fg_{-1}^*)\hbox{ }
(k\geq1)$. Observe that
$$
\fg_0 \otimes S^k\fg_{-1}^* = (\hat {\fg_0} \oplus \langle z\rangle) \otimes
S^k\fg_{-1}^*,\hbox{ where  }z = 1_{2n} \hbox{ is the center of }
\fsl(n|n).
$$
Note that
$$
\langle z\rangle\otimes S^k\fg_{-1}^*\subset \fg_{-1}\otimes S^{k+1}\fg_{-1}^*.
$$
Therefore,
$$
\fg_*(\fg_{-1},\hbox{ } \fg_0) = S^*(\fg_{-1}^*)\subplus \fg_*(\fg_{-1},
  \hat{\fg_0}).
$$
\end{proof}

%%%%%%%%%%%%%%%%%%%%%%%%%%%%%%%%%%%%%%%%%%%%%%%%%%%%%%%%%%%%%%%%%%%%%%%%%%%%%%%%%%%%%%%%%%%%%%%%%%%%%%%%%%%%%%%%
%\S2
%%%%%%%%%%%%%%%%%%%%%%%%%%%%%%%%%%%%%%%%%%%%%%%%%%%%%%%%%%%%%%%%%%%%%%%%%%%%%%%%%%%%%%%%%%%%%%%%%%%%%%%%%%%%%%%%%%%%%%%%%%%

\section*{Structure functions of Lie superalgebras
$\fvect(0|n)$ and $\fsvect(0|n)$}

Recall that $U$ is  the purely odd standard
$\fgl(n)$-module.

\ssbegin{2.1}{Theorem} If $m = 1, n >1$, then
$$
\eqalign {
&\hbox{a) } H^{k,2}_{\fg_0} = 0 \hbox{ for any } k > 0;\cr
&\hbox{b) } H^{k,2}_{\hat {\fg_0}} = \Pi
^n(\Cee)\delta_{kn}.\cr}
$$
\end{Theorem}
\ssec{2.2. Proof of part a).} Since
$\fg_*(\fg_{-1}, \fg_0) = \fvect (0|n) = \sum_{i=-1}^{n-1}\fg_i$,
where $\fg_0 = \fgl(n)$, and the $\fgl(n)$-module $\fg_i$ is
isomorphic to
  $U\otimes S^{i+1}U^*$, then for $k \geq n+2$ we have
$H^{k,2}_{\fg_0} = 0$, and for $k \leq n+1$ there exist
  the following Spencer cochain sequences:
$$
\begin{array}{ll}
   \fgl(n) \otimes U^*\buildrel {\partial^{2,1}_{\fgl(n)}} \over
   \longrightarrow U\otimes \Lambda^2U^*
   \buildrel{\partial^{1,2}_{\fgl(n)}} \over \longrightarrow 0 &(k = 1),\\
   C^{k+1,1}_{\fgl(n)}\buildrel{\partial^{k+1,
   1}_{\fgl(n)}} \over \longrightarrow C^{k,2}_{\fgl(n)}
   \buildrel{\partial^{k,2}_{\fgl(n)}}\over \longrightarrow
   C^{k-1,3}_{\fgl(n)}&(2\leq k \leq n),\\
   0\buildrel{\partial^{n+2,1}_{\fgl(n)}} \over \longrightarrow
   C^{n+1,2}_{\fgl(n)} \buildrel{\partial^{n+1,2}_{\fgl(n)}} \over
   \longrightarrow C^{n,3}_{\fgl(n)} &(k = n + 1),
\end{array}
$$
where
$$
\eqalign{
&C^{k+1,1}_{\fgl(n)} = \fg_{k-1}\otimes \fg_{-1}^*  \cong
(U\otimes S^kU^*)\otimes U^*,\cr &C^{k,2}_{\fgl(n)} =
\fg_{k-2}\otimes \Lambda^2\fg_{-1}^*\cong (U\otimes
S^{k-1}U^*)\otimes \Lambda^2U^*,\cr &C^{k-1,3}_{\fgl(n)} =
\fg_{k-3}\otimes \Lambda^3\fg_{-1}^*\cong (U\otimes
S^{k-2}U^*)\otimes \Lambda^3U^*.\cr }
$$
Recall (see sec.3.1 of
Chapter 1) that if $\fg_{-1}$ is a faithful module over a Lie
superalgebra $\fg_0$ and
$$
\fg_{k-1}\otimes \fg_{-1}^* \buildrel{\partial^{k+1,1}_{\fg_0}}
\over \longrightarrow\fg_{k-2} \otimes \Lambda^2\fg_{-1}^*
\buildrel {\partial^{k,2}_{\fg_0}} \over \longrightarrow \fg_{k-3}
\otimes \Lambda ^3\fg_{-1}^* \eqno (k \geq 1)
$$
is the  Spencer
cochain sequence  , which corresponds to the pair
$(\fg_{-1},\hbox{ } \fg_0)$, then
$$
\Im \partial^{k+1,1}_{\fg_0}
\cong (\fg_{k-1} \otimes \fg_{-1}^*) /\fg_k. \eqno(2.2.1)
$$
Let us show that $H^{1,2}_{\fg_0} = 0$. Indeed, by (2.2.1)
$$
\Im \partial^{2,1}_{\fg_0} \cong
  (\fgl(n)\otimes U^*)/(U\otimes S^2U^* )\cong
U\otimes \Lambda^2U^* = \Ker \partial^{1,2}_{\fg_0}.
$$
We will
prove that for $2 \leq k \leq n + 1$ we have $H^{k,2}_{\fgl(n)} =
0$ , using the following

\begin{Lemma} As $\fgl(n)$-modules $C^{k,2}_{\fgl(n)}, \hbox{
}C^{k+1,1}_{\fgl(n)}, \hbox{ and }\fg_k$ , where $2\leq k \leq
n+1$, are the direct sums of the irreducible submodules whose
highest weights and highest vectors are listed in Tables 10, 11,
and 12, respectively. [$r$, $s$, and $t$ denote the cyclic
permutations of $(n-k,\ldots,n)$, $(n-k+1,\ldots,n)$, and
$(n-k+2,\ldots,n)$, respectively.]
\end{Lemma}
\begin{proof}  The
proof of the Lemma  consists of

\emph{a)} a verification of
the fact that vectors $v_{\lambda_l}, v_{\beta_l},
  v_{\gamma_l}$ from the Tables 10, 11, and 12 are indeed highest with
respect to $\fgl(n)$,
i.e., $D_{i,j}v_{\lambda_l} = D_{i,j}v_{\beta_l} = D_{i,j}v_{\gamma_l} = 0$ for
$i < j$;

\emph{b)} a calculation of dimensions of given modules and dimensions
of their irreducible submodules by the formula from the Appendix.

Notice that $D_{i,j}f_p = \delta_{jp}f_i,\hbox{ } D_{i,j}\tilde{f_p} =
  -\delta_{ip}\tilde{f_j}$. According to Table 10,
$$
\lambda_1 = \delta_1 - \delta_{n-k+2} - \ldots - \delta_{n-1} - 3\delta_n,
\hbox{ }
v_{\lambda_1} = (f_1 \otimes \tilde f_{n-k+2} \wedge \tilde f_{n-k+3}
\wedge \ldots \wedge \tilde {f_n}) \otimes \tilde {f_n}^2.
$$

Then for $i < n - k + 2$  $D_{i,j}v_{\lambda_1} = 0.$ For $j > i \geq
n - k + 2$
$$
\eqalign{
&D_{i,j}(v_{\lambda_1}) = D_{i,j}(f_1\otimes \tilde f_{n-k+2}
\wedge \ldots \wedge \tilde{f_i}\wedge \ldots \wedge \tilde{f_j}
\wedge \ldots \wedge \tilde{f_n})\otimes \tilde{f_n}^2 = \cr &(f_1
\otimes \tilde f_{n-k+2} \wedge \ldots \wedge (-\tilde{f_j})
\wedge \ldots \wedge \tilde {f_j} \wedge \ldots \wedge \tilde
{f_n}) \otimes \tilde {f_n}^2 = 0.\cr }
$$
Thus, $v_{\lambda_1}$ is
highest vector. The proof of the fact that the other
vectors from Tables 10, 11, and 12 are highest with respect to
$\fgl(n)$ is similar.

Using the formula from the
Appendix we find the dimensions of the $\fgl(n)$-modules given in
Table 10:

if $2 \leq k \leq n $, then
$$
\dim V_{\lambda_1} ={n(n+3)(n+1)!\over
2(n-k+2)(k+1)(k-2)!(n-k)!},
$$
if $3 \leq k = n + 1$, then
$$
\dim V_{\lambda_1} = {(n-1)n(n+2)\over 2},
$$
if $2 \leq k\leq n + 1$, then
$$
\dim V_{\lambda_2} = {(n+1)!\over (k-2)!((n-k+1)!k},
$$
if $k = 2 \leq n$, then
$$
\dim V_{\lambda_3} = \dim V_{\lambda_2},
$$
if $3 \leq k \leq n$, then
$$
\dim V_{\lambda_3} = {(n+2)!\over 2(k-3)!(n-k+2)!k},
$$
if $3 \leq k = n$, then
$$
\dim V_{\lambda_4} = n^2 - 1,
$$
if $2 \leq k\leq n - 1$, then
$$
\dim  V_{\lambda_4} = {n(n+2)n!\over
(n-k+1)(k-1)!(n-k-1)!(k+1)}, \eqno (2.2.2)
$$
if $2 \leq k \leq n$, then
$$
\dim V_{\lambda_5} = {n!\over k!(n-k)!}, \eqno (2.2.3)
$$
if $3 \leq k \leq n-1$, then
$$
\dim V_{\lambda_6} = \dim V_{\lambda_2}.
$$
Therefore, if  $2 = k = n$, then
$$\dim V_{\lambda_1} + 2dimV_{\lambda_2} + dimV_{\lambda_5} =
{n^3(n+1)\over 2} =
\dim(U\otimes U^*) \otimes \Lambda^2U^*,$$ if $2 = k \leq n - 1$,
then
$$
\sum_{l=1}^5 \dim V_{\lambda_l} = {n^3(n+1)\over 2} =
\dim(U\otimes U^*)\otimes \Lambda^2U^*,
$$
if $3 \leq k = n + 1$,
then
$$
\sum_{l=1}^2 \dim V_{\lambda_l} = {n^2(n+1)\over 2} = \dim U\otimes
\Lambda^2U^*,
$$
if $3 \leq k = n$, then
$$
\sum _{l=1}^5 \dim V_{\lambda_l} = {n^3(n+1)\over 2} =
\dim (U\otimes S^{n-1}U^*)\otimes \Lambda^2U^*,
$$
if $ 3 \leq k \leq n - 1$, then
$$
\sum_{l=1}^6 \dim V_{\lambda_l} =
{n^2(n+1)!\over 2(n-k+1)!(k-1)!} = \dim (U\otimes
S^{k-1}U^*)\otimes \Lambda^2U^*.
$$

In order to find the dimensions
of the $\fgl(n)$-modules given in Table 11, note that if $k \geq
2$, then $\beta_1 = \lambda_4, \beta_2 = \beta_4 = \lambda_5,
\beta_3 = \lambda_2.$ Using the formula from the Appendix we get
$$
\dim V_{\beta_5} = {(n+1)!\over (n-k)(k+1)!(n-k-2)!}\hbox{ }\hbox{ for }
  2 \leq k \leq  n - 2. \eqno (2.2.4)
$$
Therefore,
if $2 \leq k = n$, then
$$
\sum_{l=1}^2 \dim V_{\beta_l} =  n^2 = \dim U\otimes U^*,
$$
if $2 \leq k = n - 1$, then
$$
\sum_{l=1}^4 \dim V_{\beta_l} = n^3 =
\dim (U\otimes S^{n-1}U^*)\otimes U^*,
$$
if $2 \leq k \leq n - 2$,
then
$$
\sum_{l=1}^5 \dim V_{\beta_l} = {n^2n!\over (n-k)!k!} =
\dim (U\otimes S^kU^*)\otimes U^*.
$$
Finally, in order to find the
dimensions of the $\fgl(n)$-modules given in Table 12, note that
if $2 \leq k \leq n - 2$, then $\gamma_1 = \beta_5$, and $\gamma_2
=\lambda_5.$ Therefore, if $2 \leq k \leq n - 2$, then
$$
\sum _{l=1}^2 \dim V_{\gamma_l} =
{nn!\over (n-k-1)!(k+1)!} = \dim (U\otimes S^{k+1}U^*) = \dim \fg_k,
$$ and if $k = n - 1$, then $\dim V_{\gamma_1} = n = \dim U
= \dim \fg_k.$ This proves  the Lemma.
\end{proof}

Let $k \geq 2$ and $\lambda = \lambda_1,\hbox{ } \lambda_2,\hbox{ }
\lambda_3.$ Then
$$
v_\lambda \not \in \Ker \partial^{k,2}_{\fgl(n)}. \eqno (2.2.5)
$$
Indeed,
$$
\partial^{k,2}_{\fgl(n)}v_\lambda (f_n, f_n, f_n) = -3v_\lambda
(f_n, f_n)(f_n),
$$
and according to Table 10, $v_\lambda (f_n, f_n)(f_n) \not= 0$.
Note that if $2 \leq k  = n$, then
$$
\Im \partial^{k+1,1}_{\fgl(n)} = V_{\lambda_2} \oplus V_{\lambda_5}.
\eqno(2.2.6)
$$
In fact, according to Table 11,
$$
C^{k+1,1}_{\fgl(n)} = V_{\lambda_2} \oplus V_{\lambda_5},
$$
  and we get (2.2.6) by (2.2.1), since $\fg_k = 0.$
Note that if $2 \leq k \leq n - 1$, then
$$
\Im \partial^{k+1,1}_{\fgl(n)} = V_{\lambda_2} \oplus V_{\lambda_4}
\oplus V_{\lambda_5}. \eqno (2.2.7)
$$
Indeed, according to Table 11,
if $2 \leq k = n - 1$, then
$$
C^{k+1,1}_{\fgl(n)} = V_{\lambda_2} \oplus
V_{\lambda_4} \oplus 2V_{\lambda_5}, \eqno (2.2.8)
$$
and if $2 \leq k \leq n - 2$, then
$$
C^{k+1,1}_{\fgl(n)} = V_{\lambda_2} \oplus
V_{\lambda_4} \oplus 2V_{\lambda_5} \oplus V_{\beta_5}. \eqno (2.2.9)
$$
Since by Table 12
$$
\eqalign{
&\fg_k = V_{\lambda_5} \hbox{ for } 2 \leq k = n - 1,\cr &\fg_k =
V_{\lambda_5} \oplus V_{\beta_5} \hbox{ for } 2 \leq k \leq n -
2,\cr }
$$
we get (2.2.7) by (2.2.1).

We will show now that  for $2 \leq k \leq n + 1$
$$
\Ker \partial^{k,2}_{\fgl(n)} = \Im \partial ^{k+1,1}_{\fgl(n)}. \eqno (2.2.10)
$$
Let
$$
\Ker \partial^{k,2}_{\fgl(n)} =\oplus_{\lambda} k_\lambda V_\lambda.
\eqno (2.2.11)
$$
Let $k = 2$. According to Table 10, if $n = 2$, then
$$C^{k,2}_{\fgl(n)} = V_{\lambda_1} \oplus 2V_{\lambda_2} \oplus
V_{\lambda_5}, \eqno (2.2.12)$$
and if $n \geq 3$, then
$$
C^{k,2}_{\fgl(n)} = V_{\lambda_1} \oplus 2V_{\lambda_2}
\oplus V_{\lambda_4} \oplus V_{\lambda_5}.
$$
Therefore, in  (2.2.11) $k_{\lambda_5} \leq 1$, and by (2.2.5)
$k_{\lambda_1} = 0,
k_{\lambda_2} \leq 1$. Note that
if $n = 2$, then $k_{\lambda_4} = 0$, and if $n \geq 3$, then $k_{\lambda_4}
\leq 1$.
Thus, by (2.2.6) and (2.2.7) we get (2.2.10).

Let $k \geq 3$. Then according to Table 10,
  if $k = n + 1$, then
$$
C^{k,2}_{\fgl(n)} = V_{\lambda_1} \oplus V_{\lambda_2}.
$$
Hence, by (2.2.5) $\Ker \partial^{k,2}_{\fgl(n)} = 0$. If $k = n$,
then
$$
C^{k,2}_{\fgl(n)} = V_{\lambda_1} \oplus 2V_{\lambda_2}
\oplus V_{\lambda_3} \oplus V_{\lambda_5}. \eqno (2.2.13)
$$
Therefore, in (2.2.11) $k_{\lambda_5} \leq 1$ and by (2.2.5) $k_{\lambda_1} =
  k_{\lambda_3} = 0, k_{\lambda_2} \leq 1$.
So from (2.2.6) we get (2.2.10).

Finally, if $k \leq n - 1$, then
$$
C^{k,2}_{\fgl(n)} = V_{\lambda_1}\oplus
2V_{\lambda_2}\oplus V_{\lambda_3}\oplus V_{\lambda_4} \oplus
V_{\lambda_5}.
$$
Therefore, in (2.2.11)  $k_{\lambda_4} \leq
1,\hbox{ } k_{\lambda_5} \leq 1$ and by (2.2.5) $k_{\lambda_1} =
k_{\lambda_3} = 0, k_{\lambda_2} \leq 1$. Thus, by (2.2.7) we get
(2.2.10). This proves part a) of  Theorem 2.1.

\ssec{2.3. Proof of part b) of Theorem 2.1.} Note that
$\fg_*(\fg_{-1},\hbox{ } \hat{\fg_0}) = \fsvect (0|n) =
\oplus_{i=-1}^{n-2}\fg_i$, where $\hat{\fg_0} = \fsl(n)$, and the
$\fsl(n)$-module $\fg_i$ is isomorphic to $U \hbox{ if }i = -1$
and to $V_{\delta_1-\delta_{n-k}-\delta_{n-k+1}-\ldots -\delta_n}
\hbox{ if } 0 \leq k \leq n - 2$.

Hence, for $k \geq
n + 1$ we have  $H^{k,2}_{\fsl(n)} = 0$ and for $1\leq k \leq n$
there exist the following Spencer cochain sequences:
$$
\begin{array}{ll}
   \fsl(n) \otimes U^*\buildrel {\partial^{2,1}_{\fsl(n)}} \over
   \longrightarrow U\otimes \Lambda^2U^*
   \buildrel{\partial^{1,2}_{\fsl(n)}} \over \longrightarrow 0 &(k =
   1),\\
   C^{k+1,1}_{\fsl(n)}\buildrel{\partial^{k+1, 1}_{\fsl(n)}}
   \over \longrightarrow C^{k,2}_{\fsl(n)}
   \buildrel{\partial^{k,2}_{\fsl(n)}}\over \longrightarrow
   C^{k-1,3}_{\fsl(n)}  &(2\leq k \leq n - 1),\\
   0\buildrel{\partial^{n+1,1}_{\fsl(n)}} \over \longrightarrow
   C^{n,2}_{\fsl(n)} \buildrel{\partial^{n,2}_{\fsl(n)}} \over
   \longrightarrow C^{n-1,3}_{\fsl(n)}&(k = n).
\end{array}
$$

First, we will
show that $H^{1,2}_{\fsl(n)} = 0$. In fact, since $\fg_1 \cong
V_{\delta_1-\delta_{n-1}-\delta_n}$, then by (2.2.1)
$$
\Im \partial^{2,1}_{\fsl(n)} \cong (\fsl(n)\otimes U^*)/
V_{\delta_1-\delta_{n-1}-\delta_n}.
$$
Since $\dim \fsl(n)\otimes U^* = (n^2 - 1)n$ and by the formula from
the Appendix
$$
\dim V_{\delta_1-\delta_{n-1}-\delta_n} = {n(n + 1)(n - 2)\over 2},
$$
then
$$
\eqalign{
&\dim \Im\partial^{2,1}_{\fsl(n)} = (n^2 - 1)n - {n(n + 1)(n -
2)\over 2} = {n^3 + n^2\over 2} = \cr &\dim U\otimes \Lambda^2U^*
= \dim \Ker \partial ^{1,2}_{\fsl(n)}.\cr }
$$
Next, we will prove
that for $2 \leq k \leq n - 1$ we have  $H^{k,2}_{\fsl(n)}
  = 0.$

\begin{Lemma} If $2 \leq k = n - 1$, then
$$C^{k+1,1}_{\fsl(n)} = V_{\lambda_4} \oplus V_{\lambda_5}. \eqno (2.3.1)$$
If $2 \leq k \leq n - 2$, then
$$
C^{k+1,1}_{\fsl(n)} = V_{\lambda_4} \oplus V_{\lambda_5} \oplus V_{\beta_5}.
\eqno(2.3.2)
$$
\end{Lemma}
\begin{proof} Let $c \in
C^{k+1,1}_{\fgl(n)}.$ Then $c\in C^{k+1,1}_{\fsl(n)}$ if and only
if
$$
c \in (\fsl(n) \otimes S^{k-1}U^*) \otimes U^*. \eqno(2.3.3)
$$
Let
$$
C^{k+1,1}_{\fsl(n)} = \oplus_{\beta} k_{\beta}V_{\beta}. \eqno (2.3.4)
$$
Using the decomposition of $C^{k+1,1}_{\fgl(n)}$ into direct sum
of irreducible $\fgl(n)$-modules given in (2.2.8) and (2.2.9), we
check condition (2.3.3) for the corresponding highest vectors:
let $\beta = \lambda_4 = \beta_1$, then up to a complex constant,
$$
v_{\beta_1}(f_n)(f_{j_1} \wedge f_{j_2} \wedge \ldots \wedge f_{j_{k-1}})
= f_1 \otimes \tilde {f}_{j_k}, \hbox { where } 2 \leq j_k \leq
n.
$$
Since $f_1 \otimes \tilde {f}_{j_k} \in \fsl(n)$, then in
(2.3.4) $k_{\lambda_4} = 1.$
Let $\beta = \lambda_2 = \beta_3$, then up to a  nonzero constant,
$$
v_{\beta_3}(f_n)(f_{n-k+2} \wedge f_{n-k+3} \wedge \ldots \wedge f_n)
=\sum_{j=1}^n f_j \otimes \tilde {f}_j.
$$
Since $\sum_{j=1}^n f_j
\otimes \tilde {f}_j \not \in \fsl(n)$, then $k_{\lambda_2} = 0.$
Let $\beta = \beta_5$. Then up to a constant,
$$
v_{\beta_5}(f_{j_1})(f_{j_2} \wedge \ldots \wedge f_{j_k}) = f_1 \otimes
\tilde {f}_{j_{k+1}}, \hbox{ where } 2 \leq j_{k+1} \leq n.
$$
Since $f_1 \otimes \tilde {f}_{j_{k+1}} \in \fsl(n)$, then
$k_{\beta_5} = 1.$ Let us show that $k_{\lambda_5} = 1.$ Indeed,
by the formula from the Appendix
$$
\dim V_{\delta_1-\delta{n-k+1}- \ldots -\delta_n} = {(n + 1)!\over
  (n - k + 1)k!(n - k -1)!}.
$$
Therefore,
$$
\dim C^{k+1,1}_{\fsl(n)} = \dim (V_{\delta_1-\delta{n-k+1}- \ldots -\delta_n}
\otimes U^*) = {n(n + 1)!\over (n - k + 1)k!(n - k -1)!}.
$$
Then by (2.2.2) and (2.2.3)
if $2 \leq k = n - 1$, then
$$
\dim C^{k+1,1}_{\fsl(n)} = \dim V_{\lambda_4} + dimV_{\lambda_5},
$$
if $2 \leq k \leq n - 2$, then by (2.2.2), (2.2.3), and (2.2.4)
$$
\dim C^{k+1,1}_{\fsl(n)} = \dim V_{\lambda_4} + dimV_{\lambda_5}
  + \dim V_{\beta_5},
$$
Thus, $k_{\lambda_5} = 1.$ This proves the Lemma.
\end{proof}

\ssbegin{2.4}{Lemma} $\Ker \partial ^{k,2}_{\fsl(n)}$ doesn't contain
irreducible $\fsl(n)$-submodules with highest weight $\lambda_2$.
\end{Lemma}
\begin{proof} Let $c \in C^{k,2}_{\fgl(n)}$. Then $c
\in C^{k,2}_{\fsl(n)}$ if and only if
$$
c \in (\fsl(n) \otimes S^{k-2}U^*) \otimes \lambda^2U^*.
$$
According to Table 10, each highest vector with weight $\lambda_2$
in $C^{k,2}_{\fgl(n)}$ is
$$
\eqalign {
&v = av_{\lambda_2} + bv_{\lambda_3}, \hbox { if } k = 2 \leq
n,\cr &v = av_{\lambda_2} + bv_{\lambda_6}, \hbox { if } 3 \leq k
\leq n - 1, \hbox { where } a,\hbox{ } b \in \Cee.\cr }
$$
Therefore,  up to a nonzero constant,
$$
v(f_n, f_n)(f_{n-k+2} \wedge \ldots \wedge f_{n-1}) =
(-1)^kaf_n \otimes \tilde {f}_n + b\sum_{j=1}^n f_j \otimes \tilde {f}_j.
$$
Note that
$$
(-1)^kaf_n\otimes \tilde {f}_n + b\sum_{j=1}^n f_j \otimes \tilde {f}_j
\in \fsl(n)
$$
if and only if  $(-1)^ka + nb = 0.$ But in this case
$v \not \in \Ker ^{k,2}_{\fsl(n)}.$ In fact,
$$
\partial^{k,2}_{\fsl(n)}(f_n, f_n, f_n) = -3v(f_n, f_n)(f_n) \not =0.
$$
This proves  Lemma 2.4.
\end{proof}

In order to prove that for $2 \leq k \leq n - 1$
$$
\Im \partial^{k+1,1}_{\fsl(n)} = \Ker \partial^{k,2}_{\fsl(n)}, \eqno (2.4.1)
$$
observe that if $2 \leq k = n - 1$, then $\fg_k = 0$ and if $2
\leq k \leq n - 2$, then $\fg_k = V_{\beta_5}.$ Thus, by (2.2.1),
(2.3.1), and (2.3.2)
$$
\Im \partial^{k+1,1}_{\fsl(n)} = V_{\lambda_4} \oplus V_{\lambda_5}.
$$
Therefore, by (2.2.7), (2.2.10), and Lemma 2.4, we get (2.4.1).

Finally, let $2 \leq k = n$. Notice that by (2.2.6) and (2.2.10)
$$
\Ker \partial^{n,2}_{\fgl(n)} = V_{\lambda_2} \oplus V_{\lambda_5}, \eqno
(2.4.2)
$$
where $V_{\lambda_5}$ is a trivial $\fsl(n)$-module.
By (2.2.12) and (2.2.13) the
multiplicity of $\lambda_5$ in $C^{n,2}_{\fgl(n)}$ is 1.
Moreover, this trivial submodule is contained in
$C^{n,2}_{\fsl(n)}$, because the $\fsl(n)$-module
$C^{n,2}_{\fsl(n)}$ is isomorphic to $\Lambda^2U \otimes
\Lambda^2U^*$, which contains a trivial $\fsl(n)$-submodule
(generated by
  $\sum_{i,j}f_if_j \otimes \tilde {f_i} \tilde {f_j}).$
Thus, $\Ker \partial^{n,2}_{\fsl(n)}$ must contain a trivial
submodule. According to (2.4.2) and Lemma 2.4, $\Ker
\partial^{n,2}_{\fsl(n)}$ is a trivial $\fsl(n)$-submodule,
generated by
$$
v_{\lambda_5} = \sum_{j=0}^{n-1} (-1)^{(n-1)j}
\sum_{i=1}^n f_i \otimes \tilde {f}_{s^j(1)} \wedge \tilde {f}_{s^j(2)}
\wedge \ldots \wedge \tilde {f}_{s^j(n-1)} \otimes \tilde {f}_{s^j(n)}
\tilde {f}_i,
$$
where $s$ is a cyclic permutation of $(1,\hbox{ } 2, \ldots ,\hbox{ } n)$.
This proves Theorem 2.1.

%%%%%%%%%%%%%%%%%%%%%%%%%%%%%%%%%%%%%%%%%%%%%%%%%%%%%%%%%%%%%%%%%%%%%%%%%%%%%%%%%%%%%%%%%%%%%%%%%%%%%%%%%%%%%%%%%
%\S3
%%%%%%%%%%%%%%%%%%%%%%%%%%%%%%%%%%%%%%%%%%%%%%%%%%%%%%%%%%%%%%%%%%%%%%%%%%%%%%%%%%%%%%%%%%%%%%%%%%%%%%%%%%%%%%%%%%%%%%%%%%
\section*{Penrose's tensors}

\ssbegin{3.1}{Theorem} If $m, n > 1$, then $H^{k,2}_{\fg_0} = 0$ for
$k > 2$ and the $\fg_0$-modules $H^{1,2}_{\fg_0}$ and
$H^{2,2}_{\fg_0}$ are the direct sums of irreducible submodules
whose highest weights are given in Table 13.

If $m =n$, then $H^{k,2}_{\hat{\fg_0}} = H^{k,2}_{\fg_0}$ for any $k$ and
if $m\not= n$, then $H^{1,2}_{\hat {\fg_0}} = H^{1,2}_{\fg_0}$
whereas $H^{2,2}_{\hat {\fg_0}} = H^{2,2}_{\fg_0} \oplus V_{\eps_1+
\eps_2-2\delta_n}\oplus V_{2\eps_1-\delta_{n-1}-\delta_n}$ if
either $m = 2$ or $n = 2$;

$H^{2,2}_{\hat {\fg_0}} = H^{2,2}_{\fg_0}\oplus
V_{\eps_1+\eps_2-2\delta_n}\oplus
V_{2\eps_1-\delta_{n-1}-\delta_n}\oplus
V_{\eps_1+\eps_2-\delta_{n-1}-\delta_n}$ if $m,\hbox{ } n
>  2.$
\end{Theorem}

\ssec{3.2. Calculation of $H^{1,2}_{\fg_0}$ and $H^{1,2}_{\hat {\fg_0}}$
for $m,\hbox{ } n \geq 2,\hbox{ } m\not= n$.}
For $k = 1$ the Spencer cochain sequence is of the form
$$
\fg_0\otimes \fg_{-1}^*\buildrel {\partial^{2,1}_{\fg_0}} \over \longrightarrow
\fg_{-1}\otimes \Lambda^2\fg_{-1}^*\buildrel
{\partial^{1,2}_{\fg_0}} \over \longrightarrow 0.
$$
Observe that
$$
\eqalign {
&\fg_{-1}\otimes \Lambda^2\fg_{-1}^* = (U\otimes V^*)\otimes
\Lambda ^2(U^*\otimes V)  \cong \cr &(U\otimes V^*)\otimes
(\Lambda^2U^*\otimes S^2V \oplus S^2U^*\otimes \Lambda^2V) \cong
\cr &(\Lambda^2U^*\otimes U)\otimes (S^2V\otimes V^*)\oplus
(S^2U^*\otimes U) \otimes (\Lambda^2V\otimes V^*),\cr
&\fg_0\otimes \fg_{-1}^* = (V\otimes V^*/\Cee \oplus U\otimes
U^*/\Cee \oplus \Cee) \otimes (U^*\otimes V),\cr &\fg_1 =
U^*\otimes V.\cr }
$$
Therefore, as $\fgl(m) \oplus\fgl(n)$-modules,
$$
\Im \partial^{2,1}_{\fg_0}\cong (V\otimes V^*/\Cee \oplus U\otimes U^*/\Cee)
\otimes (U^*\otimes V) \eqno (3.2.1)
$$
and
$$
H^{1,2}_{\fg_0}\cong (\Lambda^2U^*\otimes U/U^*)\otimes
(S^2V\otimes V^*/V) \oplus (S^2U^*\otimes U/U^*)\otimes
(\Lambda^2V\otimes V^*/V). \eqno (3.2.2)
$$
Note that
$$
\eqalign{
&S^2V\otimes V^*/V = V_{2\eps_1 - \eps_m},\cr &\Lambda^2V\otimes
V^*/V = V_{\eps_1+\eps_2-\eps_m } \hbox{ for } m>2,\cr
&\Lambda^2V\otimes V^*/V = 0 \hbox{ for } m=2. }
$$
Since $U$ is
purely odd, we deduce with the help of Table 5 of [OV] that
$$
\eqalign{
&\Lambda^2U^* \otimes U/U^* = V_{\delta_1-2\delta_n},\cr
&S^2U^*\otimes U/U^* = V_{\delta_1-\delta_{n-1}-\delta_n }
\hbox{ for } n>2,\cr
&S^2U^*\otimes U/U^* = 0
\hbox{ for } n = 2.\cr}
$$
Therefore, we have
$$
H^{1,2}_{\fg_0} = V_{2\eps_1-\eps_m+\delta_1-2\delta_n} \hbox{ if }
m = 2, \hbox{ }n > 2
$$
and
$$
H^{1,2}_{\fg_0} = V_{2\eps_1-\eps_m+\delta_1-2\delta_n}
\oplus V_{\eps_1+\eps_2-\eps_m+\delta_1-\delta_{n-1}-\delta_n}
\hbox{ if } m,\hbox{ }n > 2.
$$
By part b) of  Theorem  1.1
$\fg_*(\fg_{-1},\hbox{ }\hat{\fg_0}) = \fg_{-1} \oplus
\hat{\fg_0}$. Therefore, by (2.2.1), we have
$$
\Im \partial^{2,1}_{\hat{\fg_0}} = \hat{\fg_0} \otimes \fg_{-1}^* =
\Im \partial^{2,1}_{\fg_0}.
$$
Hence, $H^{1,2}_{\hat{\fg_0}} =
H^{1,2}_{\fg_0}.$

\ssec{3.3. Calculation of $H^{1,2}_{\hat{\fg_0}}$  for $m = n > 1$.} Since by
parts c) and d) of Theorem 1.1 the first term of the Cartan
prolongation $\fg_*(\fg_{-1}, \hat{\fg_0})$ is $U^*\otimes V$,
then by (2.2.1)
$$
\Im \partial^{2,1}_{\hat{\fg_0}} =
[(V\otimes V^*/\Cee \oplus U\otimes U^*/\Cee)\otimes (U^*\otimes
V)]/(U^*\otimes V).
$$
Therefore, by (3.2.1) and (3.2.2),
$$H^{1,2}_{\hat{\fg_0}} \cong (\Lambda^2U^*\otimes U/U^*)\otimes
(S^2V \otimes V^*/V)
\oplus(S^2U^*\otimes U/U^*) \otimes (\Lambda^2V\otimes V^*/V) \oplus
(U^*\otimes V).$$
Hence,
$$
H^{1,2}_{\hat{\fg_0}} = V_{2\eps_1-\eps_2+\delta_1-2\delta_2}
\oplus V_{\eps_1-\delta_2}  \hbox{ for } n = 2
$$
and
$$
H^{1,2}_{\hat{\fg_0}} = V_{2\eps_1-\eps_n+\delta_1-2\delta_n}
\oplus V_{\eps_1 + \eps_2-\eps_n+\delta_1-\delta_{n-1}-\delta_n}
\oplus V_{\eps_1-\delta_n}  \hbox{ for } n > 2.
$$

\ssec{3.4. Calculation of $H^{2,2}_{\fg_0}$ for $m,\hbox{
} n >1, \hbox{ } m\not=n$.} For $k = 2$ the Spencer cochain
sequence is of the form
$$
\fg_1\otimes \fg_{-1}^*\buildrel {\partial^{3,1}_{\fg_0}}\over \longrightarrow
\fg_0\otimes \Lambda^2\fg_{-1}^*
\buildrel{\partial^{2,2}_{\fg_0}}\over \longrightarrow
\fg_{-1}\otimes \Lambda^3\fg_{-1}^*.
$$
Observe that
$$
\eqalign{
&\fg_0\otimes \Lambda^2\fg_{-1}^* = (V\otimes V^*/\Cee\oplus
U\otimes U^*/\Cee \oplus \Cee) \otimes (\Lambda^2U^*\otimes
S^2V\oplus S^2U^*\otimes \Lambda^2V),\cr &\fg_1\otimes \fg_{-1}^*
= (U^*\otimes V)\otimes (U^*\otimes V),\cr &\fg_2 = 0.\cr }
$$
\begin{Lemma} As $\fgl(m) \oplus \fgl(n)$-module,
  $\fg_0 \otimes \Lambda^2\fg_{-1}^*$ is the direct sum of the
irreducible submodules whose
highest weights and highest vectors are listed in Table 14.
[$s$ and $t$ denote the cyclic permutations of $(1,\hbox{ } 2,\hbox{
} 3)$ and $(n-2,\hbox{ } n-1, \hbox{ }n)$,
respectively.]
\end{Lemma}
The proof follows from the formula given in the Appendix.

Let us show that if
$$
\eqalign{
&\lambda = 3\eps_1 - \eps_m - 2\delta_n,\hbox{ } 2\eps_1 + \eps_2
- \eps_m - 2\delta_n\hbox{ } (m > 2),\hbox{ } 2\eps_1 + \eps_2 -
\eps_m - \delta_{n-1} - \delta_n,\cr &2\eps_1 + \delta_1 -
3\delta_n,\hbox{ } 2\eps_1 + \delta_1 - \delta_{n-1} - 2\delta_n
\hbox{ }(n>2), \hbox{ or } \eps_1 + \eps_2 + \delta_1 -
\delta_{n-1} - 2\delta_n,\cr }
$$
then $v_\lambda \not\in \Ker
\partial^{2,2}_{\fg_0}$. Recall that if $v \in \fg_0\otimes
\Lambda^2\fg_{-1}^*$, then
$$
\partial^{2,2}_{\fg_0}v(g_1, g_2, g_3) = -v(g_1, g_2)g_3
- v(g_1, g_3)g_2 - v(g_2, g_3)g_1 \eqno (3.4.1)
$$
for any $g_1, g_2, g_3 \in \fg_{-1}$.

Let $\lambda = 3\eps_1 -
\eps_m - 2\delta_n$. Then
$$
\eqalign{
&\partial^{2,2}_{\fg_0}v_\lambda (f_n \otimes \tilde e_1,
f_n\otimes \tilde e_1, f_n\otimes \tilde e_1) = -3v_\lambda
(f_n\otimes \tilde e_1, f_n\otimes \tilde e_1)(f_n\otimes \tilde
e_1) = \cr &= 3A_{1,m}(f_n\otimes \tilde e_1) = -3f_n\otimes
\tilde e_m \not= 0.\cr }
$$
Let $\lambda = 2\eps_1 + \eps_2 -
\eps_m - 2\delta_n\hbox{ } (m > 2).$ Then
$$
\eqalign{
&\partial^{2,2}_{\fg_0}v_\lambda (f_n\otimes \tilde e_1,
f_n\otimes \tilde e_1, f_{n-1}\otimes \tilde e_2) =
-v_\lambda(f_n\otimes \tilde e_1, f_n\otimes \tilde
e_1)(f_{n-1}\otimes \tilde e_2) = \cr &= - A_{2,m}(f_{n-1} \otimes
\tilde e_2) = f_{n-1}\otimes \tilde e_m \not= 0.\cr }
$$
Let $\lambda = 2\eps_1 + \eps_2 - \eps_m -
\delta_{n-1} - \delta_n$. Then
$$
\eqalign{
&\partial^{2,2}_{\fg_0} v_\lambda (f_n\otimes \tilde e_1,
f_n\otimes \tilde e_1, f_{n-1}\otimes \tilde e_2) = -2v_\lambda
(f_n\otimes \tilde e_1, f_{n-1}\otimes \tilde e_2)(f_n\otimes
\tilde e_1) = \cr &= -A_{1,m}(f_n\otimes \tilde e_1) = f_n\otimes
\tilde e_m \not= 0.\cr }
$$
The proof of the fact that $v_\lambda
\not\in \Ker \partial^{2,2}_{\fg_0}$ for $\lambda = 2\eps_1 +
\delta_1 - 3\delta_n,\hbox{ } 2\eps_1 + \delta_1 - \delta_{n-1} -
2\delta_n \hbox{ }(n > 2)$, and $\eps_1 + \eps_2 + \delta_1 -
\delta_{n-1} - 2\delta_n$ is similar.

Let $\lambda =
\eps_1 + \eps_2 + \eps_3 - \eps_m - \delta_{n-1} - \delta_n \hbox{
}(m > 3).$ Let us show that if $n = 2$, then $v_\lambda \in \Ker
\partial ^{2,2}_{\fg_0}$ and if $n > 2$, then $v_\lambda \not\in
\Ker \partial ^{2,2}_{\fg_0}$. Indeed, if $n = 2$, then for $j =
0,\hbox{ } 1,\hbox{ } 2$ we have
$$
\eqalign{
&\partial^{2,2}_{\fg_0}v_\lambda(f_1\otimes \tilde e_{s^j(2)},
f_2\otimes \tilde e_{s^j(3)}, f_1\otimes \tilde e_{s^j(1)}) = \cr
&- v_\lambda(f_1\otimes \tilde e_{s^j(2)}, f_2\otimes \tilde
e_{s^j(3)}) (f_1 \otimes \tilde e_{s^j(1)})\cr &- v_\lambda
(f_2\otimes \tilde e_{s^j(3)}, f_1\otimes \tilde e_{s^j(1)})
(f_1\otimes \tilde e_{s^j(2)}) = \cr & A_{s^j(1),m}( f_1\otimes
\tilde e_{s^j(1)})/2 - A_{s^j(2),m}(f_1\otimes \tilde
e_{s^j(2)})/2 = \cr &-f_1\otimes \tilde e_m /2 + f_1\otimes \tilde
e_m/2 = 0,\cr }
$$
$$
\eqalign{
&\partial ^{2,2}_{\fg_0}v_\lambda (f_1\otimes \tilde e_{s^j(2)},
f_2\otimes \tilde e_{s^j(3)}, f_2\otimes \tilde e_{s^j(1)}) = \cr
&- v_\lambda (f_1\otimes \tilde e_{s^j(2)}, f_2\otimes \tilde
e_{s^j(3)}) (f_2\otimes \tilde e_{s^j(1)}) - \cr &- v_\lambda
(f_2\otimes \tilde e_{s^j(1)}, f_1\otimes \tilde e_{s^j(2)})
(f_2\otimes \tilde e_{s^j(3)}) = \cr &A_{s^j(1),m}(f_2\otimes
\tilde e_{s^j(1)})/2 - A_{s^j(3),m}(f_2\otimes \tilde
e_{s^j(3)})/2 = \cr &-f_2\otimes \tilde e_m/2 + f_2\otimes \tilde
e_m/2 = 0.\cr }
$$
Therefore, $v_\lambda \in \Ker
\partial^{2,2}_{\fg_0}$. If $n > 2$, then
$$
\eqalign{
&\partial^{2,2}_{\fg_0}v_\lambda (f_{n-1}\otimes \tilde e_2,
f_n\otimes \tilde e_3, f_1\otimes \tilde e_1) = - v_\lambda
(f_{n-1} \otimes \tilde e_2, f_n\otimes \tilde e_3) (f_1\otimes
\tilde e_1) = \cr &= A_{1,m}(f_1\otimes \tilde e_1)/2 =
-f_1\otimes \tilde e_m/2 \not= 0. }
$$
The proof of the fact that
if $\lambda = \eps_1 + \eps_2 + \delta_1
  - \delta_{n-2} - \delta_{n-1} - \delta_n \hbox{ }(n \geq 4)$, then
$v_\lambda \in \Ker \partial^{2,2}_{\fg_0}$ for $m = 2$ and
$v_\lambda \not\in \Ker \partial^{2,2}_{\fg_0}$ for $m > 2$ is
similar.

Finally, let us show that if
$$
\lambda = 2\eps_1 - 2\delta_n,\hbox{ } \eps_1 + \eps_2 - 2\delta_n,\hbox{ }
2\eps_1 - \delta_{n-1} - \delta_n, \hbox{ or } \eps_1 + \eps_2 -
\delta_{n-1} - \delta_n
$$
and $v_\lambda \in \Ker\partial ^{2,2}_{\fg_0}$, then
$v_\lambda \in \Im\partial^{3,1}_{\fg_0}$. Note that since $\fg_2 = 0$, then,
as $\fgl(m) \oplus \fgl(n)$-modules,
$$
\Im \partial^{3,1}_{\fg_0}
\cong \fg_1 \otimes \fg_{-1}^* = (U^*\otimes V)\otimes (U^*\otimes V).
$$
Therefore, by Table 9,
$$
\Im \partial^{3,1}_{\fg_0} = V_{2\eps_1 -2\delta_n} \oplus
V_{\eps_1+\eps_2-\delta_{n-1}-\delta_n}\oplus
V_{\eps_1+\eps_2-2\delta_n}\oplus
V_{2\eps_1-\delta_{n-1}-\delta_n}.
$$
Let $\lambda = \eps_1 +
\eps_2 - 2\delta_n$. By Table 14 $\fg_0\otimes
\Lambda^2\fg_{-1}^*$ contains two irreducible components with the
indicated highest weight, and one of the corresponding highest
vectors is $v_\lambda^1$. Observe that
$$
\eqalign {
&\partial^{2,2}_{\fg_0}v_\lambda^1(f_n\otimes \tilde e_1,
f_n\otimes \tilde e_1, f_{n-1}\otimes \tilde e_2) =
-v_\lambda^1(f_n\otimes \tilde e_1, f_n\otimes \tilde e_1)(f_{n-1}
\otimes \tilde e_2) = \cr &= -A_{2,1}(f_{n-1}\otimes \tilde e_2) =
f_{n-1}\otimes \tilde e_1 \not= 0.\cr }
$$
Therefore, $\Ker
\partial^{2,2}_{\fg_0}$ contains precisely one irreducible
submodule with highest weight $\eps_1 + \eps_2 - 2\delta_n$ and
this submodule belongs to
  $\Im \partial^{3,1}_{\fg_0}$.  Similarly, $\fg_0\otimes \Lambda^2\fg_{-1}^*$
contains two irreducible submodules with highest weight
  $2\eps_1 - \delta_{n-1} - \delta_n$, one of which belongs to $\Ker
\partial^{2,2}_{\fg_0}$ and, therefore, to $\Im \partial^{3,1}_{\fg_0}$.

Let $\lambda = 2\eps_1 - 2\delta_n$. Then by Table 14
any $\fgl(m)\oplus \fgl(n)$-highest vector of weight $\lambda$,
which belongs to $\fg_0\otimes \Lambda^2\fg_{-1}^*$, is
$$
v_\lambda = k_1v_\lambda^1 + k_2v_\lambda^2 + k_3v_\lambda^3,
\hbox{ where } k_1, k_2, k_3 \in \Cee.
$$
If $v_\lambda \in \Ker
\partial^{2,2}_{\fg_0}$, then the condition
$\partial^{2,2}_{\fg_0}v_\lambda(f_n\otimes \tilde e_1, f_n\otimes
\tilde e_1, f_n\otimes \tilde e_1) = 0$ implies
$$
k_1(m-1) - k_2(n-1) + k_3(m-n) = 0, \eqno (3.4.2)
$$
and the condition $\partial^{2,2}_{\fg_0}v_\lambda
(f_n\otimes\tilde e_2, f_n\otimes\tilde e_1, f_1\otimes \tilde
e_1) = 0$ implies that
$$
k_1m - k_2n = 0. \eqno (3.4.3)
$$
Thus, for $m \not= n$ we have
$$
k_2 = mk_1/n, k_3 = - k_1/n. \eqno (3.4.4)
$$
Therefore, $\Ker \partial^{2,2}_{\fg_0}$ contains precisely one
irreducible submodule with highest weight $2\eps_1 - 2\delta_n$
and this submodule belongs to $\Im\partial^{3,1}_{\fg_0}$.

Finally, let $\lambda = \eps_1 + \eps_2 -
\delta_{n-1} - \delta_n$. Then by Table 14 any highest vector with
weight $\lambda$, which belongs to $\fg_0\otimes
\Lambda^2\fg_{-1}^*$, is
$$
v_\lambda = k_1v_\lambda^1 + k_2v_\lambda^2 + k_3v_\lambda^3,
\hbox{ where } k_1, k_2, k_3 \in \Cee,
$$
and if $m = 2$, then $k_1
= 0$. Note that if $v_\lambda \in \Ker \partial^{2,2}_{\fg_0}$,
then
$$
\partial^{2,2}_{\fg_0}v_\lambda(f_{n-1}\otimes \tilde e_1,
f_n\otimes \tilde e_2, f_n\otimes \tilde e_1) = 0
$$
implies that
$$
k_1 + k_2 + k_3(n-m) = 0.\eqno (3.4.5)
$$
Thus, if $m = 2$, then
$$
k_2 = (2-n)k_3. \eqno (3.4.6)
$$
If $m,\hbox{ }n > 2$, then the condition
$$
\partial^{2,2}_{\fg_0}v_\lambda
(f_{n-1}\otimes \tilde e_m, f_n\otimes \tilde e_2, f_1\otimes \tilde e_1) = 0
\eqno (3.4.7)
$$
implies that $k_1 + k_2 = 0.$ Hence,
$$
k_2 = - k_1,\hbox{ } k_3 = 0. \eqno (3.4.8)
$$
Therefore, $\Ker \partial^{2,2}_{\fg_0}$ contains precisely one
highest vector of weight $\eps_1 + \eps_2 - \delta_{n-1} -
\delta_n$ which belongs to $\Im \partial^{3,1}_{\fg_0}$. Thus, we
have  the description of $H^{2,2}_{\fg_0}$ given in Table 13.

\ssec{3.5. Calculation of $H^{2,2}_{\hat{\fg_0}}$ for
$m,\hbox{ } n > 1,\hbox{ } m\not= n$.}
By part b) of Theorem 1.1
$$
\fg_*(\fg_{-1}, \hbox{ }\hat{\fg_0}) =
\fg_{-1} \oplus \hat{\fg_0}.
$$
Therefore, the Spencer cochain
sequence for $k = 2$ takes the form
$$
0\buildrel {\partial^{3,1}_{\hat{\fg_0}}}\over \longrightarrow
\hat{\fg_0}\otimes
\Lambda^2\fg_{-1}^*\buildrel{\partial^{2,2}_{\hat{\fg_0}}} \over
\longrightarrow \fg_{-1}\otimes \Lambda^3\fg_{-1}^*.
$$
Note that since $\fg_0 = \hat{\fg_0}\oplus \Cee$, then
$$
\fg_0\otimes \Lambda^2\fg_{-1}^* = \hat{\fg_0}\otimes \Lambda^2\fg_{-1}^*
\oplus V_{2\eps_1 - 2\delta_n} \oplus V_{\eps_1+\eps_2-
\delta_{n-1}-\delta_n}. \eqno(3.5.1)
$$
As we have shown in
sec.3.4, if $\lambda$ is one of the weights from Table 14, then an
irreducible module with highest weight $\lambda$ is contained in
the decomposition of $\Ker
\partial^{2,2}_{\fg_0}$ into irreducible $\fgl(m) \oplus
\fgl(n)$-modules if and only if
$$
\begin{array}{l}
   \lambda = \eps_1 + \eps_2 + \eps_3 - \eps_m - \delta_{n-1} -
\delta_n \hbox{ }(m > 3),\\
   \eps_1 + \eps_2 + \delta_1 - \delta_{n-2} - \delta_{n-1} - \delta_n
\hbox{ }(n > 3),\\
   \eps_1 + \eps_2 - 2\delta_n,\hbox{ } 2\eps_1 -
   \delta_{n-1} - \delta_n,\hbox{ } 2\eps_1 - 2\delta_n \hbox{ or }
   \eps_1 + \eps_2 - \delta_{n-1} - \delta_n
\end{array}
\eqno{(3.5.2)}
$$
and its multiplicity is 1. Therefore, by (3.5.1), if
$$
\eqalign{
&\lambda = \eps_1 + \eps_2 + \eps_3 - \eps_m - \delta_{n-1} -
\delta_n\hbox{ } (m > 3),\cr & \eps_1 + \eps_2 + \delta_1 -
\delta_{n-2} - \delta_{n-1} - \delta_n\hbox{ } (n > 3),\cr &\eps_1
+ \eps_2 - 2\delta_n, \hbox{ or } 2\eps_1 - \delta_{n-1} -
\delta_n,\cr }
$$
then the corresponding submodule is contained in $\Ker
\partial^{2,2}_{\hat{\fg_0}}$ as well.

Let $\lambda = 2\eps_1 - 2\delta_n \hbox{ and } v_\lambda \in \Ker
\partial^{2,2}_{\fg_0}$. Then  (3.4.4) where $k_3 = 0$, implies
that $v_\lambda \not\in \Ker \partial^{2,2}_{\hat{\fg_0}}$.

Let $\lambda = \eps_1 + \eps_2 - \delta_{n-1} -
\delta_n$, $v_\lambda \in \Ker \partial^{2,2}_{\fg_0}$. Then
(3.4.6) implies that $v_\lambda\not\in \Ker
\partial^{2,2}_{\hat{\fg_0}}$ for $m = 2$, and (3.4.8) implies that
$v_\lambda\in \Ker \partial^{2,2}_{\hat{\fg_0}}$ for $m,\hbox{ } n
>  2$. Thus, we have
$$
H^{2,2}_{\hat{\fg_0}} = H^{2,2}_{\fg_0} \oplus V_{\eps_1+\eps_2-
2\delta_n} \oplus V_{2\eps_1-\delta_{n-1}-\delta_n} \hbox{ if
either } m = 2 \hbox{ or } n = 2
$$
and
$$
H^{2,2}_{\hat{\fg_0}} = H^{2,2}_{\fg_0} \oplus V_{\eps_1+\eps_2-
2\delta_n} \oplus V_{2\eps_1-\delta_{n-1}-\delta_n}\oplus
V_{\eps_1+\eps_2-\delta_{n-1}-\delta_n} \hbox{ if } m,\hbox{ } n >
2.
$$
\ssec{3.6. Calculation of $H^{2,2}_{\hat{\fg_0}}$
for $m = n > 1$.} By parts c) and d) of Theorem 1.1 the
first term
  of $\fg_*(\fg_{-1},\hbox{ }\hat{\fg_0})$ is $U^*\otimes V$ and the
second one is $\Cee$ for
$n = 2$ and zero for $n > 2.$ By formula (2.2.1) we have
$$
\Im \partial^{3,1}_{\hat{\fg_0}} = (U^*\otimes V)\otimes (U^*\otimes V)
\hbox{ for } n > 2 \eqno (3.6.1)
$$
and
$$
\Im \partial^{3,1}_{\hat{\fg_0}} = (U^*\otimes V)\otimes (U^*\otimes V)/\Cee
\hbox{ for } n = 2.  \eqno(3.6.2)
$$
Therefore, by Table 9,
$$
\Im \partial^{3,1}_{\hat{\fg_0}} = V_{2\eps_1-2\delta_n} \oplus
V_{2\eps_1-\delta_{n-1}-\delta_n} \oplus
V_{\eps_1+\eps_2-2\delta_n} \hbox{ for } n= 2
$$
and
$$
\Im \partial^{3,1}_{\hat{\fg_0}} = V_{2\eps_1-2\delta_n} \oplus
V_{2\eps_1-\delta_{n-1}-\delta_n} \oplus
V_{\eps_1+\eps_2-2\delta_n} \oplus V_{\eps_1+\eps_2-
\delta_{n-1}-\delta_n} \hbox{ for } n > 2.
$$
Therefore, by (3.4.6) and (3.5.2),
$$
H^{2,2}_{\hat{\fg_0}} = 0 \hbox{ for } n = 2,\hbox{ } 3
$$
and
$$
H^{2,2}_{\hat{\fg_0}} = V_{\eps_1+\eps_2+\eps_3-\eps_n-
\delta_{n-1}-\delta_n }\oplus
V_{\eps_1+\eps_2+\delta_1-\delta_{n-2}- \delta_{n-1}-\delta_n}
\hbox{ for } n > 3.
$$

\ssec{3.7. Computation of $H^{3,2}_{\fg_0}$  for $m,\hbox{ } n > 1,\hbox{ } m
\not= n$.} For $k = 3$ the Spencer cochain sequence is of the form
$$
\fg_2\otimes \fg_{-1}^*\buildrel{\partial^{4,1}_{\fg_0}}\over \longrightarrow
\fg_1\otimes \Lambda^2\fg_{-1}^*
\buildrel{\partial^{3,2}_{\fg_0}}\over \longrightarrow\fg_0\otimes
\Lambda^3\fg_{-1}^*.
$$
Observe that
$$
\eqalign{
&\fg_1\otimes \Lambda^2\fg_{-1}^* = (U^*\otimes V)\otimes
\Lambda^2(U^*\otimes V) \cong \cr &(\Lambda^2U^*\otimes U^*)
\otimes (S^2V\otimes V)\oplus (S^2U^*\otimes U^*) \otimes
(\Lambda^2V\otimes V),\cr &\fg_2 = 0.\cr }
$$
By Table 5 from [OV]
$$
\eqalign{
&S^2V\otimes V = V_{3\eps_1}\oplus V_{2\eps_1+\eps_2},\cr
&\Lambda^2V\otimes V = V_{2\eps_1+\eps_2} \oplus
V_{\eps_1+\eps_2+\eps_3}\hbox{ for } m > 2,\cr &\Lambda^2V\otimes
V = V_{2\eps_1+\eps_2}\hbox{ for } m = 2.\cr }
$$
Since $U$ is purely odd,
$$
\eqalign{
&\Lambda^2U^*\otimes U^* = V_{-3\delta_n} \oplus
V_{-\delta_{n-1}-2\delta_n},\cr &S^2U^*\otimes U^* =
V_{-\delta_{n-1}-2\delta_n} \oplus
V_{-\delta_{n-2}-\delta_{n-1}-\delta_n}\hbox{ for } n > 2,\cr
&S^2U^*\otimes U^* = V_{-\delta_{n-1}-2\delta_n}\hbox{ for } n =
2.\cr }
$$
The above decompositions imply the following

\begin{Lemma} The $ \fgl(m) \oplus \fgl(n)$-module
$\fg_1\otimes \Lambda^2\fg_{-1}^* $ is the direct sum of
irreducible submodules whose highest weights and highest vectors
are listed in Table 15. [$s$ and $t$ denote the cyclic
permutations of $(1,\hbox{ } 2,\hbox{ } 3)$ and $(n-2,\hbox{ }
n-1, \hbox{ }n)$, respectively.]
\end{Lemma}

Let us show that $\Ker
\partial^{3,2}_{\fg_0} = 0.$ Let $\lambda = 3\eps_1 - 3\delta_n$. Then
$$
\partial^{3,2}_{\fg_0}v_\lambda(f_n\otimes\tilde e_1, f_n\otimes \tilde e_1,
f_n\otimes \tilde e_1) = 3B_{1,n}(f_n\otimes \tilde e_1) =
3(e_1\otimes \tilde e_1 + f_n\otimes \tilde f_n) \not= 0.
$$
Let $\lambda = 2\eps_1 + \eps_2 - 3\delta_n$. Then
$$
\partial^{3,2}_{\fg_0}v_\lambda(f_n\otimes \tilde e_1, f_n\otimes \tilde e_1,
f_n\otimes \tilde e_1) = -3B_{2,n}(f_n\otimes \tilde e_1) = -3e_2
\otimes \tilde e_1 \not= 0.
$$
Let $\lambda = 3\eps_1 - \delta_{n-1} - 2\delta_n$. Then
$$
\partial^{3,2}_{\fg_0}v_\lambda(f_n\otimes \tilde e_1, f_n\otimes \tilde e_1,
f_n\otimes \tilde e_1) = 3B_{1,n-1}(f_n\otimes \tilde e_1) =
3f_n\otimes \tilde f_{n-1} \not= 0.
$$
Let $\lambda = 2\eps_1 +
\eps_2 - \delta_{n-1} - 2\delta_n$. Since by Table 15
$\fg_1\otimes \Lambda^2\fg_{-1}^*$ contains two irreducible
submodules with highest weight $\lambda$, then any highest vector
of weight $\lambda$ in $\fg_1 \otimes \Lambda^2\fg_{-1}^*$ is of
the form
$$
v_\lambda = k_1 v_\lambda^1 + k_2 v_\lambda ^2,\hbox{ where }
k_1, k_2 \in \Cee.
$$
Let $v_\lambda \in \Ker
\partial^{3,2}_{\fg_0}$. If $m > 2$, then
$$
\partial^{3,2}_{\fg_0}v_\lambda (f_{n-1}\otimes \tilde e_2, f_n\otimes
\tilde e_1, f_n\otimes \tilde e_m) = -(1/2) k_2 B_{1,n} (f_n \otimes
\tilde e_m)
= -(1/2)k_2 e_1\otimes \tilde e_m = 0.
$$
Therefore, $k_2 = 0.$
Moreover,
$$
\partial^{3,2}_{\fg_0}v_\lambda (f_n\otimes \tilde e_1, f_n\otimes
\tilde e_1, f_{n-1}\otimes \tilde e_m) = k_1 B_{2,n-1}(f_{n-1}\otimes
\tilde e_m) =
k_1 e_2 \otimes \tilde e_m = 0.
$$
Hence, $k_1 = 0.$ If $n > 2$, then
$$
\partial^{3,2}_{\fg_0}v_\lambda(f_{n-1}\otimes \tilde e_2, f_n\otimes \tilde
e_1, f_1\otimes \tilde e_1) = -(1/2) k_2 B_{1,n} (f_1 \otimes \tilde e_1) =
-(1/2)k_2 f_1\otimes \tilde f_n = 0.
$$
Therefore, $k_2 = 0.$
Moreover,
$$
\partial ^{3,2}_{\fg_0}v_\lambda (f_n\otimes \tilde e_1, f_n\otimes
\tilde e_1, f_1\otimes \tilde e_2) = k_1 B_{2,n-1}(f_1 \otimes \tilde
e_2) = k_1 f_1\otimes \tilde f_{n-1} = 0.$$
Hence, $k_1 = 0.$

Let $\lambda = \eps_1 + \eps_2 + \eps_3 - \delta_{n-1} -2\delta_n$. Then
$$
\partial^{3,2}_{\fg_0}v_\lambda (f_{n-1}\otimes \tilde e_2, f_n\otimes \tilde
e_3, f_n\otimes \tilde e_2) = (1/2)B_{1,n}(f_n\otimes \tilde e_2)
= (1/2)e_1 \otimes \tilde e_2 \not= 0.
$$
Let $\lambda = 2\eps_1 +
\eps_2 - \delta_{n-2} - \delta _{n-1} - \delta_n$. Then
$$
\eqalign{
&\partial^{3,2}_{\fg_0}v_\lambda (f_{n-1} \otimes \tilde e_2,
f_n\otimes
  \tilde e_1, f_{n-1}\otimes \tilde e_1) =  \cr
&(1/2) B_{1,n-2}(f_{n-1}\otimes \tilde e_1) = (1/2) f_{n-1}\otimes
\tilde f_{n-2} \not= 0.\cr }
$$
Finally, let $\lambda = \eps_1 +
\eps_2 + \eps_3 - \delta _{n-2}
  - \delta_{n-1} - \delta_n$. Then
$$
\eqalign{
&\partial^{3,2}_{\fg_0}v_\lambda (f_{n-1} \otimes \tilde e_2,
f_n\otimes \tilde e_3, f_{n-2}\otimes \tilde e_1) = \cr
&(1/2)(B_{1,n-2}(f_{n-2}\otimes \tilde e_1) + B_{3,n}(f_n\otimes
\tilde e_3) + B_{2,n-1}(f_{n-1} \otimes \tilde e_2)) = \cr &(1/2)
(e_1 \otimes \tilde e_1 + f_{n-2} \otimes \tilde f_{n-2} +
e_3\otimes \tilde e_3 + f_n \otimes \tilde f_n + e_2\otimes \tilde
e_2 + f_{n-1} \otimes \tilde f_{n-1}) \not= 0.\cr }
$$
Thus, $H^{3,2}_{\fg_0} = 0.$

\ssec{3.8. Calculation of $H^{3,2}_{\hat{\fg_0}}$ for  $ m = n > 1$.}
By part d) of Theorem 1.1 for $n > 2$ the first term of the Cartan
prolongation of the pair $(\fg_{-1},\hbox{ } \hat{\fg_0})$ is
$U^*\otimes V$ and the second one is zero.  Therefore, by
arguments similar to those from sec.3.7 we get
$H^{3,2}_{\hat{\fg_0}} = 0.$

If $n = 2$, then by part
c) of Theorem 1.1 the first term of $\fg_*(\fg_{-1},\hbox{ } \hat
{\fg_0})$ is $U^*\otimes V$, the second one is the 1-dimensional
$\fgl(2) \oplus \fgl(2)$-module with highest weight $\eps_1 +
\eps_2 - \delta_1 - \delta_2$, and the third one is zero. Thus, by
(2.2.1),
$$
\Im \partial ^{4,1}_{\hat{\fg_0}} = V_{2\eps_1+ \eps_2-\delta_1-
2\delta_2}.
$$
By Table 15 $(U^*\otimes V)\otimes
\Lambda^2(U^*\otimes V)$ contains two irreducible $\fgl(2) \oplus
\fgl(2)$-modules with highest weight $\lambda = 2\eps_1 + \eps_2
- \delta_1 - 2\delta _2$ and one of the corresponding highest
vectors is $v_\lambda^2$. Since
$$
\partial^{3,2}_{\hat{\fg_0}}v_\lambda^2(f_1\otimes \tilde e_1, f_2\otimes
\tilde e_2, f_2\otimes \tilde e_2) = B_{1,2}(f_2\otimes \tilde
e_2) = e_1\otimes \tilde e_2 \not= 0,
$$
then $\Ker
\partial^{3,2}_{\hat{\fg_0}} = \Im \partial^{4,1}_{\hat{\fg_0}}$.
Thus, $H^{3,2}_{\hat{\fg_0}} = 0.$

\ssec{3.9.
Calculation of $H^{4,2}_{\hat{\fg_0}}$ for  $m
= n = 2$.} For $k = 4$ the Spencer cochain  sequence  is of the
form
$$
\fg_3 \otimes \fg_{-1}^* \buildrel{\partial^{5,1}_{\hat{\fg_0}}}\over
\longrightarrow \fg_2 \otimes \Lambda^2\fg_{-1}^*
\buildrel{\partial ^{4,2}_{\hat{\fg_0}}}\over \longrightarrow
\fg_1\otimes \Lambda ^3\fg_{-1}^*.
$$
By part c) of Theorem 1.1 the
second term of $\fg_*(\fg_{-1}, \hbox{ }\hat{\fg_0})$ is $\fg_2 =
V_{\eps_1+\eps_2-\delta_1-\delta_2} =\langle g\rangle$, the
1-dimensional $\fgl(2) \oplus \fgl(2)$-module,  and the third
one is zero. Since by Table 9
$$
\Lambda^2\fg_{-1}^* = V_{2\eps_1-2\delta_2} \oplus
V_{\eps_1+\eps_2-\delta_1-\delta_2},
$$
then
$$
\fg_2\otimes \Lambda^2\fg_{-1}^* = V_{3\eps_1+\eps_2-\delta_1-3
\delta_2} \oplus V_{2\eps_1+2\eps_2-2\delta_1-2\delta_2}.
$$
Let $\lambda = 3\eps_1 + \eps_2 - \delta_1 - 3\delta_2$. Then by
Table 9 $v_\lambda = g \otimes (\tilde f_2 \otimes e_1) \wedge
(\tilde f_2\otimes e_1).$ Let $v$ be an element from the basis of
$\fg_{-1}$ such that $g(v) \not=0.$ If $v = f_2 \otimes \tilde
e_1$, then
$$
\partial^{4,2}_{\hat{\fg_0}}(f_2 \otimes \tilde e_1, f_2 \otimes
\tilde e_1, v) =
-3v_{\lambda}(f_2 \otimes \tilde e_1, f_2\otimes \tilde e_1)(v) = 3g(v)\not= 0,
$$
and if $v \not= f_2\otimes \tilde e_1$, then
$$
\partial^{4,2}_{\hat{\fg_0}}v_{\lambda}(f_2\otimes \tilde e_1, f_2\otimes
\tilde e_1, v) = -v_{\lambda}(f_2 \otimes \tilde e_1, f_2 \otimes
\tilde e_1)(v) = g(v) \not= 0.
$$
Let $\lambda = 2\eps_1 + 2\eps_2
- 2\delta_1 - 2\delta_2.$ Then by Table 9
$$
\eqalign{
&v_{\lambda} = g \otimes ((\tilde f_2\otimes e_1) \wedge (\tilde
f_1 \otimes e_2) - (\tilde f_2 \otimes e_2 ) \wedge (\tilde f_1
\otimes e_1) - \cr &- (\tilde f_1 \otimes e_1) \wedge (\tilde f_2
\otimes e_2) + (\tilde f_1 \otimes e_2) \wedge (\tilde f_2 \otimes
e_1)).\cr }
$$
Let $v$ be an element of the basis of $\fg_{-1}$
such that $g(v) \not= 0.$ Then
$$
\eqalign{
&\partial ^{4,2}_{\hat{\fg_0}}v_{\lambda}(f_2 \otimes \tilde e_1,
f_1 \otimes \tilde e_2, v) = - 2v_{\lambda}(f_2 \otimes \tilde
e_1, f_1 \otimes \tilde e_2) (v) = g(v) \not = 0 \cr &\hbox{ if
either } v = f_2 \otimes \tilde e_1 \hbox{ or } v = f_1 \otimes
\tilde e_2,\cr }
$$
and
$$
\eqalign{
&\partial^{4,2}_{\hat{\fg_0}}v_{\lambda}(f_1\otimes \tilde e_1,
f_2 \otimes \tilde e_2, v) = -2v_{\lambda}(f_1 \otimes \tilde e_1,
f_2 \otimes \tilde e_2) (v) = - g(v) \not= 0\cr &\hbox{ if either
} v = f_2 \otimes \tilde e_2 \hbox{ or } v = f_1 \otimes \tilde
e_1.\cr }
$$
Therefore, $H^{4,2}_{\hat{\fg_0}} = 0.$

\ssec{3.10. Calculation of $H^{k,2}_{\fg_0}$ for
$m = n > 1, \hbox{ } k > 0$.}
\begin{Lemma}
$H^{k,2}_{\fg_0} = H^{k,2}_{\hat{\fg_0}}.$
\end{Lemma}

\begin{proof} Note that if $\fg_*(\fg_{-1},\hbox{ } \hat{\fg_0}) =
\fg_{-1} \oplus ( \oplus_{k\geq 0} \hat {\fg_k})$ is the Cartan
prolongation of the pair $(\fg_{-1},\hbox{ } \hat{\fg_0})$, then,
since $\fg_k = \hat {\fg_k}\oplus S^k(\fg_{-1}^*) \hbox{ } (k\geq
0)$, the Spencer cochain sequence is of the form
$$
\begin{array}{ll}
   (\hat{\fg_0} \oplus \Cee) \otimes
   \fg_{-1}^*\buildrel{\partial^{2,1}_{\fg_0}} \over \longrightarrow
   \fg_{-1} \otimes \Lambda^2 \fg_{-1}^*
   \buildrel{\partial^{1,2}_{\fg_0}} \over \longrightarrow 0  &\hbox{
   for } k = 1,\\

   (\hat {\fg}_{k-1} \oplus S^{k-1}(\fg_{-1}^*))\otimes \fg_{-1}^*
   \buildrel{\partial^{k+1,1}_{\fg_0}}
   \over \longrightarrow
   (\hat {\fg}_{k-2} \oplus S^{k-2}(\fg_{-1}^*)) \otimes \Lambda^2\fg_{-1}^*
   \buildrel{\partial^{k,2}_{\fg_0}}\over \longrightarrow&\\
   (\hat{\fg}_{k-3}\oplus S^{k-3}(\fg_{-1}^*)) \otimes
\Lambda^3\fg_{-1}^* &\hbox{ for } k > 1.
\end{array}
$$

Note that
since $\fg_*(\fg_{-1}, \fg_0) = S^*(\fg_{-1}^*) \subplus \fg_*
(\fg_{-1}, \hat{\fg_0})$, then the sequence
$$
\eqalign{
&S^{k-1}(\fg_{-1}^*) \otimes \fg_{-1}^*
\buildrel{\bar{\partial}^{k+1,1}_{\fg_0}} \over \longrightarrow
S^{k-2}(\fg_{-1}^*) \otimes \Lambda^2\fg_{-1}^*
\buildrel{\bar{\partial}^{k,2}_{\fg_0}} \over \longrightarrow
S^{k-3}(\fg_{-1}^*)\otimes \Lambda^3\fg_{-1}^* \hbox{ for } k\geq
1, \cr }
$$
where $\bar{\partial}^{k+1,1}_{\fg_0}$ and
$\bar{\partial}^{k,2}_{\fg_0}$ are the restrictions of the
operators $\partial^{k+1,1}_{\fg_0}$ and $\partial^{k,2}_{\fg_0}$
to $S^{k-1}(\fg_{-1}^*)\otimes \fg_{-1}^*$ and
$S^{k-2}(\fg_{-1}^*)\otimes \Lambda^2\fg_{-1}^*$, respectively,
and $S^k(\fg_{-1}^*) = 0$ for $k < 0$, is well-defined. Hence the
corresponding cohomology groups
$$
\bar{H}^{k,2}_{\fg_0} = \Ker \bar {\partial}^{k,2}_{\fg_0}/
\Im \bar{\partial}^{k+1,1}_{\fg_0}
$$
are well-defined and $H^{k,2}_{\fg_0} = H^{k,2}_{\hat{\fg_0}} \oplus
\bar{H}^{k,2}_{\fg_0}.$

Let us show that $\bar
{H}^{k,2}_{\fg_0} = 0 \hbox{ for } k>0.$ For $k = 1$ this is
obvious. Let $k = 2$. Since $S^{k-2}(\fg_{-1}^*) \otimes
\Lambda^2\fg_{-1}^* = \langle z\rangle\otimes \Lambda^2\fg_{-1}^*$, where $z$
is a generator of the center of  $\fgl(n|n)$, then
$$
\Ker \bar{\partial}^{k,2}_{\fg_0} \cong \Lambda^2 \fg_{-1}^* .
$$
By formula (2.2.1)
$$
\Im \bar{\partial}^{k+1,1}_{\fg_0} \cong \fg_{-1}^*\otimes
\fg_{-1}^*/S^2\fg_{-1}^* =
\Lambda^2\fg_{-1}^*.
$$
Therefore, $\bar{H}^{2,2}_{\fg_0} = 0.$ Let
$k = 3.$ Observe that
$$
S^2(U^*\otimes V)\otimes (U^*\otimes V) =
(S^2U^*\otimes U^*)\otimes (S^2V\otimes V)\oplus (\Lambda^2U^*\otimes U^*)
\otimes (\Lambda^2V\otimes V).
$$
By Table 5 from [OV] we get:
$$
\eqalign{
&S^2 V\otimes V = V_{3\eps_1} \oplus V_{2\eps_1+\eps_2},\hbox{ }
\Lambda^2V\otimes V = V_{2\eps_1+\eps_2} \hbox{ if } n = 2,\cr
&\Lambda^2V\otimes V = V_{2\eps_1+\eps_2}\oplus
V_{\eps_1+\eps_2+\eps_3} \hbox{ if } n > 2.\cr }
$$
Since $U$ is odd,
$$
\eqalign{
&\Lambda^2U^*\otimes U^* = V_{-3\delta_n} \oplus
V_{-\delta_{n-1}-2\delta_n},\cr
&S^2U^*\otimes U^* = V_{-\delta_{n-1}-2\delta_n}\hbox{ if } n = 2,\cr
&S^2U^*\otimes U^* = V_{-\delta_{n-1}-2\delta_n}
\oplus V_{-\delta_{n-2}-\delta_{n-1}-\delta_n} \hbox{ if } n > 2.\cr}
$$
Therefore,
$$
S^2(U^*\otimes V)\otimes (U^*\otimes V) =
V_{3\eps_1-\delta_{n-1}-2\delta_n} \oplus V_{2\eps_1+
\eps_2-3\delta_n}\oplus 2V_{2\eps_1+\eps_2-\delta_{n-1}-2\delta_n}
\hbox{ if } n = 2
$$
and
$$
\begin{array}{ll}
   S^2(U^*\otimes V)\otimes (U^*\otimes V) = &\\
   V_{3\eps_1-\delta_{n-1}-2\delta_n}\oplus V_{2\eps_1+
   \eps_2-3\delta_n}\oplus 2V_{2\eps_1+\eps_2-\delta_{n-1}-2\delta_n}
   \oplus V_{3\eps_1-\delta_{n-2}-\delta_{n-1}-\delta_n} \oplus &\\
   V_{\eps_1+\eps_2+\eps_3-3\delta_n} \oplus
   V_{2\eps_1+\eps_2-\delta_{n-2}-\delta_{n-1}-\delta_n} \oplus
   V_{\eps_1+\eps_2+\eps_3-\delta_{n-1}-2\delta_n} &\hbox{ if } n
   > 2.
\end{array}
$$

Moreover, we have
$$
\begin{array}{ll}
   S^3(U^*\otimes V) = V_{2\eps_1+\eps_2-\delta_1-2\delta_2} &\hbox{
   if } n = 2 \hbox{ and }\\
   S^3(U^*\otimes V) =
   V_{\eps_1+\eps_2+\eps_3-3\delta_n} \oplus
   V_{3\eps_1-\delta_{n-2}-\delta_{n-1}-\delta_n} \oplus
   V_{2\eps_1+\eps_2-\delta_{n-1}-2\delta_n}&\hbox{ if
   } n > 2.
\end{array}
$$
Thus, by (2.2.1)
$$
\begin{array}{ll}
   \Im\bar{\partial}^{4,1}_{\fg_0} = V_{3\eps_1-\delta_{n-1}-2\delta_n}
   \oplus V_{2\eps_1+\eps_2-3\delta_n}\oplus V_{2\eps_1+\eps_2-
   \delta_{n-1}-2\delta_n} &\hbox{ if } n = 2,\\
   \Im\bar{\partial}^{4,1}_{\fg_0} = V_{3\eps_1-\delta_{n-1}-2\delta_n}
   \oplus V_{2\eps_1+\eps_2-3\delta_n} \oplus V_{2\eps_1+\eps_2-
   \delta_{n-1}-2\delta_n}\oplus &\\
   \oplus V_{2\eps_1+\eps_2-\delta_{n-2}-\delta_{n-1}-\delta_n} \oplus
   V_{\eps_1+\eps_2+\eps_3-\delta_{n-1}-2\delta_n} &\hbox{ if } n  > 2.
\end{array}
$$

Finally, the decomposition of the $\fgl(n)\oplus
\fgl(n)$-module $\fg_{-1}^* \otimes \Lambda^2\fg_{-1}^*$ into the
direct sum of irreducible components is given in Table 15.
Checking the action of $\bar{\partial}^{3,2}_{\fg_0}$ on the
highest vectors  we get:
$$
\Im \bar{\partial}^{4,1}_{\fg_0} = \Ker \bar {\partial}^{3,2}_{\fg_0}.
$$
Note that for $k > 3$ the cohomology groups $\bar
{H}^{k,2}_{\fg_0}$ coincide with the Spencer cohomology groups
$H^{k-2,2}_{\fo(n^2)}$ corresponding to the Cartan prolongation
  $\fg_*(V(0|n^2),\hbox{ } \fo(n^2)) = \fh(0|n^2)$, where
$\fo(n^2)$ is the orthogonal Lie algebra and $V(0|n^2)$ is the
standard odd $\fo(n^2)$-module.
These groups are vanishing for $k > 3$ (see Theorem 1.3 of Chapter 3).

\end{proof}

%%%%%%%%%%%%%%%%%%%%%%%%%%%%%%%%%%%%%%%%%%%%%%%%%%%%%%%%%%%%%%%%%%%%%%%%%%%%%%%%%%%%%%%%%%%%%%%%%%%%%%%%%%%%%%%%%%%%%%%%%%%%%%%%%%%%%%%%%%%%%%%%%%%%%%%%%%%%%%%%%%%%%%%%%%%%%%%%%%%%%%%%%%%%%%
%Chapter 3
%%%%%%%%%%%%%%%%%%%%%%%%%%%%%%%%%%%%%%%%%%%%%%%%%%%%%%%%%%%%%%%%%%%%%%%%%%%%%%%%%%%%%%%%%%%%%%%%%%%%%%%%%%%%%%%%%%%%%%%%%%%%%%%%%%%%%%%%%%%%%%%%%%%%%%%%%%%%%%%%%%%%%%%%%%%%%%%%%%%%%%%%%%%%%%%%%%%%%%%
\chapter[The analogues of the Riemann--Weyl tensors]
{The analogues of the Riemann--Weyl tensors for
classical superspaces}

Recall that $\Zee$-grading of
depth 1 of a Lie (super)algebra $\fg$ is the $\Zee$-grading of the
form $\fg = \oplus_{i\geq {-1}}\fg_i$. All such $\Zee$-gradings of
simple finite-dimensional complex Lie superalgebras are listed in
[S2]. Denote by $V_\lambda$ the irreducible module over a Lie
superalgebra with highest weight $\lambda$ and an even highest
vector.

%%%%%%%%%%%%%%%%%%%%%%%%%%%%%%%%%%%%%%%%%%%%%%%%%%%%%%%%%%%%%%%%%%%%%%%%%%%%%%%%%%%%%%%%%%%%%%%%%%%%%%%%%%%
%\S1
%%%%%%%%%%%%%%%%%%%%%%%%%%%%%%%%%%%%%%%%%%%%%%%%%%%%%%%%%%%%%%%%%%%%%%%%%%%%%%%%%%%%%%%%%%%%%%%%%%%%%%%%%%%%%%%%%%%%%%%%%%%%%%%%

\section*{Spencer cohomology of $\fsl(m|n)$ and
$\fpsl(n|n)$}

\ssec{1.1. Description of  the $\Zee$-gradings of depth 1.}

Let $V(m-p|q)$ and $U(p|n-q)$ be the standard $\fsl(m-p|q)$ and
$\fsl(p|n-q)$- modules, respectively.

All $\Zee$-gradings of depth 1 of $\fg = \fsl(m|n)$ and
$\fpsl(n|n)$ are of the form $\fg = \fg_{-1} \oplus \fg_0 \oplus
\fg_1$, where $\fg_1 = \fg_{-1}^* = V(m-p|q) \otimes U(p|n-q)^*.$

A) For $\fsl(m|n)$, where $m\not=n$, there are the
following possible values of $\fg_0$ for the $\Zee$-gradings of
depth 1:

a) $\fc(\fsl(m) \oplus \fsl(n))$;

b) $\fc(\fsl(m|q) \oplus \fsl(n-q))$, if $p = 0, q
\not= 0, n-q \not= 0$;

c) $\fc(\fsl(m-p) \oplus
\fsl(p|n))$, if $q = 0, p \not= 0, m-p\not= 0$;

d) $\fc(\fsl(m-p|q) \oplus \fsl(p|n-q))$, if $p \not= 0, q\not= 0$.

B) For $\fsl(n|n)$ there are the following possible
values of $\fg_0$ for the $\Zee$-gradings of depth 1:

a) $\fc(\fsl(n) \oplus \fsl(n))$;

b) $\fc(\fsl(n|q)
\oplus \fsl(n-q))$, if $p = 0, q \not= 0, n-q \not= 0$;

c) $\fc(\fsl(n-p) \oplus \fsl(p|n))$, if $q = 0, p
\not= 0, n-p\not= 0$;

d) $\fc(\fsl(n-p|q) \oplus
\fsl(p|n-q))$, if $p \not= 0, q\not= 0$.

  C) The $\Zee$-gradings of $\fpsl(n|n)$ are similar to those of
$\fsl(n|n)$, only $\fg_0$ is centerless.

\ssbegin{1.2}{Theorem}[Cartan prolongations] For the cases of sec.1.1. we have:

A) $\fg = \fsl(m|n)$, where $m\not=n$. Then
$\fg_*(\fg_{-1},\hbox{ } \fg_0) = \fg$, except for the following
cases:

a) if $n = 1$, then $\fg_* = \fvect(0|m)$, if $m = 1$, then $\fg_* =
\fvect(0|n)$;

b) if $n - q = 1$, then $\fg_* = \fvect(q|m)$, if $m
= 0, q = 1$, then $\fg_* = \fvect(n-1|0)$;

c) if $m - p = 1$, then $\fg_* = \fvect(p|n)$, if $n = 0, p =
1$, then $\fg_* = \fvect(m-1|0)$;

d) if $n - q = 0, p = 1$, then $\fg_* = \fvect(m-1|n)$, if $m - p = 0, q
= 1$, then $\fg_* = \fvect(n-1|m)$.

B) $\fg =\fsl(n|n)$. Then $\fg_*(\fg_{-1},\hbox{ } \fg_0) =
S^*(\fg_{-1}^*)\subplus \fpsl(n|n)$, except for the following cases:

a) if $n = 2$, then $\fg_* = S^*(\fg_{-1}^*)\subplus
\fh(0|4)$;

b) if $n - q = 1$, then $\fg_* =
\fvect(q|n)$;

c) if $n - p = 1$, then $\fg_* =\fvect(p|n)$;

d) if $n - q = 0, p = 1$ or $n - p = 0, q = 1$ then $\fg_* = \fvect(n-1|n)$.

C) $\fg =\fpsl(n|n)$. Then $\fg_*(\fg_{-1},\hbox{ } \fg_0) = \fg$, except
for the following cases:

a) if $n = 2$, then $\fg_* =\fh(0|4)$;

b) if $n - q = 1$, then $\fg_* = \fsvect(q|n)$;

c) if $n - p = 1$, then $\fg_* = \fsvect(p|n)$;

d) if $n - q = 0, p = 1$ or $n - p = 0, q = 1$, then $\fg_* = \fsvect(n-1|n)$.
\end{Theorem}

Let $\fg = \fsl(m|n)$, where $m\not=n$, or $\fpsl(n|n)$. We will
describe the Spencer cohomology groups for all $\Zee$-gradings of
depth 1 listed in sec.1.1.

First consider the cases easiest to formulate. Let $\langle \pi
_i\rangle$ be the $i$-th fundamental
weight of $\fg_0$.

\ssbegin{1.3}{Theorem}
\emph{1)} For $\fg_* = \fvect(m|n)$, $\fsvect(m|n)$ SFs vanish except for
$\fsvect(0|n)$, when SFs are of order $n$ and constitute the
$\fg_0$-module $\Pi ^n(\langle 1\rangle)$.

\emph{2)} For $\fg_*$ of series $\fh(0|n)$, nonzero SFs are
of order 1. For $n > 3$ SFs constitute $\fg_0$-module $\Pi
(V_{3\pi _1} \oplus V_{\pi _1})$.

\emph{3)} For $\fg_* =\fsh(0|n)$, nonzero SFs are the same as for
$\fh(0|n)$ and
additionally $\Pi^{n-1}(V_{\pi _1})$ of order $n - 1$.
\end{Theorem}

Consider  the $\Zee$-gradings of depth 1 of $\fg = \fsl(m|n)$
($m\not=n$) and $\fpsl(n|n)$ listed in sec.1.1 for which
$\fg_*(\fg_{-1}, \hbox{ }\fg_0) = \fg$. Describe the corresponding
SFs.

Case a) was discussed in Chapter 2. Consider case b).

\ssbegin{1.4}{Theorem} The nonzero SFs are of
orders 1 and 2. The $\fg_0$-module $H^{2,2}_{\fg_0}$ splits into
the direct sum of irreducible components whose weights are given
in Table 16. Table 16 also contains the highest weights (with
respect to the bases $\eps_1, \ldots ,\eps_{m+q}$ and $\delta_1,
\ldots , \delta_{n-q}$ of the dual spaces to the maximal tori of
$\fsl(m|q)$ and $\fsl(n-q)$, respectively) of irreducible
components of $H^{1,2}_{\fg_0}$ for the cases when
$H^{1,2}_{\fg_0}$ does split into the direct sum of irreducible
$\fg_0$-modules.
\end{Theorem}

Exceptional cases are as follows: if $m = q - 1, m > 1, n-q \geq 3$, then
$H^{1,2}_{\fg_0} =
V_{\eps_1+\eps_2-\eps_{m+q}+\delta_1-\delta_{n-q-1} -\delta_{n-q}}
\oplus X$, where $X$ is given by the nonsplit exact
sequence of $\fg_0$-modules
$$
0\longrightarrow V_{2\eps_1 - \eps_{m+q}+\delta_1-2\delta{n-q}}
\longrightarrow X  \longrightarrow \Pi (V_{\eps_1
+\delta_1-2\delta_{n-q}}) \longrightarrow 0;\eqno(1)
$$
if $m = q -
1, m > 1, n-q = 2$, then $H^{1,2}_{\fg_0} = X$, where $X$ is given
by (1);

if $m = 1, q = 2, n-q \geq 3$, then
$H^{1,2}_{\fg_0} =
\Pi(V_{\eps_1+\eps_2-\eps_3+\delta_1-\delta_{n-q-1}-\delta_{n-q}})
\oplus X$, where $X$ is given by the nonsplit exact sequence of
$\fg_0$-modules
$$
\begin{array}{l}
   0\longrightarrow V_{2\eps_1 - \eps_3+\delta_1-2\delta_{n-q}}
   \oplus \Pi (V_{-\eps_1+\eps_2 +\eps_3+\delta_1-2\delta_{n-q}})\\
   \longrightarrow X \longrightarrow \Pi (V_{\eps_1
   +\delta_1-2\delta_{n-q}}) \longrightarrow 0;
\end{array}
\eqno{(2)}
$$
if $m = 1,
q = 2, n-q = 2$, then $H^{1,2}_{\fg_0} = X$, where $X$ is given by
(2);
if $m = q + 1, n - q \geq 3$,
then $H^{1,2}_{\fg_0} = V_{2\eps_1-\eps_{m+q}+\delta_1
-2\delta_{n-q}} \oplus X$, where $X$ is given by the nonsplit
exact sequence of $\fg_0$-modules
$$
\begin{array}{ll}
   0\longrightarrow

V_{\eps_1+\eps_2-\eps_{m+q}+\delta_1-\delta_{n-q-1}-\delta_{n-q}}\longrightarrow
X&\\
   \longrightarrow\Pi
   (V_{\eps_1+\delta_1-\delta_{n-q-1}-\delta_{n-q}})\longrightarrow
   0&(q \geq 2),
\end{array}
$$

$$
\begin{array}{ll}
   0\longrightarrow V_{\eps_1-\eps_2+\eps_3+\delta_1-\delta_{n-2}
   -\delta_{n-1}}\oplus V_{\eps_1+\eps_2-\eps_3+\delta_1-\delta_{n-2}
   -\delta_{n-1}}\longrightarrow X&\\
   \longrightarrow \Pi
   (V_{\eps_1+\delta_1-\delta_{n-2} -\delta_{n-1}})\longrightarrow 0
   &(q = 1).
\end{array}
$$

Case c) is similar to  case b). Consider case d).

\ssbegin{1.5}{Theorem} The nonzero SFs are of orders 1 and 2. The
$\fg_0$-module $H^{2,2}_{\fg_0}$ splits into the direct sum of
irreducible components whose weights are given in Table 17. Table
17 also contains the highest weights (with respect to the bases
$\eps_1, \ldots ,\eps_{m-p+q}$ and $\delta_1, \ldots ,
\delta_{p+n-q}$ of the dual spaces to the maximal tori of
$\fsl(m-p|q)$ and $\fsl(p|n-q)$, respectively) of
irreducible components of $H^{1,2}_{\fg_0}$ for the cases when
$H^{1,2}_{\fg_0}$ does split into the direct sum of irreducible
$\fg_0$-modules.
\end{Theorem}

Exceptional cases are $m = p + q \pm
1$ and $n = p + q \pm 1$. More precisely: if $m = p +
q + 1$, $n \not= p + q \pm 1$, $q$, then $H^{1,2}_{\fg_0} =
V_{2\eps_1-\eps_{m-p+q}+\delta_1-2\delta_{p+n-q}} \oplus Y$, where
$Y$ is given by the nonsplit exact sequence of $\fg_0$-modules
$$
\begin{array}{ll}
   0\longrightarrow V_{\eps_1+\eps_2-\eps_3+\delta_1-\delta_{p+n-2}
   -\delta_{p+n-1}}\oplus

V_{\eps_1-\eps_2+\eps_3+\delta_1-\delta_{p+n-2}-\delta_{p+n-1}}\longrightarrow
Y&\\
   \longrightarrow V_{\eps_1+\delta_1-\delta_{p+n-2}
   -\delta_{p+n-1}}\longrightarrow 0 &(q = 1),
\end{array}
$$
$$
\begin{array}{ll}
   0\longrightarrow

V_{\eps_1+\eps_2-\eps_{m-p+q}+\delta_1-\delta_{p+n-q-1}-\delta_{p+n-q}}\longrightarrow
Y&\\
   \longrightarrow V_{\eps_1+\delta_1-
   \delta_{p+n-q-1}-\delta_{p+n-q}}\longrightarrow 0 &(q \geq 2);
\end{array}
$$
if $m = p + q + 1$, $n = p + q - 1$, then $H^{1,2}_{\fg_0} = X \oplus
Y$ , where $X$ is given by the nonsplit exact sequence of
$\fg_0$-modules
$$
\begin{array}{ll}
   0\longrightarrow
   V_{2\eps_1-\eps_{m-2+q}+\delta_1-2\delta_3}\oplus
   V_{2\eps_1-\eps_{m-2+q}-\delta_1-\delta_2+\delta_3}
   \longrightarrow&\\
   X\longrightarrow
   V_{2\eps_1-\eps_{m-2+q}-\delta_3} \longrightarrow 0 &(p = 2),
\end{array}
$$

$$
\begin{array}{ll}
   0\longrightarrow
   V_{2\eps_1-\eps_{m-p+q}+\delta_1-2\delta_{p+n-q}} \longrightarrow X&\\
   \longrightarrow V_{2\eps_1-\eps_{m-p+q}-\delta_{p+n-q}}
   \longrightarrow 0 &(p \geq 3),
\end{array}
$$
and $Y$ is given by the
nonsplit exact sequence of $\fg_0$-modules
$$
\begin{array}{ll}
   0\longrightarrow
   V_{\eps_1+\eps_2-\eps_3+\delta_1-\delta_{p+n-2}-\delta_{p+n-1}}\oplus

V_{\eps_1-\eps_2+\eps_3+\delta_1-\delta_{p+n-2}-\delta_{p+n-1}}\longrightarrow
&\\
   Y\longrightarrow
   V_{\eps_1+\delta_1-\delta_{p+n-2}-\delta_{p+n-1}}\longrightarrow 0
   & (q = 1),
\end{array}
$$
$$
\begin{array}{ll}
   0\longrightarrow
   V_{\eps_1+\eps_2-\eps_{m-p+q}+\delta_1-\delta_{p+n-q-1}-\delta_{p+n-q}}
   \longrightarrow Y &\\
   \longrightarrow
   V_{\eps_1+\delta_1-\delta_{p+n-q-1}-\delta_{p+n-q}}\longrightarrow
   0 &(q \geq 2);
\end{array}
$$
if $m = p + q + 1$, $n =  q $, then
$H^{1,2}_{\fg_0} =
V_{2\eps_1-\eps_{m-p+q}+\delta_1-\delta_{p-1}-\delta_p }\oplus Y$
$(p \geq 3)$ or $H^{1,2}_{\fg_0} = Y$ $(p = 2)$, where $Y$ is
given by the nonsplit exact sequence of $\fg_0$-modules
$$
0\longrightarrow V_{\eps_1+\eps_2-\eps_3+\delta_1-2\delta_p}
\oplus V_{\eps_1-\eps_2+\eps_3+\delta_1-2\delta_p} \longrightarrow
Y\longrightarrow V_{\eps_1+\delta_1-2\delta_p} \longrightarrow 0
\eqno (q = 1),
$$
$$
0\longrightarrow V_{\eps_1+\eps_2-\eps_{m-p+q}+\delta_1-2\delta_p}
\longrightarrow Y\longrightarrow V_{\eps_1+\delta_1-2\delta_p}
\longrightarrow 0 \eqno (q \geq 2);
$$
if $n = p + q + 1$, $m\not=
p + q \pm 1$, $p$, then $H^{1,2}_{\fg_0} =
V_{2\eps_1-\eps_{m-p+q}+\delta_1-2\delta_{p+n-q}} \oplus Y$, where
$$
\begin{array}{ll}
   0\longrightarrow V_{\eps_1+\eps_2-\eps_{m-1+q}+\delta_1-\delta_2-
   \delta_3}\oplus V_{\eps_1+\eps_2-\eps_{m-1+q}-\delta_1+\delta_2-
   \delta_3} \longrightarrow &\\
   Y\longrightarrow
   V_{\eps_1+\eps_2-\eps_{m-1+q}- \delta_3} \longrightarrow 0 & (p =1),
\end{array}
$$
$$
\begin{array}{ll}
   0\longrightarrow
   V_{\eps_1+\eps_2-\eps_{m-p+q}+\delta_1-\delta_{p+n-q-1}-\delta_{p+n-q}}
   \longrightarrow Y &\\
   \longrightarrow
   V_{\eps_1+\eps_2-\eps_{m-p+q}- \delta_{p+n-q}} \longrightarrow 0 &
   (p \geq 2);
\end{array}
$$
if $n = p + q + 1$, $m = p + q - 1 $, then
$H^{1,2}_{\fg_0} = X \oplus Y$, where
$$
\begin{array}{ll}
   0\longrightarrow V_{2\eps_1-\eps_3+\delta_1-2\delta_{p+n-2}}\oplus
   V_{-\eps_1+\eps_2+\eps_3+\delta_1-2\delta_{p+n-2}}\longrightarrow X&\\
   \longrightarrow V_{\eps_1+\delta_1-2\delta_{p+n-2}}
   \longrightarrow 0 & (q = 2),
\end{array}
$$
$$
0\longrightarrow V_{2\eps_1-\eps_{m-p+q}+\delta_1-2\delta_{p+n-q}}
\longrightarrow X \longrightarrow
V_{\eps_1+\delta_1-2\delta_{p+n-q}} \longrightarrow 0 \eqno (q
\geq 3),
$$
and $Y$ is given by the nonsplit exact sequence of
$\fg_0$-modules
$$
\begin{array}{ll}
   0\longrightarrow
   V_{\eps_1+\eps_2-\eps_{m-1+q}+\delta_1-\delta_2-\delta_3} \oplus
   V_{\eps_1+\eps_2-\eps_{m-1+q}-\delta_1+\delta_2-\delta_3}&\\
   \longrightarrow Y \longrightarrow
   V_{\eps_1+\eps_2-\eps_{m-1+q}-\delta_3} \longrightarrow 0 &(p = 1),
\end{array}
$$
$$
\begin{array}{ll}
   0\longrightarrow
   V_{\eps_1+\eps_2-\eps_{m-p+q}+\delta_1-\delta_{p+n-q-1}-\delta_{p+n-q}}
   \longrightarrow Y &\\
   \longrightarrow
   V_{\eps_1+\eps_2-\eps_{m-p+q}-\delta_{p+n-q}} \longrightarrow 0
   &(p \geq 2);
\end{array}
$$
if $n = p + q + 1$, $m = p $, then $H^{1,2}_{\fg_0} =
V_{\eps_1+\eps_2-\eps_q+\delta_1-2\delta_{p+n-q} }\oplus Y$ $(q
\geq 3)$, and $H^{1,2}_{\fg_0} = Y$ $(q = 2)$, where
$$
0\longrightarrow V_{2\eps_1-\eps_q+\delta_1-\delta_2-\delta_3}\oplus
V_{2\eps_1-\eps_q-\delta_1+\delta_2-\delta_3}\longrightarrow
Y\longrightarrow V_{2\eps_1-\eps_q-\delta_3}
\longrightarrow 0 \eqno (p = 1),
$$
$$
0\longrightarrow
V_{2\eps_1-\eps_q+\delta_1-\delta_{p+n-q-1}-\delta_{p+n-q}}\longrightarrow
Y\longrightarrow V_{2\eps_1-\eps_q-\delta_{p+n-q}}
\longrightarrow 0 \eqno (p \geq 2);
$$
if $m = p + q - 1$, $n \not=
p + q \pm 1$, $q$, then $H^{1,2}_{\fg_0} = X\oplus
V_{\eps_1+\eps_2-\eps_{m-p+q}+\delta_1-
\delta_{p+n-q-1}-\delta_{p+n-q}}$, where $X$ is given by the
nonsplit exact sequence of $\fg_0$-modules
$$
\begin{array}{ll}
   0\longrightarrow V_{2\eps_1-\eps_3+\delta_1-2\delta_{p+n-2}}
   \oplus V_{-\eps_1+\eps_2+\eps_3+\delta_1-2\delta_{p+n-2}}
   \longrightarrow &\\
   X\longrightarrow
   V_{\eps_1+\delta_1-2\delta_{p+n-2}} \longrightarrow 0 & (q = 2),
\end{array}
$$
$$
0\longrightarrow V_{2\eps_1-\eps_{m-p+q}+\delta_1-2\delta_{p+n-q}}
\longrightarrow X\longrightarrow
V_{\eps_1+\delta_1-2\delta_{p+n-q}} \longrightarrow 0 \eqno (q
\geq 3);
$$
if $m = p + q - 1$, $n = q $, then $H^{1,2}_{\fg_0} =
X\oplus V_{\eps_1+\eps_2-\eps_{m-p+q}+\delta_1- 2\delta_p}$
$(p\geq 3)$, where $X$ is given by the nonsplit exact sequence of
$\fg_0$-modules
$$
\begin{array}{ll}
   0\longrightarrow
   V_{2\eps_1-\eps_3+\delta_1-\delta_{p-1}-\delta_p} \oplus
   V_{-\eps_1+\eps_2+\eps_3+\delta_1-\delta_{p-1}-\delta_p}
   \longrightarrow &\\
   X\longrightarrow
   V_{\eps_1+\delta_1-\delta_{p-1}-\delta_p} \longrightarrow 0 &(q =2),
\end{array}
$$
$$
0\longrightarrow V_{2\eps_1-\eps_{m-p+q}+\delta_1-\delta_{p-1}-\delta_p}
\longrightarrow X\longrightarrow V_{\eps_1+\delta_1-\delta_{p-1}-\delta_p}
\longrightarrow 0 \eqno (q \geq 3);
$$
if $n = p + q - 1$, $m\not =
p + q \pm 1$, $p$, then $H^{1,2}_{\fg_0} = X \oplus
V_{\eps_1+\eps_2-\eps_{m-p+q}+\delta_1-
\delta_{p+n-q-1}-\delta_{p+n-q}}$, where $X$ is given by the
nonsplit exact sequence of $\fg_0$-modules
$$
\begin{array}{ll}
   0\longrightarrow V_{2\eps_1-\eps_{m-2+q}+\delta_1-2\delta_3}
   \oplus V_{2\eps_1-\eps_{m-2+q}-\delta_1-\delta_2+\delta_3}
   \longrightarrow &\\
   X\longrightarrow
   V_{2\eps_1-\eps_{m-2+q}-\delta_3} \longrightarrow 0 & (p = 2),
\end{array}
$$
$$
0\longrightarrow V_{2\eps_1-\eps_{m-p+q}+\delta_1-2\delta_{p+n-q}}
\longrightarrow X\longrightarrow
V_{2\eps_1-\eps_{m-p+q}-\delta_{p+n-q}} \longrightarrow 0\eqno (p
\geq 3);
$$
if $n = p + q - 1$, $m = p$, then $H^{1,2}_{\fg_0} = X
\oplus V_{2\eps_1-\eps_q+\delta_1-
\delta_{p+n-q-1}-\delta_{p+n-q}}$ $(q\geq 3)$, where $X$ is given
by the nonsplit exact sequence of $\fg_0$-modules
$$
\begin{array}{ll}
   0\longrightarrow V_{\eps_1+\eps_2-\eps_q+\delta_1-2\delta_3}
   \oplus V_{\eps_1+\eps_2-\eps_q-\delta_1-\delta_2+\delta_3}
   \longrightarrow &\\
   X\longrightarrow
   V_{\eps_1+\eps_2-\eps_q-\delta_3} \longrightarrow 0& (p = 2),
\end{array}
$$
$$
0\longrightarrow V_{\eps_1+\eps_2-\eps_q+\delta_1-2\delta_{p+n-q}}
\longrightarrow X\longrightarrow
V_{\eps_1+\eps_2-\eps_q-\delta_{p+n-q}} \longrightarrow 0 \eqno (p
\geq 3);
$$
if $m = n = p + q + 1$, then $H^{1,2}_{\fg_0} = X
\oplus V_{2\eps_1-\eps_{m-p+q}+\delta_1- 2\delta_{p+n-q}}\oplus
V_{\eps_1- \delta_{p+n-q}}$, where $X$ is given by the nonsplit
exact sequence of $\fg_0$-modules
$$
0\longrightarrow Y\longrightarrow X\longrightarrow
V_{\eps_1-\delta_{p+n-q}}\longrightarrow 0,
$$
and $Y$ is given by
the nonsplit exact sequence of $\fg_0$-modules
$$
\begin{array}{ll}
   0\longrightarrow

V_{\eps_1+\eps_2-\eps_{m-p+q}+\delta_1-\delta_{p+n-q-1}-\delta{p+n-q}}\longrightarrow
   Y \longrightarrow
   V_{\eps_1+\delta_1-\delta_{p+n-q-1}-\delta{p+n-q}}\oplus&\\
   V_{\eps_1+\eps_2-\eps_{m-p+q}-\delta{p+n-q}} \longrightarrow 0 &
   (p\geq 2, q\geq 2),
\end{array}
$$
$$
\begin{array}{ll}
   0\longrightarrow
   V_{\eps_1+\eps_2-\eps_{m-1+q}+\delta_1-\delta_2-\delta_3}\oplus
   V_{\eps_1+\eps_2-\eps_{m-1+q}-\delta_1+\delta_2-\delta_3}\longrightarrow
   Y \longrightarrow V_{\eps_1+\eps_2-\eps_{m-1+q}-\delta_3}\oplus &\\
   V_{\eps_1+\delta_1-\delta_2-\delta_3}\oplus
   V_{\eps_1-\delta_1+\delta_2-\delta_3} \longrightarrow 0 & (p = 1,
   q\geq 2),
\end{array}
$$
$$
\begin{array}{ll}
   0\longrightarrow
   V_{\eps_1+\eps_2-\eps_3+\delta_1-\delta_{p+n-2}-\delta_{p+n-1}}\oplus

V_{\eps_1-\eps_2+\eps_3+\delta_1-\delta_{p+n-2}-\delta_{p+n-1}}\longrightarrow
   Y &\\
   \longrightarrow
   V_{\eps_1+\eps_2-\eps_3-\delta_{p+n-1}}
   \oplus V_{\eps_1+\delta_1-\delta_{p+n-2}-\delta_{p+n-1}}\oplus
   V_{\eps_1-\eps_2+\eps_3-\delta_{p+n-1}} \longrightarrow 0&\\
   \hbox{ } &(p \geq 2, q = 1),
\end{array}
$$
$$
\begin{array}{ll}
0\longrightarrow
V_{\eps_1+\eps_2-\eps_3+\delta_1-\delta_2-\delta_3}\oplus
V_{\eps_1-\eps_2+\eps_3+\delta_1-\delta_2-\delta_3}\oplus
V_{\eps_1+\eps_2-\eps_3-\delta_1+\delta_2-\delta_3}\oplus
V_{\eps_1-\eps_2+\eps_3-\delta_1+\delta_2-\delta_3}&\\
\longrightarrow Y\longrightarrow
V_{\eps_1+\eps_2-\eps_3-\delta_3}\oplus
V_{\eps_1-\eps_2+\eps_3-\delta_3}\oplus
V_{\eps_1+\delta_1-\delta_2-\delta_3}\oplus&\\
V_{\eps_1-\delta_1+\delta_2-\delta_3} \longrightarrow 0 & (p = 1,
q = 1)
\end{array}
$$
if $m = n = p + q - 1$, then $H^{1,2}_{\fg_0} = X
\oplus V_{\eps_1+\eps_2-\eps_{m-p+q}+\delta_1-
\delta_{p+n-q-1}-\delta_{p+n-q}}\oplus V_{\eps_1-
\delta_{p+n-q}}$, where $X$ is given by the nonsplit exact
sequence of $\fg_0$-modules
$$
0\longrightarrow Y\longrightarrow X\longrightarrow
V_{\eps_1-\delta_{p+n-q}}\longrightarrow 0,
$$
and $Y$ is given by
the nonsplit exact sequence of $\fg_0$-modules
$$
\begin{array}{ll}
0\longrightarrow
V_{2\eps_1-\eps_{m-p+q}+\delta_1-2\delta_{p+n-q}} \longrightarrow
Y &\\
\longrightarrow V_{2\eps_1-\eps_{m-p+q}-\delta_{p+n-q}}
\oplus V_{\eps_1+\delta_1-2\delta_{p+n-q}} \longrightarrow 0&\\
\hbox{ }& (p \geq 3, q\geq 3),
\end{array}
$$
$$
\begin{array}{ll}
0\longrightarrow
V_{2\eps_1-\eps_{m-2+q}+\delta_1-2\delta_3}\oplus
V_{2\eps_1-\eps_{m-2+q}-\delta_1-\delta_2+\delta_3}
\longrightarrow Y &\\
\longrightarrow
V_{2\eps_1-\eps_{m-2+q}-\delta_3} \oplus
V_{\eps_1+\delta_1-2\delta_3}\oplus &\\
V_{\eps_1-\delta_1-\delta_2+\delta_3} \longrightarrow 0 &(p = 2,
q\geq 3),
\end{array}
$$
$$
\begin{array}{ll}
0\longrightarrow
V_{2\eps_1-\eps_3+\delta_1-2\delta_{p+n-2}}\oplus
V_{-\eps_1+\eps_2+\eps_3+\delta_1-2\delta_{p+n-2}} \longrightarrow
Y&\\
\longrightarrow V_{\eps_1+\delta_1-2\delta_{p+n-2}}\oplus
V_{2\eps_1-\eps_3-\delta_{p+n-2}}\oplus&\\
V_{-\eps_1+\eps_2+\eps_3-\delta_{p+n-2}} \longrightarrow 0 &(p
\geq 3, q = 2),
\end{array}
$$
$$
\begin{array}{ll}
0\longrightarrow V_{2\eps_1-\eps_3+\delta_1-2\delta_3}\oplus
V_{-\eps_1+\eps_2+\eps_3+\delta_1-2\delta_3}\oplus
V_{2\eps_1-\eps_3-\delta_1-\delta_2+\delta_3}\oplus&\\
V_{-\eps_1+\eps_2+\eps_3-\delta_1-\delta_2+\delta_3}
\longrightarrow Y\longrightarrow V_{2\eps_1-\eps_3-\delta_3}\oplus
V_{-\eps_1+\eps_2+\eps_3-\delta_3}\oplus&\\
V_{\eps_1+\delta_1-2\delta_3}\oplus
V_{\eps_1-\delta_1-\delta_2+\delta_3} \longrightarrow 0 &(p = 2, q
= 2).
\end{array}
$$
\begin{Remark} The irreducible $\fg_0$-modules in the
above listed nonsplit exact sequences are given regardless of
their parity, which can be easily recovered from the corresponding
highest weights.
\end{Remark}

%%%%%%%%%%%%%%%%%%%%%%%%%%%%%%%%%%%%%%%%%%%%%%%%%%%%%%%%%%%%%%%%%%%%%%%%%%%%%%%%%%%%%%%%%%%%%%%%%%%%%%%%%%%%
%\S2
%%%%%%%%%%%%%%%%%%%%%%%%%%%%%%%%%%%%%%%%%%%%%%%%%%%%%%%%%%%%%%%%%%%%%%%%%%%%%%%%%%%%%%%%%%%%%%%%%%%%%%%%%%%%%%%%%%%%%%%%%%%%%%%

\section*{Spencer cohomology of  $\fpsq(n)$}

\ssec{2.1. Definition of $\fpsq(n)$}

Set $\fq(n) = \lbrace X \in
\fgl(n|n)\hbox{ } |\hbox{ } [X, J_{2n}] = 0\rbrace$
  for an odd nondegenerate form $J_{2n}$, such that
$J_{2n}^2 = -1_{2n}$.  The usual choice for $J_{2n}$
is $J_{2n} =  \antidiag (1_n,\hbox{ } -1_n)$. Then we have:

$\fq(n) = \lbrace X \in \fgl(n|n)\hbox{ } | \hbox{ }X
= \diag(A,\hbox{ } A) + \antidiag(B,\hbox{ } B),
\hbox{ where } A,\hbox{ } B \in \fgl(n)\rbrace.$ Let $otr X =
tr B$. Set $\fsq(n) = \lbrace X \in \fq(n)\hbox{ }
|\hbox{ } otrX = 0 \rbrace$,

  $\fpsq(n) = \fsq(n)/\langle 1_{2n} \rangle$,

  $\fs(\fq(p)\oplus \fq(n-p)) = \lbrace X \in
\fq(p)\oplus \fq(n-p)\hbox{ }|\hbox{ } otr\fq(p) + otr\fq(n-p) = 0
\rbrace$,
$\fp\fs(\fq(p)\oplus \fq(n-p)) = \fs(\fq(p)\oplus
\fq(n-p))/\langle  1_p + 1_{n-p} \rangle$

\ssec{2.2. $\Zee$-gradings of depth 1 of $\fpsq(n)$.}

Let $V(n\hbox{ }|\hbox{ }n)$ be the standard $\fq(n)$-module.
  All $\Zee$-gradings of depth 1 of $\fg = \fpsq(n)$ are of the form
$\fg_{-1} \oplus \fg_0 \oplus \fg_{1}$, where $\fg_0 = \fp\fs (\fq(p)
\oplus \fq(n-p))$,  $p > 0$, and as $\fg_0$-modules $\fg_1 \cong
\fg_{-1}^*$, where $\fg_{-1}$ is either one of the two irreducible
$\fg_0$-modules in $V(p\hbox{ }|\hbox{ }p)^* \otimes V(n-p\hbox{
}|\hbox{ }n-p)$. Explicitly:

$\fg_{-1} = \langle  (x \pm
\Pi(x)) \otimes (y \pm \Pi(y)) \rangle, \hbox { where } x\in V(p\hbox{
}|\hbox{ }p)^*, y\in V(n-p\hbox{ }|\hbox{ }n-p)$.

Let $\eps_1, \ldots , \eps_p$ and $\delta_1, \ldots , \delta_{n-p}$ be
the standard bases of the dual spaces to the spaces of diagonal
matrices in $\fq(p)$ and $\fq(n-p)$, respectively.

\ssbegin{2.3}{Theorem}
\emph{1)}
$\fg_*(\fg_{-1},\hbox{ } \fg_0) = \fg$,

\emph{2)} all $SFs$
are of order 1 and split into the direct sum of two irreducible
$\fg_0$-submodules with highest weights $2\eps_1 - \eps_p +
\delta_1 - 2\delta_{n-p}$ and $\eps_1 - \delta_{n-p}$.
\end{Theorem}
%%%%%%%%%%%%%%%%%%%%%%%%%%%%%%%%%%%%%%%%%%%%%%%%%%%%%%%%%%%%%%%%%%%%%%%%%%%%%%%%%%%%%%%%%%%%%%%%%%%%%%%%%%%%%%%%%
%\S3
%%%%%%%%%%%%%%%%%%%%%%%%%%%%%%%%%%%%%%%%%%%%%%%%%%%%%%%%%%%%%%%%%%%%%%%%%%%%%%%%%%%%%%%%%%%%%%%%%%%%%%%%%%%%%%%%%%%%%%%%%%

\section*{Spencer cohomology of $\fosp(m|2n)$}

\ssec{3.1. Definition of $\fosp(m|2n)$ }. $\fosp(m|2n)$ is a Lie superalgebra,
which preserves a nondegenerate supersymmetric even bilinear form
on a superspace $V$, $\dim V = (m|2n)$.

\ssec{3.2.} Consider the $\Zee$-grading of depth 1 of $\fg =
\fosp(m|2n)$, which
is defined as follows: $\fg = \fg_{-1} \oplus \fg_0 \oplus \fg_1$,
where $\fg_0 = \fc \fosp(m-2|2n)$ is the central extension of
$\hat{\fg_0} = \fosp(m-2|2n)$, $\fg_1 \cong \fg_{-1}$ is the
standard $\fg_0$-module.

Let $m = 2r +2$ or $m = 2r +3$, $n > 0$. Let $\eps_1, \ldots ,
\eps_r$ and $\delta_1, \ldots , \delta_n$ be the standard bases of
the dual spaces to the spaces of diagonal matrices in $\fo(m-2)$
and $\fsp(n)$, respectively.

\ssbegin{3.3}{Theorem} \emph{1)}
$\fg_*(\fg_{-1}, \hbox{ }\fg_0) = \fg$, $\fg_*(\fg_{-1},\hbox{ }
\hat{\fg_0}) = \fg_{-1}\oplus \hat{\fg_0}.$

\emph{2)} If $k\not= 2$, then $H^{k,2}_{\fg_0} =
H^{k,2}_{\hat{\fg_0}} = 0.$ As a
$\hat{\fg_0}$-module $H^{2,2}_{\hat{\fg_0}}$ is isomorphic to
$S^2(E^2(\fg_{-1}))/E^4(\fg_{-1})$ and splits into the direct sum
of three irreducible components (analogues of the Weyl tensor, the
traceless Ricci tensor, and  the scalar curvature). The highest
weights of these components are listed in Table 18. As
$\hat{\fg_0}$-modules, $H^{2,2}_{\hat{\fg_0}} \cong
H^{2,2}_{\fg_0} \oplus S^2(\fg_{-1})$. The $\fg_0$-module
$H^{2,2}_{\fg_0}$ is irreducible.
\end{Theorem}

\section*{Spencer cohomology of $D(\alpha)$}

\ssec{4.1. Definition of $D(\alpha)$.}
$D(\alpha)$, where $\alpha \in
\Cee\backslash\lbrace 0, -1\rbrace$, is a one-parameter family
consisting of all simple Lie superalgebras for which $D(\alpha)_0
= \fsl(2)\oplus \fsl(2)\oplus \fsl(2)$  and its representation on
$D(\alpha)_1$ is $\fsl(2)\otimes \fsl(2)\otimes \fsl(2)$.

\ssec{4.2. $\Zee$-gradings of depth 1 of $D(\alpha)$.}
Let $\eps_1, \eps_2, \eps_3$ be
the standard basis of the dual space to the space of diagonal
matrices in $\fgl(1|2), V_{\lambda}$ be the irreducible
$\fsl(1|2)$-module with highest weight $\lambda$ and an even
highest vector.
All $\Zee$-gradings of depth 1 of $D(\alpha)$ are of the form
$\fg = \fg_{-1} \oplus \fg_0 \oplus \fg_1$, where $\fg_0\cong
\fgl(1|2)$. There are the following possible values of $\fg_1$
and $\fg_{-1}$ for the $\Zee$-gradings of depth 1:

a) $\fg_1 = V_{(1+\alpha )\eps_1},\hbox{ } \fg_{-1} =
V_{-\alpha\eps_1},$

b) $\fg_1 = V_{({1+\alpha\over
\alpha})\eps_1},\hbox{ } \fg_{-1} = V_{-{1\over\alpha} \eps_1},$

c) $\fg_1 = V_{({\alpha\over 1+\alpha})\eps_1},\hbox{
} \fg_{-1} = V_{{1\over 1+\alpha}\eps_1}.$

More explicitly, let $e^i_1, e^i_2$ be the basis of the standard
$\fsl(2)_i$-module $V_i$, where $i = 1,\hbox{ }2,\hbox{ } 3.$ Then
the $\Zee$-grading in case a) can be described as follows:
$\fg_0 = (\fg_0)_0 \oplus (\fg_0)_1$, where
$(\fg_0)_0 = \fsl(2)_1 \oplus \langle  {1 \hbox{ }0\choose
0-1}_2 \rangle \oplus \langle  {1 \hbox{ }0\choose 0-1}_3
\rangle, (\fg_0)_1 = V_1\otimes e_1^2\otimes e_2^3 \oplus V_1\otimes
e_2^2\otimes e_1^3$;

$\fg_1 = (\fg_1)_0 \oplus
(\fg_1)_1$, where $(\fg_1)_0 = \langle  {0 \hbox{ }1\choose
0\hbox{ }0}_2 \rangle \oplus \langle  {0 \hbox{ }1\choose 0\hbox{
}0}_3 \rangle,
(\fg_1)_1 = V_1\otimes e_1^2\otimes e_1^3$; $\fg_{-1}
= (\fg_{-1})_0 \oplus (\fg_{-1})_1$, where
$(\fg_{-1})_0 = \langle  {0 \hbox{ }0\choose 1\hbox{ }0}_2 \rangle
\oplus \langle  {0
\hbox{ }0\choose 1\hbox{ }0}_3 \rangle, (\fg_{-1})_1 = V_1\otimes
e_2^2\otimes e_2^3$; The $\Zee$-gradings in cases b)
and c) can be described similarly.

\ssbegin{4.3}{Theorem} For all $\Zee$-gradings of depth 1 of $\fg =
D(\alpha)$ we
have

\emph{1)} $\fg_*(\fg_{-1},\hbox{ } \fg_0) = \fg.$

\emph{2)} The nonzero SFs are of order 2, and for the cases
considered in sec.4.2, the $\fg_0$-module $H^{2,2}_{\fg_0}$ is
isomorphic to

a) $\Pi (V_{(2\alpha + 1)\eps_1 + \eps_2}),$

b) $\Pi (V_{({\alpha + 2\over \alpha})\eps_1 + \eps_2}),$

c) $\Pi (V_{({\alpha -
1\over \alpha + 1})\eps_1 + \eps_2}),$ respectively.
\end{Theorem}

\section*{Spencer cohomology of $AB_3$}

\ssec{5.1. Definition of $AB_3$} $AB_3$ is a simple Lie
superalgebra for which $(AB_3)_0 = \fsl(2) \oplus \fo(7)$ and its
representation on $(AB_3)_1$ is $\fsl(2) \otimes spin_7.$

\ssec{5.2. $\Zee$-grading of depth 1 of  $AB_3$.}
Let $\eps_1, \delta_1, \delta_2$ be the standard basis of the dual
space to the space of diagonal matrices in  $\fosp(2|4)$,
$V_{\lambda}$ be an irreducible $\fosp(2|4)$-module with highest
weight $\lambda$ and an even highest vector.

There is only one $\Zee$-grading of depth 1 in $\fg =  AB_3$,
namely, $\fg = \fg_{-1} \oplus \fg_0 \oplus \fg_1$, where
$\fg_{-1} = V_{-\eps_1+\delta_1+\delta_2}$,
$\fg_0 = \fc\fosp(2|4) $, $\fg_1 = V_{3\eps_1}$.

Note that $\fo(7) = V_1 \oplus \fo(5)
\oplus \Cee \oplus V_2$, where $V_1$, $V_2$ are standard
$\fo(5)$-modules. The space of the representation $spin_7$ after
restriction of $\fo(7)$  to $\fo(5)$ decomposes into  the direct sum
of two irreducible subspaces, which we denote by $U_1$ and $U_2$.
Let $e_1, e_2$ be the basis of the standard $\fsl(2)$-module $V$.
Then
$\fg_0 = (\fg_0)_0 \oplus (\fg_0)_1$, where

$(\fg_0)_0 = \langle  {1 \hbox{ }0\choose 0-1} \rangle \oplus
\fo(5) \oplus \Cee, (\fg_0)_1 = e_1\otimes U_1 \oplus e_2\otimes
U_2;$

$\fg_1 = (\fg_1)_0 \oplus (\fg_1)_1$, where
$(\fg_1)_0 = \langle  {0 \hbox{ }1\choose 0\hbox{ }0} \rangle
\oplus V_2, (\fg_1)_1 = e_1 \otimes U_2$;

$\fg_{-1} =
(\fg_{-1})_0 \oplus (\fg_{-1})_1$, where
$(\fg_{-1})_0 = \langle  {0 \hbox{ }0\choose 1\hbox{ }0} \rangle \oplus V_1,
(\fg_{-1})_1 = e_2 \otimes U_1$.

\ssbegin{5.3}{Theorem}
\emph{1)} $\fg_*(\fg_{-1},\hbox{ } \fg_0) = \fg.$

\emph{2)} The nonzero SFs are of order 1. The $\fg_0$- module $H^{1,2}_{\fg_0}$
is given by the nonsplit exact sequence of $\fg_0$- modules
$$
0\longrightarrow X\longrightarrow H^{1,2}_{\fg_0} \longrightarrow
V_{\eps_1 + 2\delta_1}\longrightarrow 0,
$$
where $X$
is given by the nonsplit exact sequence of $\fg_0$-modules
$$
0\longrightarrow \Pi(V_{4\eps_1 + 2\delta_1 +
\delta_2}) \longrightarrow X\longrightarrow V_{3\eps_1 +
2\delta_1}\longrightarrow 0.
$$
\end{Theorem}

\appendix

\renewcommand{\thechapter}{\empty}

\chapter[Appendix. The formula of dimensions]{Appendix. The formula
of dimensions of
irreducible $\fsl(n)$-modules}

Let $\eps_1,\ldots,\eps_n$ be the standard basis of the dual space
to the space of  diagonal matrices in $\fgl(n), V_\lambda$ be the
irreducible $\fsl(n)$-module with highest weight $\lambda = k_1\eps_1
+ k_2\eps_2 + \ldots + k_n\eps_n$, where $k_i\in \Zee$. Then
$$
\dim V_\lambda = \Pi^{n-1}_{i=1}\Pi^{n-i}_{j=1}
(1 + {k_i - k_{i+j}\over j}).
$$
\begin{proof} A weight $\lambda$ is
the highest weight of an irreducible $\fsl(n)$-module if and only
if $\lambda$ is a dominant integer form, i.e., if
$$
2(\lambda, \alpha_i)/(\alpha_i, \alpha_i) \in \Zee_{+}.
$$
It is known  [GG] that the inner products of the weights $\eps_i$
and of weight $\rho$, where $\rho = (\sum_{\beta \in
\Delta_{+}}\beta)/2$, with fundamental weights $\alpha_j$ are:
$$
\begin{array}{ll}
   (\eps_1, \alpha_1) = 1/(2n),\hbox{ } (\eps_1, \alpha_j) = 0
   &\hbox{for } 2\leq j\leq {n-1}; \\
   (\eps_i, \alpha_{i-1})
   = - 1/(2n),\hbox{ } (\eps_i, \alpha_i) = 1/(2n),\hbox{ }
   (\eps_i, \alpha_j) = 0 \hbox{ }(j\not= i-1, i) &\\
   \hbox{ }&\hbox{ for } 2\leq i\leq {n-1};\\
   (\eps_n, \alpha_{n-1}) = -
   1/(2n),\hbox{ } (\eps_n, \alpha_j) = 0 &\hbox{for } 1\leq j\leq
   {n-2};\\
   (\rho, \alpha_i) = 1/(2n)  &\hbox{for } 1\leq
   i\leq{n-1}.
\end{array}
$$
Thus,
$$
(\lambda, \alpha_i) = {k_i - k_{i+1}\over 2n} \hbox{ and } k_i\geq k_{i+1}.
$$
By Weyl's character formula [GG]
$$
\dim V_\lambda = \Pi_{\beta\in \Delta_{+}}(1 + {(\lambda, \beta)\over
(\rho, \beta)}).
$$
For $\fsl(n)$ we have $\Delta_{+} =
\lbrace\alpha_i + \alpha_{i+1} + \ldots + \alpha_j,\hbox{ where }
1\leq i\leq {n-1}, j\geq i\rbrace.$ Since
$$
(\lambda, \alpha_i + \ldots + \alpha_j) = 1/(2n) ((k_i - k_{i+1}) +
(k_{i+1} - k_{i+2}) + \ldots + (k_j - k_{j+1})) = {k_i - k_{j+1}\over 2n},
$$
we have
$$
\eqalign{
&\Pi_{\beta\in \Delta_{+}}(1 + {(\lambda, \beta)\over (\rho, \beta)}) =
\Pi ^{n-1}_{i=1}\Pi^{n-1}_{j=i}
(1 + {({k_i - k_{j+1}\over 2n})\over ({j-i+1\over 2n})}) = \cr
& = \Pi ^{n-1}_{i=1}\Pi^{n-i}_{j=1} (1 + {k_i - k_{i+j}\over j}).\cr}
$$
\end{proof}
\vfill \newpage

 \footnotesize

\section*{References}

\begin{itemize}
%%%%%%%%%%%%%%%%%%%%%%%%%%%%%%%%%%%%%%%%%%%%%%%%%%%%%%%%%%%%%%%%%%%%%%%%%%%%%%%%%%%%%%%%%%%%%%%%%%%%%%%%%%%%%%%%%%%%%%%%%%%%%%%%%%%%%%%%%%%%%%%%%%%%%%%%%%%%%%%%

\item[{[ALV]}] D. Alekseevskii, V. Lychagin, A. Vinogradov,
  {\it{Basic ideas and concepts of differential geometry}},
   Current Problems in Mathematics, vol.~28,
  VINITI,
  Moscow, 1988 (Russian) English translation in: Geometry, I, 1--264, Encyclopaedia Math. Sci., 28, Springer, Berlin, 1991.

%%%%%%%%%%%%%%%%%%%%%%%%%%%%%%%%%%%%%%%%%%%%%%%%%%%%%%%%%%%%%%%%%%%%%%%%%%%%%%%%%%%%%%%%%%%%%%%%%%%%%%%%%%%%%%%%%%%%%%%%%%%%%%%%%%%%%%%%%%%%%%%%%%%%%%%%%%%%%%%%

\item[{[AHS]}] M. F. Atiyah, N. J. Hitchin, and I. M. Singer,
Self-duality in four-dimensional Riemannian geometry, Proc. Roy.
Soc. London. Ser. A. 362 (1978), 425--461.

%%%%%%%%%%%%%%%%%%%%%%%%%%%%%%%%%%%%%%%%%%%%%%%%%%%%%%%%%%%%%%%%%%%%%%%%%%%%%%%%%%%%%%%%%%%%%%%%%%%%%%%%%%%%%%%%%%%%%%%%%%%%%%%%%%%%%%%%%%%%%%%%%%%%%%%%%%%%%%%%

\item[{[Fu]}] D. B. Fuks,
{\it{Cohomology of infinite-dimensional Lie algebras}}, Consultants
Bureau, New York and London, 1986.

%%%%%%%%%%%%%%%%%%%%%%%%%%%%%%%%%%%%%%%%%%%%%%%%%%%%%%%%%%%%%%%%%%%%%%%%%%%%%%%%%%%%%%%%%%%%%%%%%%%%%%%%%%%%%%%%%%%%%%%%%%%%%%%%%%%%%%%%%%%%%%%%%%%%%%%%%%%%%%%%

\item[{[GM]}]  S. I. Gelfand, Yu. I. Manin, {\em Methods of
homological algebra}. Second edition. Springer Monographs in
Mathematics. Springer-Verlag, Berlin, 2003. xx+372 pp.

%%%%%%%%%%%%%%%%%%%%%%%%%%%%%%%%%%%%%%%%%%%%%%%%%%%%%%%%%%%%%%%%%%%%%%%%%%%%%%%%%%%%%%%%%%%%%%%%%%%%%%%%%%%%%%%%%%%%%%%%%%%%%%%%%%%%%%%%%%%%%%%%%%%%%%%%%%%%%%%%

\item[{[Gi]}] S.G. Gindikin,
{\it{Integral geometry and twistors}},
  Lecture Notes in Math., vol. 970, Springer-Verlag,
Berlin, Heidelberg and New York, 1982, pp. 2--42.

%%%%%%%%%%%%%%%%%%%%%%%%%%%%%%%%%%%%%%%%%%%%%%%%%%%%%%%%%%%%%%%%%%%%%%%%%%%%%%%%%%%%%%%%%%%%%%%%%%%%%%%%%%%%%%%%%%%%%%%%%%%%%%%%%%%%%%%%%%%%%%%%%%%%%%%%%%%%%%%%

\item[{[G1]}] A. Goncharov, Infinitesimal structures related to
Hermitian symmetric spaces, Functional. Anal. Appl., 15 (1981),
221--223.

%%%%%%%%%%%%%%%%%%%%%%%%%%%%%%%%%%%%%%%%%%%%%%%%%%%%%%%%%%%%%%%%%%%%%%%%%%%%%%%%%%%%%%%%%%%%%%%%%%%%%%%%%%%%%%%%%%%%%%%%%%%%%%%%%%%%%%%%%%%%%%%%%%%%%%%%%%%%%%%%

\item[{[G2]}] A. Goncharov, Generalized conformal structures on
manifolds, Selecta Math. Sov., 6 (1987), 307--340.

%%%%%%%%%%%%%%%%%%%%%%%%%%%%%%%%%%%%%%%%%%%%%%%%%%%%%%%%%%%%%%%%%%%%%%%%%%%%%%%%%%%%%%%%%%%%%%%%%%%%%%%%%%%%%%%%%%%%%%%%%%%%%%%%%%%%%%%%%%%%%%%%%%%%%%%%%%%%%%%%

\item[{[GG]}] M. Goto and F. Grosshans,
{\it{Semisimple Lie algebras}},
Lecture Notes in  Pure and Appl. Math., vol.38, M. Dekker, New York, 1978.

%%%%%%%%%%%%%%%%%%%%%%%%%%%%%%%%%%%%%%%%%%%%%%%%%%%%%%%%%%%%%%%%%%%%%%%%%%%%%%%%%%%%%%%%%%%%%%%%%%%%%%%%%%%%%%%%%%%%%%%%%%%%%%%%%%%%%%%%%%%%%%%%%%%%%%%%%%%%%%%%

\item[{[Gu]}] V. Guillemin, The integrability problem for
$G$-structures, Trans. Amer. Math. Soc. 116 (1964), 544--560.

%%%%%%%%%%%%%%%%%%%%%%%%%%%%%%%%%%%%%%%%%%%%%%%%%%%%%%%%%%%%%%%%%%%%%%%%%%%%%%%%%%%%%%%%%%%%%%%%%%%%%%%%%%%%%%%%%%%%%%%%%%%%%%%%%%%%%%%%%%%%%%%%%%%%%%%%%%%%%%%%

\item[{[He]}] S. Helgason, {\it{Differential geometry, Lie groups and
symmetric spaces}}, Academic Press, New York, 1978.
%%%%%%%%%%%%%%%%%%%%%%%%%%%%%%%%%%%%%%%%%%%%%%%%%%%%%%%%%%%%%%%%%%%%%%%%%%%%%%%%

\item[{[K1]}] V. Kac, Lie superalgebras, Adv. Math., 26 (1977),
8--96.

%%%%%%%%%%%%%%%%%%%%%%%%%%%%%%%%%%%%%%%%%%%%%%%%%%%%%%%%%%%%%%%%%%%%%%%%%%%%%%%%%%%%%%%%%%%%%%%%%%%%%%%%%%%%%%%%%%%%%%%%%%%%%%%%%%%%%%%%%%%%%%%%%%%%%%%%%%%%%%%%

\item[{[K2]}] V. Kac, Classification of simple $\Zee$ -graded Lie
superalgebras and simple Jordan superalgebras, Comm. Algebra 13
(1977), 1375--1400.

%%%%%%%%%%%%%%%%%%%%%%%%%%%%%%%%%%%%%%%%%%%%%%%%%%%%%%%%%%%%%%%%%%%%%%%%%%%%%%%%%%%%%%%%%%%%%%%%%%%%%%%%%%%%%%%%%%%%%%%%%%%%%%%%%%%%%%%%%%%%%%%%%%%%%%%%%%%%%%%%

\item[{[Kob]}] S. Kobayashi, {\it{Transformation groups in
differential geometry}},
  Springer-Verlag, Berlin and New York, 1972.

%%%%%%%%%%%%%%%%%%%%%%%%%%%%%%%%%%%%%%%%%%%%%%%%%%%%%%%%%%%%%%%%%%%%%%%%%%%%%%%%%%%%%%%%%%%%%%%%%%%%%%%%%%%%%%%%%%%%%%%%%%%%%%%%%%%%%%%%%%%%%%%%%%%%%%%%%%%%%%%%

\item[{[KN]}] S. Kobayashi and T. Nagano,
 On filtered Lie algebras and geometric structures, III,
J. Math. Mech., 14 (1965), 679--706.

%%%%%%%%%%%%%%%%%%%%%%%%%%%%%%%%%%%%%%%%%%%%%%%%%%%%%%%%%%%%%%%%%%%%%%%%%%%%%%%%%%%%%%%%%%%%%%%%%%%%%%%%%%%%%%%%%%%%%%%%%%%%%%%%%%%%%%%%%%%%%%%%%%%%%%%%%%%%%%%%

\item[{[Kos]}] B. Kostant, Lie algebra cohomology and the
generalized Borel-Weil theorem, Ann. of Math., 74 (1961),
329--387.

%%%%%%%%%%%%%%%%%%%%%%%%%%%%%%%%%%%%%%%%%%%%%%%%%%%%%%%%%%%%%%%%%%%%%%%%%%%%%%%%%%%%%%%%%%%%%%%%%%%%%%%%%%%%%%%%%%%%%%%%%%%%%%%%%%%%%%%%%%%%%%%%%%%%%%%%%%%%%%%%

\item[{[L1]}] D. Leites, Introduction to Supermanifold theory,
  Russian Math. Surveys,
  33 (1980), 1--55.

%%%%%%%%%%%%%%%%%%%%%%%%%%%%%%%%%%%%%%%%%%%%%%%%%%%%%%%%%%%%%%%%%%%%%%%%%%%%%%%%%%%%%%%%%%%%%%%%%%%%%%%%%%%%%%%%%%%%%%%%%%%%%%%%%%%%%%%%%%%%%%%%%%%%%%%%%%%%%%%%

\item[{[L2]}] D. Leites, {\it{Supermanifold theory}},
Karelia Branch of the USSR Acad. of Sci., Petrozavodsk, 1983 (Russian).

%%%%%%%%%%%%%%%%%%%%%%%%%%%%%%%%%%%%%%%%%%%%%%%%%%%%%%%%%%%%%%%%%%%%%%%%%%%%%%%%%%%%%%%%%%%%%%%%%%%%%%%%%%%%%%%%%%%%%%%%%%%%%%%%%%%%%%%%%%%%%%%%%%%%%%%%%%%%%%%%

\item[{[L3]}] D. Leites (ed.), {\it{Seminar on supermanifolds}},
Reports of Dept. of Math. of Stockholm Univ., 1--34 (1986--89).

%%%%%%%%%%%%%%%%%%%%%%%%%%%%%%%%%%%%%%%%%%%%%%%%%%%%%%%%%%%%%%%%%%%%%%%%%%%%%%%%%%%%%%%%%%%%%%%%%%%%%%%%%%%%%%%%%%%%%%%%%%%%%%%%%%%%%%%%%%%%%%%%%%%%%%%%%%%%%%%%

\item[{[L4]}] D. Leites, Selected problems of supermanifold
theory, Duke Math.J., 54 (1987), 649--656.

%%%%%%%%%%%%%%%%%%%%%%%%%%%%%%%%%%%%%%%%%%%%%%%%%%%%%%%%%%%%%%%%%%%%%%%%%%%%%%%%%%%%%%%%%%%%%%%%%%%%%%%%%%%%%%%%%%%%%%%%%%%%%%%%%%%%%%%%%%%%%%%%%%%%%%%%%%%%%%%%

\item[{[LPS]}] D. Leites, E. Poletaeva, V. Serganova, On Einstein
equations on manifolds and supermanifolds, J. Nonlinear Math.
Physics, v.  9, 2002, no.  4, 394--425; math.DG/0306209

%%%%%%%%%%%%%%%%%%%%%%%%%%%%%%%%%%%%%%%%%%%%%%%%%%%%%%%%%%%%%%%%%%%%%%%%%%%%%%%%%%%%%%%%%%%%%%%%%%%%%%%%%%%%%%%%%%%%%%%%%%%%%%%%%%%%%%%%%%%%%%%%%%%%%%%%%%%%%%%%

\item[{[LSV]}] D. Leites, V.Serganova, G. Vinel, Classical
superspaces and related structures. In: Diff. Geom. Methods in
Theor. Physics, Proc., Rapallo, Italy 1990 (C. Bartocci, U.
Bruzzo, R. Cianci, eds.), Lecture Notes in Phys., vol.375,
Springer-Verlag, Berlin, Heidelberg and New York, 1991, pp.
286-297.

%%%%%%%%%%%%%%%%%%%%%%%%%%%%%%%%%%%%%%%%%%%%%%%%%%%%%%%%%%%%%%%%%%%%%%%%%%%%%%%%%%%%%%%%%%%%%%%%%%%%%%%%%%%%%%%%%%%%%%%%%%%%%%%%%%%%%%%%%%%%%%%%%%%%%%%%%%%%%%%%

\item[{[LRC]}] V. Lychagin, V.  Rubtsov, I. Chekalov, A
classification of Monge-Amp\`ere equations, Seminar on
supermanifolds (D. Leites, ed.), Reports of Dept. of Math. of
Stockholm Univ., 28 (1988), 16--58; Ann. Sci. \'Ecole Norm. Sup.
(4) 26 (1993), no. 3, 281--308.

%%%%%%%%%%%%%%%%%%%%%%%%%%%%%%%%%%%%%%%%%%%%%%%%%%%%%%%%%%%%%%%%%%%%%%%%%%%%%%%%%%%%%%%%%%%%%%%%%%%%%%%%%%%%%%%%%%%%%%%%%%%%%%%%%%%%%%%%%%%%%%%%%%%%%%%%%%%%%%%%

\item[{[M]}] Yu.I. Manin, {\it{Gauge field theory and complex
geometry}}, Second edition. Grundlehren der Mathematischen
Wissenschaften [Fundamental Principles of Mathematical Sciences],
289. Springer-Verlag, Berlin, 1997. xii+346 pp

%%%%%%%%%%%%%%%%%%%%%%%%%%%%%%%%%%%%%%%%%%%%%%%%%%%%%%%%%%%%%%%%%%%%%%%%%%%%%%%%%%%%%%%%%%%%%%%%%%%%%%%%%%%%%%%%%%%%%%%%%%%%%%%%%%%%%%%%%%%%%%%%%%%%%%%%%%%%%%%%

\item[{[OV]}] A.L. Onishchik, E.B. Vinberg,
{\it{Lie groups and algebraic groups}},
Springer-Verlag, Berlin and New York, 1990.

%%%%%%%%%%%%%%%%%%%%%%%%%%%%%%%%%%%%%%%%%%%%%%%%%%%%%%%%%%%%%%%%%%%%%%%%%%%%%%%%%%%%%%%%%%%%%%%%%%%%%%%%%%%%%%%%%%%%%%%%%%%%%%%%%%%%%%%%%%%%%%%%%%%%%%%%%%%%%%%%

\item[{[Pe]}] R. Penrose, Nonlinear gravitons and curved twistor
theory, Gen. Relativity and Gravitation, 7 (1976), 31--52.

%%%%%%%%%%%%%%%%%%%%%%%%%%%%%%%%%%%%%%%%%%%%%%%%%%%%%%%%%%%%%%%%%%%%%%%%%%%%%%%%%%%%%%%%%%%%%%%%%%%%%%%%%%%%%%%%%%%%%%%%%%%%%%%%%%%%%%%%%%%%%%%%%%%%%%%%%%%%%%%%

\item[{[Sh]}] W. Shmid, Die Randwerte holomorpher Functionen auf
hermitesch
  symmetrischen R\"{a}umen,
Invent. Math. 9 (1969), 61--80.

%%%%%%%%%%%%%%%%%%%%%%%%%%%%%%%%%%%%%%%%%%%%%%%%%%%%%%%%%%%%%%%%%%%%%%%%%%%%%%%%%%%%%%%%%%%%%%%%%%%%%%%%%%%%%%%%%%%%%%%%%%%%%%%%%%%%%%%%%%%%%%%%%%%%%%%%%%%%%%%%

\item[{[S1]}] V.Serganova, Classification of real simple Lie
superalgebras and symmetric superspaces, Funct. Anal. Appl.  17
(1983), 200--207.

%%%%%%%%%%%%%%%%%%%%%%%%%%%%%%%%%%%%%%%%%%%%%%%%%%%%%%%%%%%%%%%%%%%%%%%%%%%%%%%%%%%%%%%%%%%%%%%%%%%%%%%%%%%%%%%%%%%%%%%%%%%%%%%%%%%%%%%%%%%%%%%%%%%%%%%%%%%%%%%%

\item[{[S2]}] V. Serganova, Gradings of depth 1 of simple
finite-dimensional
  Lie superalgebras,
Proc. of the 19-th All-Union Algebraic Conference, Part II,
Lvov (1987), 256 (Russian).

%%%%%%%%%%%%%%%%%%%%%%%%%%%%%%%%%%%%%%%%%%%%%%%%%%%%%%%%%%%%%%%%%%%%%%%%%%%%%%%%%%%%%%%%%%%%%%%%%%%%%%%%%%%%%%%%%%%%%%%%%%%%%%%%%%%%%%%%%%%%%%%%%%%%%%%%%%%%%%%%

\item[{[St]}] S. Sternberg,
{\it{Lectures on differential geometry}},
2nd ed., Chelsea Pub. Co., New York, 1983.

%%%%%%%%%%%%%%%%%%%%%%%%%%%%%%%%%%%%%%%%%%%%%%%%%%%%%%%%%%%%%%%%%%%%%%%%%%%%%%%%%%%%%%%%%%%%%%%%%%%%%%%%%%%%%%%%%%%%%%%%%%%%%%%%%%%%%%%%%%%%%%%%%%%%%%%%%%%%%%%%

\end{itemize}

\chapter{Tables}

\section*{Table 1.  Irreducible $\fgl(n)$-submodules of $\fcpe(n)\otimes V^*$}

%{\scriptsize
\begin{tabular}{|c|c|c|}

\hline
\bf $\fgl(n)$-submodule&\bf  Highest weight&\bf   Highest vector\\

\hline

$E^2V_0^*\otimes V_0$&$\eps_1-\eps_{n-1}-\eps_n$
&$f_{n-1}\wedge \tilde f_n\otimes \tilde e_1$\\

&$-\eps_n$&$\sum_{i=1}^nf_n\wedge \tilde f_i\otimes \tilde e_i$\\

\hline

$E^2V_0^*\otimes V_0^*$&$-\eps_{n-1}-2\eps_n$&$f_{n-1}\wedge \tilde
f_n\otimes \tilde f_n$\\

& $-\eps_{n-2}-\eps_{n-1}-\eps_n$&$f_{n-2}\wedge \tilde
f_{n-1}\otimes \tilde f_n
+f_{n-1}\wedge \tilde f_n\otimes \tilde f_{n-2}$\\

&&$+f_n\wedge \tilde f_{n-2}\otimes \tilde f_{n-1}$\\

\hline

$V_0^*\wedge V_0\otimes V_0$&$2\eps_1-\eps_{n}$&$f_{n}\wedge \tilde
e_1\otimes \tilde e_1$\\

&$\eps_1+\eps_2-\eps_n$&$f_n\wedge \tilde e_1\otimes \tilde
e_2-f_n\wedge \tilde e_2\otimes e_1$\\

&$\eps_1$&$\sum_{i=1}^nf_i\wedge\tilde e_i\otimes \tilde e_1$\\

&$\eps_1$&$\sum_{i=1}^nf_i\wedge \tilde e_1\otimes \tilde e_i$\\

\hline

$V_0^*\wedge V_0\otimes V_0^*$&$\eps_1-2\eps_n$&$f_n\wedge \tilde
e_1\otimes \tilde f_n$\\

&$\eps_1-\eps_{n-1}-\eps_n$&$f_{n-1}\wedge \tilde e_1\otimes \tilde
f_n-f_n\wedge \tilde e_1\otimes \tilde
f_{n-1}$\\

&$-\eps_n$&$\sum_{i=1}^nf_i\wedge \tilde e_i\otimes \tilde f_n$\\

&$-\eps_n$&$\sum_{i=1}^nf_n\wedge \tilde e_i\otimes \tilde f_i$\\

\hline

$S^2V_0 \otimes V_0$&$3\eps_1$&$e_1 \tilde e_1 \otimes\tilde e_1$\\

&$2\eps_1+\eps_2$&$e_1\tilde e_2 \otimes \tilde e_1-
e_1\tilde e_1 \otimes \tilde e_2$\\

\hline

$S^2V_0\otimes V_0^*$&$2\eps_1-\eps_n$&$e_1\tilde e_1\otimes \tilde f_n$\\

&$\eps_1$&$\sum_{i=1}^ne_1\tilde e_i\otimes \tilde f_i$\\

\hline

$z\otimes V_0$&$\eps_1$&$\sum_{i=1}^ne_i\tilde f_i\otimes \tilde e_1$\\

\hline

$z\otimes V_0^*$&$-\eps_n$&$\sum_{i=1}^ne_i\tilde f_i\otimes \tilde f_n$\\

\hline

\end{tabular}
%}

\section*{Table 2. Irreducible $\fgl(n)$-submodules of $H^q(V_0,  \fg_*)$}

%{\scriptsize
\begin{tabular}{|c|c|c|}

\hline

\bf  q&\bf Highest weight&\bf Highest vector\\

\hline

$0$&$-\eps_n$&$f_n$\\

&$2\eps_1$&$e_1\tilde e_1$\\

&$\eps_1$&$e_1$\\

&$0$&$\tau - z$\\

\hline

$1$&$-2\eps_n$&$f_n\otimes \tilde f_n$\\

&$2\eps_1-\eps_n$&$(e_1\tilde e_1)\otimes \tilde f_n$\\

&$\eps_1-2\eps_n$&$(e_1\wedge \tilde f_n)\otimes \tilde f_n$\\

&$-\eps_n$&$(\tau -z)\otimes \tilde f_n$\\

\hline

$2$&$-2\eps_{n-1}-2\eps_n$&$(f_{n-1}\wedge \tilde f_n)\otimes \tilde
f_{n-1}\wedge \tilde f_n$\\

&$2\eps_1-\eps_{n-1}-\eps_n$&$e_1\tilde e_1\otimes \tilde
f_{n-1}\wedge \tilde f_n$\\

&$\eps_1-\eps_{n-1}-2\eps_n$&$(e_1\wedge \tilde f_n)\otimes \tilde
f_{n-1}\wedge \tilde f_n$\\

&$-\eps_{n-1}-\eps_n$&$(\tau -z)\otimes \tilde f_{n-1}\wedge \tilde f_n$\\

\hline

\end{tabular}
%}

\section*{Table 3. Irreducible $\fgl(n)$-submodules of $E^{p,0}_1$}

%{\scriptsize
\begin{tabular}{|c|c|c|c|}

\hline

\bf weight&$\bf E^{1,0}_1=\oplus_\lambda V_\lambda \otimes V_0$&
$\bf E^{2,0}_1=\oplus_\lambda V_\lambda \otimes S^2V_0$&
$\bf E^{3,0}_1=\oplus_\lambda V_\lambda \otimes S^3V_0$\\

$\bf \lambda$&$\bf V_\lambda\otimes V_0$&$\bf V_\lambda \otimes S^2V_0$&
$\bf V_\lambda\otimes S^3V_0$\\

\hline

$-\eps_n$&$\eps_1-\eps_n$&$2\eps_1-\eps_n$&$3\eps_1-\eps_n$\\

&$0$&$\eps_1$&$2\eps_1$\\

\hline

$2\eps_1$&$3\eps_1$&$4\eps_1$&$5\eps_1$\\

&$2\eps_1+\eps_2$&$2\eps_1+2\eps_2$&$4\eps_1+\eps_2$\\

&&$3\eps_1+\eps_2$&$3\eps_1+2\eps_2$\\

\hline

$\eps_1$&$2\eps_1$&$3\eps_1$&$4\eps_1$\\

&$\eps_1+\eps_2$&$2\eps_1+\eps_2$&$3\eps_1+\eps_2$\\

\hline

$0$&$\eps_1$&$2\eps_1$&$3\eps_1$\\

\hline

\end{tabular}
%}

\section*{Table 4. Irreducible $\fgl(n)$-submodules of $E^{p,1}_1$}

%{\scriptsize
\begin{tabular}{|c|c|c|c|}

\hline

\bf weight&$\bf E^{0,1}_1 = \oplus_\lambda V_\lambda$
&$\bf E^{1,1}_1 = \oplus_ \lambda V_\lambda \otimes V_0$
&$\bf E^{2,1}_1 = \oplus_ \lambda V_\lambda \otimes S^2V_0$\\

$\bf \lambda$&$\bf V_\lambda$&$\bf V_\lambda \otimes V_0$
&$\bf V_\lambda \otimes S^2V_0$\\

\hline

$-2\eps_n$&$-2\eps_n$&$\eps_1-2\eps_n$&$2\eps_1-2\eps_n$\\

&&$-\eps_n$&$\eps_1-\eps_n$\\

&&&$0$\\

\hline

$2\eps_1-\eps_n$&$2\eps_1-\eps_n$&$3\eps_1-\eps_n$&$4\eps_1-\eps_n$\\

&&$2\eps_1$&$3\eps_1$\\

&&$2\eps_1+\eps_2-\eps_n$&$3\eps_1+\eps_2-\eps_n$\\

&&&$2\eps_1+2\eps_2-\eps_n$\\

&&&$2\eps_1+\eps_2$\\

\hline

$\eps_1-2\eps_n$&$\eps_1-2\eps_n$&$2\eps_1-2\eps_n$&$3\eps_1-2\eps_n$\\

&&$\eps_1+\eps_2-2\eps_n$&$2\eps_1-\eps_n$\\

&&$\eps_1-\eps_n$&$\eps_1$\\

&&&$2\eps_1+\eps_2-2\eps_n$\\

&&&$\eps_1+\eps_2-\eps_n$\\

\hline

$-\eps_n$&$-\eps_n$&$\eps_1-\eps_n$&$2\eps_1-\eps_n$\\

&&$0$&$\eps_1$\\

\hline

\end{tabular}
%}

\section*{Table 5. Irreducible $\fgl(n)$-submodules of $E^{p,2}_1$}

\begin{tabular}{|c|c|c|}

\hline

\bf weight&$\bf E^{0,2}_1=\oplus_\lambda V_\lambda$
&$\bf E^{1,2}_1=\oplus_\lambda V_\lambda\otimes V_0$\\

$\bf \lambda$&$\bf V_\lambda$&$\bf V_\lambda\otimes V_0$\\

\hline

$-2\eps_{n-1}-2\eps_n$&$-2\eps_{n-1}-2\eps_n$&$\eps_1-2\eps_{n-1}-2\eps_n$\\

&&$-\eps_{n-1}-2\eps_n$\\

\hline

$2\eps_1-\eps_{n-1}-\eps_n$&$2\eps_1-\eps_{n-1}-\eps_n$&$3\eps_1-\eps_{n-1}-\eps_n$\\

&&$2\eps_1-\eps_n$\\

\hline

$\eps_1-\eps_{n-1}-2\eps_n$&$\eps_1-\eps_{n-1}-2\eps_n$&$2\eps_1-\eps_{n-1}-2\eps_n$\\

&&$\eps_1-2\eps_n$\\

&&$\eps_1-\eps_{n-1}-\eps_n$\\

\hline

$-\eps_{n-1}-\eps_n$&$-\eps_{n-1}-\eps_n$&$\eps_1-\eps_{n-1}-\eps_n$\\

&&$-\eps_n$\\

\hline

\end{tabular}

\section*{Table 6. Irreducible $\fgl(n)$-submodules of $E^{p,0}_0$}

\begin{tabular}{|c|c|c|c|}

\hline

&$\bf E^{1,0}_0=\oplus_U U\otimes V_0$&$\bf E^{2,0}_0=\oplus_UU\otimes S^2V_0$
&$\bf E^{3,0}_0=\oplus_UU\otimes S^3V_0$\\

$\bf U$&$\bf U\otimes V_0$&$\bf U\otimes S^2V_0$&$\bf U\otimes S^3V_0$\\

\hline

$V_0 (\mult 2)$&$2\eps_1$&$3\eps_1$&$4\eps_1$\\

&$\eps_1+\eps_2$&$2\eps_1+\eps_2$&$3\eps_1+\eps_2$\\

\hline

$V_0^* (\mult 2)$&$\eps_1-\eps_n$&$2\eps_1-\eps_n$&$3\eps_1-\eps_n$\\

&$0$&$\eps_1$&$2\eps_1$\\

\hline

$S^2V_0$&$3\eps_1$&$4\eps_1$&$5\eps_1$\\

&$2\eps_1+\eps_2$&$3\eps_1+\eps_2$&$4\eps_1+\eps_2$\\

&&$2\eps_1+2\eps_2$&$3\eps_1+2\eps_2$\\

\hline

$E^2V_0^*$&$\eps_1-\eps_{n-1}-\eps_n$&$2\eps_1-\eps_{n-1}-\eps_n$
&$3\eps_1-\eps_{n-1}-\eps_n$\\

&$-\eps_n$&$\eps_1-\eps_n$&$2\eps_1-\eps_n$\\

\hline

$\fsl(n)$&$2\eps_1-\eps_n$&$3\eps_1-\eps_n$&$4\eps_1-\eps_n$\\

&$\eps_1$&$2\eps_1$&$3\eps_1+\eps_2-\eps_n$\\

&$\eps_1+\eps_2-\eps_n$&$2\eps_1+\eps_2-\eps_n$&$3\eps_1$\\

&&$\eps_1+\eps_2$&$2\eps_1+\eps_2$\\

\hline

$\Cee(\mult 3)$&$\eps_1$&$2\eps_1$&$3\eps_1$\\

\hline

\end{tabular}

\section*{Table 7. Irreducible $\fgl(n)$-submodules of  $E^3V$}

\begin{tabular}{|c|c|c|}

\hline

\bf Space &\bf Highest vectors & \bf Highest weights\\

\hline

$E^3V_0$ &$e_1\wedge e_2\wedge e_3$&$0$\\

\hline

$(E^2V_0)(V_0^*)$&$(e_1\wedge e_2)f_3$&$-2\eps_3$\\

&$\sum_{i=1}^3 (e_1\wedge e_i)f_i$&$\eps_1$\\

\hline

$V_0(S^2V_0^*)$&$e_1f_3^2$&$\eps_1 - 2\eps_3$\\

&$\sum_{i=1}^3 e_if_if_3$&$-\eps_3$\\

\hline

$S^3V_0^*$&$f_3^3$&$-3\eps_3$\\

\hline

\end{tabular}

\section*{Table 8.  Irreducible $\fgl(m)\oplus \fgl(n)$-submodules of
$\fg_0\otimes \fg_{-1}^*$}

\begin{tabular}{|c|c|c|}

\hline

${\bf \fgl(m)\oplus \fgl(n)}${\bf -module}&
{\bf Highest weight}&{\bf Highest vector}\\

\hline

$(V \otimes V^*)/\Cee \otimes (U^*\otimes V)$&$2\eps_1-\eps_m-\delta_n$
&$(e_1\otimes \tilde e_m)\otimes (\tilde f_n\otimes e_1)$\\

&$\eps_1-\delta_n $&$v^1_{\lambda}=\sum_{i=1}^m(e_i \otimes \tilde e_i)
\otimes (\tilde f_n\otimes e_1)-$\\

&&$-m\sum_{i=1}^m(e_1 \otimes \tilde e_i) \otimes (\tilde f_n\otimes e_i)$\\

&$\eps_1+\eps_2-\eps_m-\delta_n$
&$(e_1\otimes \tilde e_m) \otimes (\tilde f_n\otimes e_2)-$\\

&(if $m\geq 3$)&$-(e_2 \otimes \tilde e_m)\otimes (\tilde f_n \otimes e_1)$\\

\hline

$(U\otimes U^*)/\Cee\otimes (U^*\otimes V)$&$\eps_1+\delta_1-2\delta_n$
&$(f_1\otimes \tilde f_n)\otimes (\tilde f_n\otimes e_1)$\\

&$\eps_1-\delta_n$&$v^2_{\lambda}=\sum_{i=1}^n(f_i\otimes \tilde f_i)\otimes
(\tilde f_n\otimes e_1)-$\\

&&$-n\sum_{i=1}^n(f_i\otimes \tilde f_n)\otimes (\tilde f_i\otimes e_1)$\\

&$\eps_1+\delta_1-\delta_{n-1}-\delta_n$
&$(f_1\otimes \tilde f_{n-1})\otimes (\tilde f_n\otimes e_1)-$\\

&(if $n\geq 3$)&$-(f_1\otimes \tilde f_n)\otimes (\tilde f_{n-1}\otimes e_1)$\\

\hline

$\Cee\otimes (U^*\otimes V)$&$\eps_1-\delta_n$
&$v^3_{\lambda}=(n\sum_{i=1}^me_i\otimes\tilde e_i +$\\

&&$+m\sum_{i=1}^nf_i\otimes \tilde f_i)\otimes (\tilde f_n\otimes e_1)$\\

\hline

\end{tabular}

\section*{Table 9. Irreducible $\fgl(m)\oplus \fgl(n)$-submodules of
  $(U^*\otimes V)\otimes (U^*\otimes V)$}

\begin{tabular}{|c|c|c|}

\hline

$\bf \fgl(m)\oplus \fgl(n)\hbox{-module}$&
\bf Highest weight&\bf Highest vector\\

\hline

$\Lambda^2U^*\otimes S^2V$&$2\eps_1-2\delta_n$
&$(\tilde f_n\otimes e_1)\otimes (\tilde f_n\otimes e_1)$\\

\hline

$S^2U^*\otimes \Lambda^2V$&$\eps_1+\eps_2-\delta_{n-1}-\delta_n$
&$(\tilde f_n\otimes e_1)\otimes (\tilde f_{n-1}\otimes e_2)-$\\

&&$-(\tilde f_n\otimes e_2)\otimes (\tilde f_{n-1}\otimes e_1)-$\\

&&$-(\tilde f_{n-1}\otimes e_1)\otimes (\tilde f_n\otimes e_2)+$\\

&&$+(\tilde f_{n-1}\otimes e_2)\otimes (\tilde f_n\otimes e_1)$\\

\hline

$\Lambda^2U^* \otimes \Lambda^2V$&$\eps_1+\eps_2-2\delta_n$&
$(\tilde f_n\otimes e_1)\otimes (\tilde f_n\otimes e_2)-$\\

&&$-(\tilde f_n\otimes e_2)\otimes (\tilde f_n\otimes e_1)$\\

\hline

$S^2U^*\otimes S^2 V$&$2\eps_1-\delta_{n-1}-\delta_n$
&$(\tilde f_{n-1}\otimes e_1)\otimes (\tilde f_n\otimes e_1)-$\\

&&$-(\tilde f_n\otimes e_1)\otimes (\tilde f_{n-1}\otimes e_1)$\\

\hline

\end{tabular}

\section*{Table 10. Irreducible $\fgl(n)$-submodules of  $C^{k,2}_{\fgl(n)}$}

{\scriptsize
\begin{tabular}{|c|c|c|}

\hline

\bf  k&\bf  Highest weight&\bf  Highest vector\\

\hline

$2$&$\lambda_1 = \delta_1-3\delta_n$
&$v_{\lambda_1} = (f_1\otimes \tilde f_n)\otimes \tilde f_n^2$\\

&$\lambda_2 = -2\delta_n$
&$v_{\lambda_2} = \sum_{j=1}^n(f_j \otimes \tilde f_n) \otimes \tilde
f_j\tilde f_n$\\

&$\lambda_3 = -2\delta_n$
&$v_{\lambda_3} = \sum_{j=1}^n(f_j\otimes \tilde f_j)\otimes \tilde f_n^2$\\

&$\lambda_4 = \delta_1-\delta_{n-1}-2\delta_n (\hbox{ if } n\geq 3)$
&$v_{\lambda_4} = (f_1\otimes \tilde f_{n-1})\otimes \tilde f_n^2-
(f_1\otimes \tilde f_n)\otimes \tilde f_{n-1}\tilde f_n$\\

&$\lambda_5 =  -\delta_{n-1}-\delta_n$
&$v_{\lambda_5} = \sum_{j=1}^n((f_j \otimes \tilde f_{n-1}) \otimes
\tilde f_j\tilde f_n-(f_j \otimes \tilde f_n)\otimes \tilde f_j\tilde
f_{n-1})$\\

\hline

$3\leq k \leq n + 1$&$\lambda_1 = \delta_1-\delta_{n-k+2}-\ldots
-\delta_{n-1}-3\delta_n$
&$v_{\lambda_1} = (f_1\otimes \tilde f_{n-k+2}\wedge \ldots \wedge
\tilde f_n) \otimes \tilde f_n^2$\\

&$\lambda_2 = -\delta_{n-k+2}-\ldots -\delta_{n-1}-2\delta_n$
&$v_{\lambda_2} = \sum_{j=1}^n(f_j \otimes \tilde f_{n-k+2}\wedge\ldots \wedge
\tilde f_n) \otimes \tilde f_j\tilde f_n$\\

\hline

$3\leq k \leq n$&$\lambda_3 = -\delta_{n-k+3}-\ldots -\delta_n -2\delta_n$
&$v_{\lambda_3} = \sum_{j=1}^n(f_j \otimes \tilde f_j\wedge \tilde f_{n-k+3}
\wedge \ldots \wedge \tilde f_n)\otimes \tilde f_n^2$\\

&$\lambda_4 = \delta_1-\delta_{n-k+1}-\ldots -\delta_{n-1}-2\delta_n$
&$v_{\lambda_4} = \sum_{j=0}^{k-1}
(-1)^{(k-1)j}(f_1\otimes \tilde f_{s^j(n-k+1)}\wedge\ldots $\\

&&$\wedge \tilde f_{s^j(n-1)})\otimes \tilde f_{s^j(n)}\tilde f_n$\\

&$\lambda_5 = -\delta_{n-k+1}-\delta_{n-k+2} -\ldots -\delta_n$
&$v_{\lambda_5} = \sum_{j=0}^{k-1}
(-1)^{(k-1)j}\sum_{i=1}^n(f_i \otimes \tilde f_{s^j(n-k+1)}\wedge \ldots$\\

&&$ \wedge \tilde f_{s^j(n-1)})\otimes \tilde f_i\tilde f_{s^j(n)}$\\

$3\leq k \leq n - 1$&$\lambda_6 = -\delta_{n-k+2}-\ldots
-\delta_{n-1}-2\delta_n$
&$v_{\lambda_6} = \sum_{j=0}^{k-2}
(-1)^{(k-2)j}\sum_{i=1}^n(f_i \otimes \tilde f_i\wedge \tilde
f_{t^j(n-k+2)}\wedge\ldots$\\

&&$\wedge \tilde f_{t^j(n-1)}) \otimes \tilde f_{t^j(n)}\tilde f_n$\\

\hline

\end{tabular}
}

\section*{Table 11. Irreducible
$\fgl(n)$-submodules of $C^{k+1,1}_{\fgl(n)}$}

{\scriptsize
\begin{tabular}{|c|c|c|}

\hline

\bf  k&\bf Highest weight&\bf Highest vector\\

\hline

$2\leq k \leq n$&$\beta_1 = \delta_1-\delta_{n-k+1}-\ldots
-\delta_{n-1}-2\delta_n$
&$v_{\beta_1} = (f_1\otimes \tilde f_{n-k+1} \wedge \ldots \wedge
\tilde f_{n-1}\wedge\tilde f_n)\otimes \tilde f_n$\\

&$\beta_2 = -\delta_{n-k+1}-\ldots -\delta_{n-1}-\delta_n$
&$v_{\beta_2} = \sum_{j=1}^nf_j \otimes \tilde f_{n-k+1} \wedge \ldots \wedge
\tilde f_{n-1}\wedge \tilde f_n)\otimes \tilde f_j$\\

$2\leq k\leq n - 1$&$\beta_3 = -\delta_{n-k+2}-\ldots -\delta_{n-1}-2\delta_n$
&$v_{\beta_3} = \sum_{j=6}^n(f_j \otimes \tilde f_j \wedge \tilde
f_{n-k+2} \wedge \ldots
\wedge \tilde f_{n-1}\wedge \tilde f_n)\otimes \tilde f_n$\\

&$\beta_4 = -\delta_{n-k+1}-\ldots -\delta_{n-1}-\delta_n$
&$v_{\beta_4} = \sum_{j=0}^{k-1}
(-1)^{(k-3)j}\sum_{i=1}^n(f_i \otimes \tilde f_i \wedge
\tilde f_{s^j(n-k+1)}\wedge \ldots$\\

&&$ \wedge f_{s^j(n-1)})\otimes \tilde f_{s^j(n)}$\\

$2\leq k \leq n - 2$&$\beta_5 = \delta_1-\delta_{n-k}-\ldots
-\delta_{n-1}-\delta_n$
&$v_{\beta_5} = \sum_{j=0}^k (-1)^{kj}(f_1 \otimes \tilde
f_{r^j(n-k)}\wedge \ldots$\\

&&$ \wedge \tilde f_{r^j(n-1)}) \otimes \tilde f_{r^j(n)}$\\

\hline

\end{tabular}
}

\section*{Table 12. Irreducible
$\fgl(n)$-submodules of  $\fg_k$}

\begin{tabular}{|c|c|c|}

\hline

\bf  k&\bf Highest weight&\bf Highest vector\\

\hline

$2\leq k \leq n - 1$&$\gamma_1 =
\delta_1-\delta_{n-k}-\delta_{n-k+1}-\ldots -\delta_n$
&$v_{\gamma_1} = f_1\otimes\tilde f_{n-k}\wedge \ldots \wedge \tilde f_n$\\

$2\leq k \leq n - 2$&$\gamma_2 =
-\delta_{n-k+1}-\delta_{n-k+2}-\ldots -\delta_n$
&$v_{\gamma_2} = \sum_{j=1}^nf_j\otimes \tilde f_j\wedge \tilde f_{n-k+1}
\wedge \ldots \wedge \tilde f_n$\\

\hline

\end{tabular}

\section*{Table 13. Structure functions of  $\fsl(m|n)$ endowed with
the standard $\Zee$-grading}

\begin{tabular}{|c|c|c|c|}

\hline

\bf m&\bf n&\bf $H^{1,2}_{\fg_0}$&\bf $H^{2,2}_{\fg_0}$\\

\hline

$3$&$2$&$2\eps_1-\eps_3+\delta_1-2\delta_2$&------\\

\hline

$2$&$3$&$2\eps_1-\eps_2+\delta_1-2\delta_3$&------\\

\hline

$\geq 4$&$2$&$2\eps_1-\eps_m+\delta_1-2\delta_2$
&$\eps_1+\eps_2+\eps_3-\eps_m-\delta_1-\delta_2$\\

\hline

$2$&$\geq 4$&$2\eps_1-\eps_2+\delta_1-2\delta_n$
&$\eps_1+\eps_2+\delta_1-\delta_{n-2}-\delta_{n-1}-\delta_n$\\

\hline

&&$2\eps_1-\eps_m+\delta_1-2\delta_n$&\\

$\geq 3$&$\geq
3$&$\eps_1+\eps_2-\eps_m+\delta_1-\delta_{n-1}-\delta_n$&------\\

&&$\eps_1-\delta_n \hbox{ }(\hbox{if } m=n) $&\\

\hline

$2$&$2$&$2\eps_1-\eps_2+\delta_1-2\delta_2$&\\
&&&------\\
&&$\eps_1-\delta_2$&\\

\hline

\end{tabular}

\section*{Table 14. Irreducible $\fgl(m)\oplus \fgl(n)$submodules of
$\fg_0\otimes \Lambda^2\fg_{-1}^*$}

%\begin{center}
{\scriptsize
\begin{tabular}{|c|c|c|}
\hline
$\bf \fgl(m)\oplus \fgl(n)\hbox{-module}$&\bf Highest weight&\bf
Highest vector\\

\hline
$V\otimes V^*/\Cee\otimes \Lambda^2U^*\otimes S^2V$&$3\eps_1-\eps_m-2\delta_n$
&$(e_1\otimes \tilde e_m)\otimes (\tilde f_n \otimes e_1) \wedge
(\tilde f_n\otimes e_1)$\\

&$2\eps_1-2\delta_n$&$v^1_{\lambda}=\sum_{i=1}^m
(e_i\otimes \tilde e_i)\otimes (\tilde f_n\otimes e_1)\wedge (\tilde
f_n\otimes e_1)-$\\

&&$m\sum_{i=1}^m(e_1\otimes \tilde e_i) \otimes (\tilde f_n\otimes
e_i)\wedge (\tilde f_n\otimes e_1)$\\

&$\eps_1+\eps_2-2\delta_n$& $v^1_{\lambda}=\sum_{i=1}^m((e_1\otimes \tilde e_i)
\otimes (\tilde f_n\otimes e_i)\wedge (\tilde f_n\otimes e_2)-$\\

&&$(e_2\otimes \tilde e_i)\otimes (\tilde f_n\otimes e_i)\wedge
(\tilde f_n\otimes e_1))$\\

&$2\eps_1+\eps_2-\eps_m-2\delta_n \hbox{ }(m\geq 3)$
&$(e_1\otimes \tilde e_m)\otimes (\tilde f_n\otimes e_1)\wedge
(\tilde f_n\otimes e_2)-$\\

&&$(e_2\otimes \tilde e_m) \otimes (\tilde f_n\otimes e_1) \wedge
(\tilde f_n \otimes e_1)$\\

\hline

$V\otimes V^*/\Cee\otimes S^2U^*\otimes \Lambda^2V$
&$2\eps_1+\eps_2-\eps_m-\delta_{n-1}-\delta_n$
&$(e_1\otimes \tilde e_m) \otimes (\tilde f_{n-1}\otimes e_1) \wedge
(\tilde f_n \otimes e_2)-$\\

&&$(e_1\otimes \tilde e_m) \otimes (\tilde f_n\otimes e_1) \wedge
(\tilde f_{n-1} \otimes e_2)$\\

&$2\eps_1-\delta_{n-1}-\delta_n \hbox{ }(m\geq 3)$ &$v^1_{\lambda}=\sum_{i=1}^m
(e_1\otimes \tilde e_i) \otimes (\tilde f_{n-1}\otimes e_i) \wedge
(\tilde f_n \otimes e_1)-$\\

&&$(e_1\otimes \tilde e_i) \otimes (\tilde f_n\otimes e_i) \wedge (\tilde f_{n-1}
\otimes e_1)$\\

&$\eps_1+\eps_2-\delta_{n-1}-\delta_n$
&$v^1_{\lambda}=\sum_{i=1}^m((e_1 \otimes \tilde e_i)\otimes (\tilde
f_{n-1} \otimes e_i)
\wedge (\tilde f_n \otimes e_2)-$\\

&$(m\geq 3)$&$(e_1 \otimes \tilde e_i)\otimes (\tilde f_{n} \otimes e_i)
\wedge (\tilde f_{n-1} \otimes e_2)-$\\

&&$(e_2 \otimes \tilde e_i)\otimes (\tilde f_{n-1} \otimes e_i)
\wedge (\tilde f_{n} \otimes e_1)+$\\

&&$(e_2 \otimes \tilde e_i)\otimes (\tilde f_{n} \otimes e_i)
\wedge (\tilde f_{n-1} \otimes e_1))$\\

&$\eps_1+\eps_2+\eps_3-\eps_m-\delta_{n-1}-\delta_n$
&$\sum_{j=0}^2((e_{s^j(1)} \otimes \tilde e_m)\otimes (\tilde f_{n-1} \otimes
e_{s^j(2)})\wedge (\tilde f_n \otimes e_{s^j(3)})-$\\

&$(m \geq 4)$
&$(e_{s^j(1)} \otimes \tilde e_m)\otimes (\tilde f_{n} \otimes
e_{s^j(2)})\wedge (\tilde f_{n-1} \otimes e_{s^j(3)}))$\\

\hline

$U\otimes U^*/\Cee\otimes \Lambda^2U^*\otimes
S^2V$&$2\eps_1+\delta_1-3\delta_n$
&$(f_1\otimes \tilde f_n)\otimes (\tilde f_n\otimes e_1)\wedge
(\tilde f_n\otimes e_1)$\\

&$2\eps_1-2\delta_n$
&$v^2_{\lambda}=\sum_{i=1}^n(f_i\otimes \tilde f_i)\otimes (\tilde
f_n\otimes e_1)
\wedge (\tilde f_n\otimes e_1)-$\\

&&$n\sum_{i=1}^n(f_i\otimes \tilde f_n)\otimes (\tilde f_i\otimes e_1)\wedge
(\tilde f_n\otimes e_1)$\\

&$2\eps_1-\delta_{n-1}-\delta_n$
&$v^2_{\lambda}=\sum_{i=1}^n((f_i\otimes \tilde f_{n-1})\otimes
(\tilde f_i\otimes e_1)
\wedge (\tilde f_n\otimes e_1)-$\\

&&$(f_i\otimes \tilde f_n)\otimes (\tilde f_i\otimes e_1)\wedge
(\tilde f_{n-1}\otimes e_1))$\\

&$2\eps_1+\delta_1-\delta_{n-1}-2\delta_n$
&$(f_1\otimes \tilde f_n)\otimes (\tilde f_{n-1}\otimes e_1)\wedge
(\tilde f_n\otimes e_1)-$\\

&$(n \geq 3)$&$(f_1\otimes \tilde f_{n-1})\otimes (\tilde f_n\otimes
e_1)\wedge (\tilde f_n\otimes e_1)$\\

\hline

$U\otimes U^*/\Cee\otimes S^2U^*\otimes \Lambda^2V$
&$\eps_1+\eps_2+\delta_1-\delta_{n-1}-2\delta_n$
&$(f_1\otimes \tilde f_n)\otimes (\tilde f_{n-1}\otimes e_1)\wedge
(\tilde f_n\otimes e_2)-$\\

&&$(f_1\otimes \tilde f_n)\otimes (\tilde f_{n}\otimes e_1)\wedge (\tilde
  f_{n-1}\otimes e_2)$\\

&$\eps_1+\eps_2-2\delta_n \hbox{ }(n \geq 3)$
&$v^2_{\lambda}=\sum_{i=1}^n((f_i \otimes \tilde f_n) \otimes (\tilde
f_i\otimes e_1)\wedge (\tilde f_n\otimes e_2)-$\\

&&$(f_i \otimes \tilde f_n) \otimes (\tilde f_n\otimes e_1)\wedge
(\tilde f_i\otimes e_2))$\\

&$\eps_1+\eps_2-\delta_{n-1}-\delta_n \hbox{ }(n \geq 3)$
&$v^2_{\lambda}=\sum_{i=1}^n((f_i \otimes \tilde f_n) \otimes (\tilde
f_i\otimes e_1)\wedge (\tilde f_{n-1}\otimes e_2)-$\\

&&$(f_i \otimes \tilde f_n) \otimes (\tilde f_i\otimes e_2)\wedge
(\tilde f_{n-1}\otimes e_1)-$\\

&&$(f_i \otimes \tilde f_{n-1}) \otimes (\tilde f_i\otimes e_1)\wedge
(\tilde f_{n}\otimes e_2)+$\\

&&$(f_i \otimes \tilde f_{n-1}) \otimes (\tilde f_i\otimes e_2)\wedge
(\tilde f_{n}\otimes e_1))$\\

&$\eps_1+\eps_2+\delta_1-\delta_{n-2}-\delta_{n-1}-\delta_n$
&$\sum_{j=0}^2((f_1\otimes \tilde f_{t^j(n-2)})\otimes (\tilde
f_{t^j(n-1)}\otimes e_1)\wedge (\tilde f_{t^j(n)}\otimes e_2)-$\\

&$(n \geq 4)$&$(f_1\otimes \tilde f_{t^j(n-2)})\otimes (\tilde
f_{t^j(n-1)}\otimes e_2)\wedge (\tilde f_{t^j(n)}\otimes e_1))$\\

\hline

$\Cee\otimes \Lambda^2U^*\otimes S^2V$&$2\eps_1-2\delta_n$
&$v^3_{\lambda}=(n\sum_{i=1}^me_i\otimes \tilde e_i+$\\

&&$+m\sum_{i=1}^nf_i\otimes \tilde f_i)
\otimes (\tilde f_n\otimes e_1)\wedge (\tilde f_n\otimes e_1)$\\

\hline

$\Cee\otimes S^2U^*\otimes \Lambda^2V$&$\eps_1+\eps_2-\delta_{n-1}-
\delta_n$
&$v^3_{\lambda}=(n\sum_{i=1}^me_i\otimes \tilde e_i+$\\

&&$+m\sum_{i=1}^nf_i\otimes \tilde f_i)
\otimes ((\tilde f_{n-1}\otimes e_1)\wedge (\tilde f_n\otimes e_2)-$\\

&&$(\tilde f_{n}\otimes e_1)\wedge (\tilde f_{n-1}\otimes e_2))$\\

\hline

\end{tabular}
}
%\end{center}

\section*{Table 15. Irreducible $\fgl(m)\oplus \fgl(n)$-submodules of
$\fg_1\otimes \Lambda^2\fg_{-1}^*$}

{\scriptsize
\begin{tabular}{|c|c|c|}

\hline

$\bf \fgl(m)\oplus \fgl(n)\hbox{-module}$&\bf Highest weight&\bf
Highest vector\\

\hline

$(\Lambda^2U^*\otimes U^*)\otimes (S^2V\otimes V)$
&$3\eps_1-3\delta_n$
&$(\tilde f_n\otimes e_1)\otimes (\tilde f_n\otimes e_1)\wedge
(\tilde f_n\otimes e_1)$\\

&$2\eps_1+\eps_2-3\delta_n$
&$(\tilde f_n\otimes e_1)\otimes (\tilde f_n\otimes e_2)\wedge
(\tilde f_n\otimes e_1)-$\\

&&$(\tilde f_n\otimes e_2)\otimes (\tilde f_n\otimes e_1)\wedge
(\tilde f_n\otimes e_1)$\\

&$3\eps_1-\delta_{n-1}-2\delta_n$
&$(\tilde f_{n-1}\otimes e_1)\otimes (\tilde f_n\otimes e_1)\wedge
(\tilde f_n\otimes e_1)-$\\

&&$(\tilde f_n\otimes e_1)\otimes (\tilde f_{n-1}\otimes e_1)\wedge
(\tilde f_n\otimes e_1)$\\

&$2\eps_1+\eps_2-\delta_{n-1}-2\delta_n$
&$v^1_{\lambda}=(\tilde f_n\otimes e_1)\otimes (\tilde f_{n-1}\otimes
e_1)\wedge (\tilde f_n\otimes e_2)-$\\

&&$(\tilde f_n\otimes e_2)\otimes (\tilde f_{n-1}\otimes e_1)\wedge
(\tilde f_n\otimes e_1)-$\\

&&$(\tilde f_{n-1}\otimes e_1)\otimes (\tilde f_n\otimes e_1)\wedge
(\tilde f_n\otimes e_2)+$\\

&&$(\tilde f_{n-1}\otimes e_2)\otimes (\tilde f_n\otimes e_1)\wedge
(\tilde f_n\otimes e_1)$\\

\hline

$(S^2U^*\otimes U^*)\otimes (\Lambda^2V\otimes  V)$
&$2\eps_2+\eps_2-\delta_{n-1}-2\delta_n$
&$v^2_{\lambda}=(\tilde f_n\otimes e_1)\otimes (\tilde f_{n-1}\otimes
e_1)\wedge (\tilde f_n\otimes e_2)-$\\

&&$(\tilde f_n\otimes e_1)\otimes (\tilde f_{n-1}\otimes e_2)\wedge
(\tilde f_n\otimes e_1)$\\

&$\eps_1+\eps_2+\eps_3-\delta_{n-1}-2\delta_n$
&$\sum_{j=0}^2((\tilde f_n\otimes e_{s^j(1)})\otimes (\tilde
f_{n-1}\otimes e_{s^j(2)})\wedge (\tilde f_n\otimes e_{s^j(3)})-$\\

&$(m \geq 3)$
&$(\tilde f_n\otimes e_{s^j(1)})\otimes (\tilde f_{n}\otimes
e_{s^j(2)})\wedge (\tilde f_{n-1}\otimes e_{s^j(3)})$\\

&$2\eps_1+\eps_2-\delta_{n-2}-\delta_{n-1}-\delta_n$
&$\sum_{j=0}^2((\tilde f_{t^j(n-2)}\otimes e_1)\otimes (\tilde
f_{t^j(n-1)}\otimes e_2)\wedge (\tilde f_{t^j(n)}\otimes e_1)-$\\

&$(n \geq 3)$
&$(\tilde f_{t^j(n-2)}\otimes e_1)\otimes (\tilde f_{t^j(n-1)}\otimes
e_1)\wedge (\tilde f_{t^j(n)}\otimes e_2))$\\

&$\eps_1+\eps_2+\eps_3-\delta_{n-2}-\delta_{n-1}-\delta_n$
&$\sum_{i=0}^2\sum_{j=0}^2
((\tilde f_{t^j(n-2)}\otimes e_{s^i(1)})\otimes (\tilde f_{t^j(n-1)}\otimes
  e_{s^i(2)})\wedge$\\

&$(m,\hbox{ } n \geq 3)$&$(\tilde f_{t^j(n)}\otimes
e_{s^i(3)})-(\tilde f_{t^j(n-2)}\otimes e_{s^i(1)})\otimes $\\

&&$(\tilde f_{t^j(n)}\otimes
  e_{s^i(2)})\wedge (\tilde f_{t^j(n-1)}\otimes e_{s^i(3)}))$\\

\hline

\end{tabular}
}

\section*{Table 16. Spencer cohomology of  $\bf \fsl(m|n)$
endowed with a  $\Zee$-grading, where  $ \fg_0 = \fc(\fsl(m|q)\oplus
\fsl(n-q))$}

\begin{tabular}{|c|c|c|c|c|}

\hline

$\bf m$&$\bf q$&$\bf n - q$&$\bf H^{1,2}_{\fg_0}$&$\bf H^{2,2}_{\fg_0}$\\

\hline

$\geq 2$&$\geq 1$&$\geq 3$&$m \not= q \pm 1$&\\

&&&$2\eps_1-\eps_{m+q}+\delta_1-2\delta_{n-q}$&\\

&&&&------\\

&&&$\eps_1+\eps_2-\eps_{m+q}+\delta_1-\delta_{n-q-1}-\delta_{n-q}$&\\

&&&$\eps_1-\delta_{n-q} \hbox{ (if } m = n)$&\\

\hline

$\geq 3$&$\geq 1$&$2$&$m \not= q - 1$&\\

&&&$2\eps_1-\eps_{m+q}+\delta_1-2\delta_2$&$\eps_1+
\eps_2+\eps_3-\eps_{m+q}-\delta_1-\delta_2$\\

&&&$\eps_1-\delta_2 \hbox{ (if } m = n)$&\\

\hline

$2$&$\geq 2$&$2$&$q \not= 3$&$\eps_1+\eps_2+\eps_3
-\eps_{q+2}-\delta_1-\delta_2$\\

&&&$2\eps_1-\eps_{q+2}+\delta_1-2\delta_2$&\\

\hline

$2$&$1$&$2$&$2\eps_1-\eps_3+\delta_1-2\delta_2$&$\eps_1+\eps_3-\delta_1-\delta_2$\\

\hline

$1$&$\geq 1$&$\geq 3$&$q \not= 2$&\\

&&&$2\eps_1-\eps_{q+1}+\delta_1-2\delta_{n-q}$&------\\

&&&$\eps_1+\eps_2-\eps_{q+1}+\delta_1-\delta_{n-q-1}-\delta_{n-q}$&\\

\hline

$1$&$\geq 1$&$2$&$q\not=
2$&$\eps_1+2\eps_2-\eps_{q+1}-\delta_1-\delta_2 \hbox{ }(q \not= 1)$\\

&&&$2\eps_1-\eps_{q+1}+\delta_1-2\delta_2$&$2\eps_2-\delta_1-\delta_2
\hbox{ }(q = 1)$\\

\hline

$0$&$2$&$2$&&$3\eps_1-\eps_2-\delta_1-\delta_2$\\
&&&------&\\
&&&&$\eps_1+\eps_2+\delta_1-3\delta_2$\\

\hline

$0$&$2$&$\geq
3$&$2\eps_1-\eps_2+\delta_1-\delta_{n-3}-\delta_{n-2}$&$\eps_1+\eps_2+\delta_1-3\delta_{n-2}$\\

\hline

$0$&$\geq 3$&$2$&$\eps_1+\eps_2-\eps_q+\delta_1-2\delta_2$&
$3\eps_1-\eps_q-\delta_1-\delta_2$\\

\hline

$0$&$\geq 3$&$\geq 3$&$2\eps_1-\eps_q+\delta_1-\delta_{n-q-1}-
\delta_{n-q}$&\\

&&&&------\\

&&&$\eps_1+\eps_2-\eps_q+\delta_1-2\delta_{n-q}$&\\

\hline

\end{tabular}

\section*{Table 17.  Spencer cohomology of $\fsl(m|n)$ endowed with a
$\Zee$-grading, where  $\fg_0 = \fc(\fsl(m-p|q)\oplus \fsl(p|n-q))$}

{\scriptsize
\begin{tabular}{|c|c|c|c|c|c|}

\hline

$\bf m-p$&$\bf q$&$\bf p$&$\bf n-q$&$\bf H^{1,2}_{\fg_0}$&$\bf
H^{2,2}_{\fg_0}$\\

\hline

$0$&$2$&$2$&$0$&$2\eps_1-\eps_2+\delta_1-2\delta_2$&\\

&&&&&------\\

&&&&$\eps_1-\delta_2$&\\

\hline

$0$&$2$&$3$&$0$&$2\eps_1-\eps_2+\delta_1-2\delta_3$&------\\

\hline

$0$&$3$&$2$&$0$&$2\eps_1-\eps_3+\delta_1-2\delta_2$&------\\

\hline

$0$&$2$&$\geq 4$&$0$&$2\eps_1-\eps_2+\delta_1-2\delta_p$
&$\eps_1+\eps_2+\delta_1-\delta_{p-2}-\delta_{p-1}-\delta_p$\\

\hline

$0$&$\geq 4$&$2$&$0$&$2\eps_1-\eps_q+\delta_1-2\delta_2$
&$\eps_1+\eps_2+\eps_3-\eps_q-\delta_1-\delta_2$\\

\hline

$0$&$\geq 3$&$\geq 3$&$0$&$2\eps_1-\eps_q+\delta_1-2\delta_p$&\\

&&&&$\eps_1+\eps_2-\eps_q+\delta_1-\delta_{p-1}-\delta_p$&------\\

&&&&$\eps_1-\delta_p \hbox{ }(m=n)$&\\

\hline

$0$&$2$&$\geq 1$&$\geq 1$&$n \not= p+q+1$&\\

&&&&$2\eps_1-\eps_2+\delta_1-\delta_{p+n-3}-\delta_{p+n-2}$&$\eps_1+\eps_2+\delta_1-3\delta_{p+n-2}$\\

&&&&$\eps_1-\delta_{p+n-2} \hbox{ }(m=n)$&\\

\hline

$\geq 1$&$\geq 1$&$2$&$0$&$m \not= p+q+1$&\\

&&&&$\eps_1+\eps_2-\eps_{m-p+q}+\delta_1-2\delta_2$&$3\eps_1-\eps_{m-p+q}-\delta_1-\delta_2$\\

&&&&$\eps_1-\delta_2 \hbox{ }(m=n)$&\\

\hline

$0$&$\geq 3$&$\geq 1$&$\geq 1$&$n \not= p+q \pm 1$&\\

&&&&$\eps_1+\eps_2-\eps_q+\delta_1-2\delta_{p+n-q}$&\\

&&&&&------\\

&&&&$2\eps_1-\eps_q+\delta_1-\delta_{p+n-q-1}-\delta_{p+n-q}$&\\

&&&&$\eps_1-\delta_{p+n-q} \hbox{ }(m=n)$&\\

\hline

$\geq 1$&$\geq 1$&$\geq 3$&$0$&$m \not= p+q \pm 1$&\\

&&&&$2\eps_1-\eps_{m-p+q}+\delta_1-\delta_{p-1}-\delta_p$&\\

&&&&&------\\

&&&&$\eps_1+\eps_2-\eps_{m-p+q}+\delta_1-2\delta_p$&\\

&&&&$\eps_1-\delta_p \hbox{ }(m=n)$&\\

\hline

$\geq 1$&$\geq 1$&$\geq 1$&$\geq 1$&$m, n \not= p+q \pm 1$&\\

&&&&$2\eps_1-\eps_{m-p+q}+\delta_1-2\delta_{p+n-q}$&\\

&&&&&------\\

&&&&$\eps_1+\eps_2-\eps_{m-p+q}+\delta_1-\delta_{p+n-q-1}-\delta_{p+n-q}$&\\

&&&&$\eps_1-\delta_{p+n-q} \hbox{ }(m=n \geq 3)$&\\

\hline

\end{tabular}
}

\section*{Table 18. Spencer cohomology of $\fosp(m|2n)$}
%\small

\begin{tabular}{|c|c|c|c|}

\hline

$\bf r$&$\bf n$&$\bf H^{2,2}_{\fg_0}$&$\bf S^2(\fg_{-1})$\\

\hline

&&&$0\hbox{ (if }m = 2)$\\

$0$&$1$&------&\\

&&&$0,\hbox{ } \delta_1\hbox{ (if }m =3)$\\

\hline

$0$&$\geq 2$&$2\delta_1+2\delta_2$&$\delta_1+\delta_2, \hbox{ }0$\\

\hline

$1$&$1$&$\eps_1+\delta_1$&$2\eps_1,\hbox{ } 0$\\

\hline

$1$&$\geq 2$&$2\eps_1+\delta_1+\delta_2$&$2\eps_1,
\hbox{ }0$\\

\hline

$\geq 2$&$\geq 1$&$2\eps_1+2\eps_2$&$2\eps_1,\hbox{ } 0$\\

\hline

\end{tabular}

\end{document}